\documentclass[a4paper,12pt]{article}

\usepackage{mathtools}
\usepackage{amsmath,amssymb,amsthm, bm}
\usepackage{enumerate}
\usepackage[noadjust]{cite}

\usepackage[spanish, english]{babel}

\usepackage[ansinew]{inputenc}
\usepackage{here}
\usepackage{hyperref}

\usepackage{xcolor}

\usepackage{graphicx}

\usepackage[normalem]{ulem}

\newcounter{thm_compactness}
\newcounter{specialization_compactness}

\newcommand{\Cut}{\operatorname{Cut}}
\newcommand{\Aut}{\operatorname{Aut}}
\newcommand{\tp}{\operatorname{tp}}

\newcommand{\BB}{\mathcal{B}}
\newcommand{\CC}{\mathcal{C}}

\newcommand{\FF}{\mathcal{F}}
\newcommand{\GG}{\mathcal{G}}

\newcommand{\LL}{\mathcal{L}}
\newcommand{\MM}{\mathcal{M}}
\newcommand{\NN}{\mathcal{N}}

\newcommand{\SSS}{\mathcal{S}}

\newcommand{\XX}{\mathcal{X}}

\newcommand{\BBB}{\mathcal{B}}
\newcommand{\CCC}{\mathcal{C}}
\newcommand{\DD}{\mathcal{D}}
\newcommand{\U}{\mathcal{U}}

\newcommand{\Mon}{\mathcal{U}}

\newcommand{\ar}{\rangle}
\newcommand{\al}{\langle}
\newcommand{\exR}{M_{\pm\infty}}

\newcommand{\badS}{\mathbb{S}}

\newcommand{\R}{\mathbb{R}}

\newcommand{\lleq}{\preccurlyeq}
\newcommand{\lle}{\prec}

\newcommand{\rel}{ rel }

\newcommand{\vc}{\text{vc}}
\newcommand{\VC}{\text{VC}}

\newcommand{\prop}[2]{Prop(#1,#2)}

\newcommand{\mytau}{\tilde{\tau}}

\newcommand{\defn}{Def(M^n)}

\theoremstyle{definition}
\newtheorem{definition}{Definition}[section]
\newtheorem{remark}[definition]{Remark}
\newtheorem{example}[definition]{Example}
\newtheorem{fact}[definition]{Fact}

\theoremstyle{theorem}
\newtheorem{theorem}[definition]{Theorem}
\newtheorem{thmx}{Theorem}
 
\newtheorem{lemma}[definition]{Lemma}
\newtheorem{proposition}[definition]{Proposition}
\newtheorem{corollary}[definition]{Corollary}
\newtheorem{claim}[definition]{Claim}

\newtheorem*{question*}{Question}
\newtheorem{question}[definition]{Question}

\title{Types, transversals and definable compactness in o-minimal structures}
\author{Pablo And\'ujar Guerrero}
\date{}

\begin{document}
\maketitle 

\noindent
{\small \emph{Mathematics Subject Classification 2020.} 03C64 (Primary); 54A05, 54D30 (Secondary). }

\begin{abstract}
Through careful analysis of types inspired by~\cite{atw1} we characterize a notion of definable compactness for definable topologies in general o-minimal structures, generalizing results from~\cite{pet_pillay_07} about closed and bounded definable sets in o-minimal expansions of ordered groups. Along the way we prove a parameter version for o-minimal theories of the connection between dividing and definable types known in the more general dp-minimal context~\cite{simon_star_14}, through an elementary proof that avoids the use of existing forking and $\VC$ literature. In particular we show that, if an $A$-definable family of sets has the $(p,q)$-property, for some $p\geq q$ with $q$ large enough, then the family admits a partition into finitely many subfamilies, each of which extends to an $A$-definable type. 
\end{abstract}

\section*{Introduction}\label{section:intro}
We study the relationship between (complete) types, intersection properties, and the topological notion of definable compactness in general o-minimal structures. 

In~\cite{pet_pillay_07} Peterzil and Pillay showed that, in an o-minimal structure with definable choice functions, if a definable family of closed and bounded sets is consistent (i.e. has the finite intersection property), then it has a finite transversal (i.e. there exists a finite set that intersects each set in the family). They derived this from the work of Dolich~\cite{dolich04} on forking in o-minimal theories --in particular Dolich's proof that, in an o-minimal expansion of a field\footnote{More specifically, Dolich works in the slightly more general setting of a \textquotedblleft nice" o-minimal theory.}, a formula $\varphi(x,b)$ does not fork (or does not divide) over a set $A$ if and only if it is \textquotedblleft good over $A$". It is easy to see that being good over $A$ implies that $\varphi(x,b)$ extends to a (global) $A$-definable type. This equivalence between formulas not forking over a set $A$ and extending to $A$-definable types was conjectured for all dp-minimal theories by Simon~\cite[Conjecture 5.2]{simon15}, in the case where $A$ is a model. The most general result up to date in this direction is due to Simon and Starchenko~\cite{simon_star_14}, who proved that the conjecture holds for a large class of dp-minimal theories (see Theorem~\ref{thm:forking_definable_types}), including those that are linearly ordered (e.g. o-minimal) or have definable Skolem functions (e.g. the theory of the p-adics). The aforementioned result in~\cite{pet_pillay_07} highlights how this line of research on forking in tame theories can produce results about intersecting definable families of closed sets in some definable topology with some appropriate notion of topological compactness.  

Other results about intersecting definable families of sets in o-minimal structures include the work of Aschenbrenner and Fischer~\cite{aschen_fischer_11}, who proved a definable Helly theorem in definably complete expansions of real closed fields $\mathcal{R}$. They showed that, if a definable family of closed and bounded convex subsets of $R^n$ satisfies that any subfamily of size $n+1$ has nonempty intersection, then the whole family has nonempty intersection. 

Moreover, in the setting of NIP theories (a class that includes dp-minimal theories), the Alon-Kleitman-Matou\v{s}ek $(p,q)$-theorem~\cite{matousek04} implies the following. Any definable family of sets in an NIP structure with the $(p,q)$-property (the property that, for every $p$ sets in the family, some $q$ intersect), for some $p\geq q$ with $q$ large enough (i.e. $q$ greater than the $\VC$-codensity of the family), admits a partition into finitely many consistent subfamilies. The definable $(p,q)$-conjecture of Simon~\cite[Conjecture 2.15]{simon15} states that these subfamilies can be chosen to be definable (see~\cite{boxall_kestner_18} for recent work in this direction). 

In the present paper we study types and intersecting definable families of sets in connection with a notion of definable compactness that applies to definable topological spaces. A definable topological space in a structure $\MM$ is a topological space $(X,\tau)$, with $X\subseteq M^n$, such that there exists a basis for $\tau$ that is (uniformly) definable. The study of tame definable topological spaces is connected to the development of o-minimality and dates back to Pillay~\cite{pillay87}. A prominent example in o-minimal structures, beyond the canonical o-minimal euclidean topology, is given by definable groups~\cite{pillay88}, which admit a definable manifold structure in the sense of van den Dries~\cite[Chapter 10]{dries98} making them topological groups. Further examples within o-minimality are given by definable families of functions with definable norms~\cite{thomas12}, definable metric spaces~\cite{walsberg15}\cite{dries03}, and some classical non-metrizable spaces in topology~\cite[Appendix A]{andujar_thesis}. The first notion of definable compactness in o-minimal structures was introduced for definable manifold spaces, albeit the definition generalizes to any o-minimal definable topology, by Peterzil and Steinhorn~\cite{pet_stein_99}, and proved to be equivalent in the euclidean topology to being closed and bounded. This is the property that every definable curve converges, which we refer to in this paper as \emph{curve-compactness}. A more recent approach with wider applicability (including outside linearly ordered structures) to the notion, and the one we use in this paper, is as follows: every downward directed definable family of nonempty closed sets has nonempty intersection (a family is downward directed if for any two sets in the family there is a third one that is a subset of their intersection). This property has been studied in recent years in a general model theoretic context by Fornasiero~\cite{fornasiero}. In the context of definable topologies in o-minimal structures, it has been approached by Johnson~\cite{johnson14}, and by Thomas, Walsberg and the author~\cite{atw1}. Johnson~\cite{johnson14} noted that, within o-minimality, it is equivalent to curve-compactness for euclidean sets. In~\cite{atw1}, we proved this equivalence for all definable topologies in o-minimal expansions of ordered groups. Additionally we proved, in the case of an o-minimal expansion of an ordered field, that definable compactness is equivalent to \emph{type-compactness:} for every definable type that inhabits the topological space there is a point contained in every closed set in the type. From this we derived the equivalence with classical topological compactness in o-minimal expansions of the field of reals.  

We seek to characterize definable compactness in the general o-minimal setting, in particular in terms of intersecting definable families of closed sets, and also by showing the equivalence of the various compactness notions. Our study of o-minimal types allows us to set many of our questions and results outside a topological context, and later apply these to reach the desired characterization (similarly to how the work from Dolich~\cite{dolich04} on forking was used in~\cite{pet_pillay_07}). From this perspective we improve, within o-minimality, all the results mentioned above, in a number of ways. Our focus throughout is on definable families of sets, and our proofs are mostly elementary, in that they rely solely on o-minimal cell decomposition.

The structure and main results of the paper are as follows. 

In Section~\ref{section:prelim} we introduce the necessary preliminaries. 

In Section~\ref{section:types} we study the relationship between types and definable downward directed families of sets. Observe that a type basis is always downward directed. We say that a type has a uniform basis if it has a basis given by its restriction to some formula (this includes every $1$-type in an o-minimal structure). If the uniform basis is definable we simply call it a definable basis. Our main result in the section is the following.

\begin{thmx}[Theorem~\ref{thm:types} (definable version) and Remark~\ref{remark:parameter_thm_type}] \label{thm:intro_types}
Let $\SSS\subseteq \mathcal{P}(M^n)$ be a definable family of sets. The following are equivalent.
\begin{enumerate}[(1)]
\item \label{itm1: intro_thm_types}$\SSS$ extends to a definable $n$-type.
\item \label{itm2: intro_thm_types} There exists a definable downward directed family finer than $\SSS$. 
\item \label{itm3: intro_thm_types}$\SSS$ extends to an $n$-type with a definable basis of cells.  
\end{enumerate}
Moroever, the equivalence \eqref{itm1: intro_thm_types}$\Leftrightarrow$\eqref{itm2: intro_thm_types} still holds if we write ``$A$-definable" in place of ``definable", for any $A\subseteq M$.  
\end{thmx}

The implication \eqref{itm1: intro_thm_types}$\Rightarrow$\eqref{itm3: intro_thm_types} in Theorem~\ref{thm:intro_types} can be seen as a strong density result for types with a definable basis. The theorem expands on~\cite{atw1}, where explicit constructions for these definable bases are given in the cases where the underlying structure expands an o-minimal ordered group or field (see Theorem $27$ and Corollary $32$). This approach to studying types appears to be novel.

In Section~\ref{section:transversals} we investigate the $(p,q)$-property among definable families of sets. Using elementary o-minimal methods and Theorem~\ref{thm:intro_types} we prove the following. 

\begin{thmx}[Theorem~\ref{thm:finite_tame_transversals}] \label{thm:intro_FTT}
Let $\SSS$ be a definable family of nonempty sets in an o-minimal structure $\MM$. Let $q>\dim \cup\SSS$. The following are equivalent.   
\begin{enumerate}[(1)]
\item \label{itm:thm_transversals_1} There exists some $p\geq q$ such that $\SSS$ has the $(p,q)$-property.
\item \label{itm:thm_transversals_2} $\SSS$ can be partitioned into finitely many subfamilies, each of which extends to a definable type.    
\end{enumerate}  
\end{thmx}
We also use results from VC Theory to derive (Corollary~\ref{prop:vc_tame_transversals}) a version of Theorem~\ref{thm:intro_FTT} where the lower bound $\dim \cup\SSS$ is substituted by the VC-codensity of the family $\SSS$. 
 
In Subsection~\ref{subsection:types} we observe how Theorem~\ref{thm:intro_FTT} is equivalent to the aforementioned result about the equivalence between a formula non-forking over a model $M$ and extending to an $M$-definable type~\cite{simon_star_14}. We prove this last equivalence within o-minimality for forking over any set (not just a model); see Theorem~\ref{thm:FTT_forking}.

Finally, in Section~\ref{section:compactness} we use Theorems~\ref{thm:intro_types} and~\ref{thm:intro_FTT}, as well as Corollary~\ref{prop:vc_tame_transversals}, to characterize definable compactness as follows. 


\begin{thmx}[Theorem~\ref{thm:compactness}]\label{thm:intro_compactness} Let $\MM$ be an o-minimal structure. Let $(X,\tau)$ be a definable topological space in $\MM$. The following are equivalent. 
\begin{enumerate}[(1)]
\item \label{itm:intro_compactness_1} $(X,\tau)$ is definably compact, i.e. every downward directed definable family of nonempty closed sets has nonempty intersection. 
\item \label{itm:intro_compactness_2} Every definable type in $X$ has a limit, i.e. a point contained in every closed set in the type. 
\item \label{itm:intro_compactness_2.5} Any definable family of closed sets that extends to a definable type has nonempty intersection.
\item \label{itm:intro_compactness_3} Any consistent definable family of closed sets has a finite transversal, i.e. there exists a finite set that intersects every set in the family.
\item \label{itm:intro_compactness_5} Any definable family $\CC$ of closed sets with the $(p,q)$-property, where $p\geq q > \dim \cup\CC$, has a finite transversal.
\item \label{itm:intro_compactness_4} Any definable family $\CC$ of closed sets with the $(p,q)$-property, where $p\geq q$ and $q$ is greater than the $\VC$-codensity of $\CC$, has a finite transversal.
\setcounter{thm_compactness}{\value{enumi}}
\end{enumerate}
Moreover all the above imply and, if $\tau$ is Hausdorff or $\MM$ has definable choice, are equivalent to:
\begin{enumerate}[(1)]
\setcounter{enumi}{\value{thm_compactness}}
\item \label{itm:intro_compactness_6} $(X,\tau)$ is curve-compact, i.e. every definable curve in $X$ converges. 
\end{enumerate}
\end{thmx}

The precise role of Theorems~\ref{thm:intro_types} and \ref{thm:intro_FTT} in proving Theorem~\ref{thm:intro_compactness} is as follows. Theorem~\ref{thm:intro_types} gives the equivalence between definable compactness and type-compactness (i.e.~\eqref{itm:intro_compactness_1}$\Leftrightarrow$\eqref{itm:intro_compactness_2}), as well as the equivalence with condition~\eqref{itm:intro_compactness_2.5}. Parting from the equivalence \eqref{itm:intro_compactness_1}$\Leftrightarrow$\eqref{itm:intro_compactness_2.5} in Theorem~\ref{thm:intro_compactness}, it is easy to see that the characterizations of definable compactness given by \eqref{itm:intro_compactness_3}, \eqref{itm:intro_compactness_5} and \eqref{itm:intro_compactness_4} can be restated outside a topological setting, as the question of whether a ``sufficiently intersecting" definable family of sets admits a finite partition into subfamilies, each of which extends to a definable type. This motivated the work in Section~\ref{section:transversals}, where Theorem~\ref{thm:intro_FTT} yields the implications \eqref{itm:intro_compactness_2.5}$\Rightarrow$\eqref{itm:intro_compactness_3} and \eqref{itm:intro_compactness_2.5}$\Rightarrow$\eqref{itm:intro_compactness_5}, and moreover the VC version of Theorem~\ref{thm:intro_FTT} (Corollary~\ref{prop:vc_tame_transversals}) provides the implication \eqref{itm:intro_compactness_2.5}$\Rightarrow$\eqref{itm:intro_compactness_4}. Lastly, Theorem~\ref{thm:intro_types} is also used to prove the connection between definable compactness and curve-compactness (\eqref{itm:intro_compactness_1} and \eqref{itm:intro_compactness_6} respectively). 

We derive the following results from Theorem~\ref{thm:intro_compactness}. 
We use the equivalence \eqref{itm:intro_compactness_1}$\Leftrightarrow$\eqref{itm:intro_compactness_2}, together with the Marker-Steinhorn Theorem (Theorem~\ref{thm:marker_steinhorn}), to prove (Corollary~\ref{cor:definably_compact_iff_compact}) the equivalence between definable compactness and classical compactness in o-minimal expansions of $(\mathbb{R},<)$.
In~\cite[Theorem 2.1]{pet_pillay_07} it is shown that, in an o-minimal theory with definable choice functions, if a non-forking formula over a model $M$ defines (in a monster model) a closed and bounded set $C$, then $C\cap M^n \neq \emptyset$. And from this it is derived that every consistent definable family of closed and bounded sets has a finite transversal. We generalize this result (see Remark~\ref{remark:peterzil_pillay_aschenbrenner}) to any formula that defines a definably compact closed set in some $M$-definable topology, also dropping the assumption of having choice functions. In the same remark we explain how our work can be used to expand on the definable Helly's Theorem proved in~\cite{aschen_fischer_11} by showing in particular that, in an o-minimal expansion of an ordered field $\MM$, any definable family of convex subsets of $M^n$ with the property that every subfamily of size $n+1$ has nonempty intersection extends to a definable type. 



We end the paper with two appendices. In Appendix~\ref{section:appendix A} we use the approach of various proofs in this paper to present a shortened proof of the Marker-Steinhorn Theorem. Appendix~\ref{section:appendix B} is due to Will Johnson, and describes a counterexample to a parameter version of Lemma~\ref{lemma:existence_complete_directed_families}. 

\subsection*{Acknowledgements}

We thank Sergei Starchenko for the idea of a characterization of definable compactness in terms of $\VC$ theory and finite transversals, which resulted in Theorem~\ref{thm:compactness}, and has ultimately motivated much of the work in this paper. 

We thank Margaret Thomas for careful reading of early drafts and for providing extensive feedback. Similarly, we thank Matthias Aschenbrenner for his helpful comments on the paper, and in particular for suggesting the use of terminology ``$n$-consistent" and ``$n$-inconsistent". We also thank Antongiulio Fornasiero, who shared with us his unpublished notes on definable compactness~\cite{fornasiero}.     

Finally, a big thank you goes to Will Johnson for some very helpful conversations on the contents of Section~\ref{section:types}, for sharing with us an alternative proof of Theorem~\ref{thm:types}, and for agreeing to write Appendix~\ref{section:appendix B}, a counterexample to a parameter version of Lemma~\ref{lemma:existence_complete_directed_families}.

During the writing of this paper the author was supported by the Institute Henri Poincare (during the 2018 Program on Model Theory, Combinatorics and Valued Fields), the Canada Natural Sciences and Engineering Research
Council (NSERC) Discovery Grant RGPIN-06555-2018, the Graduate School and Department of Mathematics at Purdue University, and by the Fields Institute for Research in Mathematical Sciences (during the 2021 Thematic Program on Trends in Pure and Applied Model Theory).

\section{Preliminaries}\label{section:prelim}

We fix a first order o-minimal $\LL$-structure $\MM=(M,<,\ldots)$ expanding a dense linear order without endpoints. Throughout unless otherwise specified \textquotedblleft definable" means \textquotedblleft definable in $\MM$ with parameters". Let $n$, $m$ and $k$ denote natural numbers. Throughout $\defn$ denotes the algebra of definable subsets of $M^n$. 

All variables and parameters $u$, $x$, $a$ ... are $n$-tuples for some $n$. Let $l(u)$ denote the length of any given variable or parameter $u$. 

Every formula we consider is in $\LL$. Let $\varphi(u,v)$ and $\phi(v)$ be formulas in $m+n$ and $n$ free variables respectively. For any set $S\subseteq M^n$ and $a\in M^m$, let $\varphi(a, S)=\{b\in S, \MM\models \varphi(a,b)\}$ and $\phi(S)=\{b\in S : \MM\models \phi(b)\}$. 

A family of sets $\SSS\subseteq \defn$ is uniform\footnote{We choose this terminology in part to justify the description of a \textquotedblleft uniform topology" among types in Remark~\ref{remark:type_topology}, in connection with the uniform convergence topology among function spaces.} if there is a formula $\varphi(u,v)$  and a set $\Omega\subseteq M^{l(u)}$ such that $\SSS=\{\varphi(u,M^{l(v)}) : u\in \Omega\}$. We call $\Omega$ the index set. The family $\SSS$ is definable whenever $\Omega$ can be chosen definable. For example if $\MM$ denotes the real algebraic numbers, then the family of intervals $(s,t)$ for $s<\pi<t$ is uniform but not definable.     
 
For a given $n$, let $\pi$ denote the projection $M^{n+1}\rightarrow M^n$ onto the first $n$ coordinates, where $n$ will often be omitted and clear from context. For a family $\SSS\subseteq Def(M^{n+1})$ let $\pi(\SSS)=\{\pi(S) :S\in\SSS\}$. 

Given a function $f$ we denote its domain $dom(f)$ and its graph by $graph(f)$. We adopt the convention of saying that a function $f:M^n\rightarrow \exR$ is definable if its restriction $f|_{f^{-1}(M)}$ and the sets $f^{-1}(\{+\infty\})$ and $f^{-1}(\{-\infty\})$ are definable, and similarly say that a family $\FF=\{f_u : u\in \Omega\}$ of said functions is definable if the families $\{ f_u|_{f_u^{-1}(M)} : u\in \Omega\}$, $\{ f^{-1}_u(\{+\infty\}) : u\in \Omega\}$ and  $\{ f^{-1}_u(\{-\infty\}) : u\in \Omega\}$ are. 

For a definable set $X$ we denote its frontier and boundary (in the canonical \textquotedblleft euclidean" o-minimal topology) by $\partial C$ and $bd(X)$ respectively.

\subsection{Cell decomposition}

We call $\DD$ a cell decomposition of $S\subseteq M^n$ if $\DD$ is a finite partition of $S$ into cells and, for any $0<m<n$, the family $\{\pi_{\leq m}(D): D\in\DD\}$ is a cell partition of $\pi_{\leq m}(S)$, where $\pi_{\leq m}:M^n\rightarrow M^m$ denotes the projection to the first $m$ coordinates.  By o-minimality every definable set admits a cell decomposition. 

\textbf{Uniform cell decomposition:} If $S\subseteq M^{n+m}$ is a definable set and $\DD$ is a cell decomposition of $S$, then, for any fiber $S_u=\{v: \al u, v\ar \in S\}$, $u\in M^n$, the family of fibers $\{ D_u : D\in\DD, D_u\neq \emptyset\}$ is a cell decomposition of $S_u$.

\subsection{Type preliminaries}\label{subsection:types}

Let $p$ be an $n$-type (over a set $A$). We refer to $p$ indistinctly as a consistent family of formulas (with parameters from $A$) with $n$-free variables $\{ \varphi(v) : \varphi\in p\}$ and as a consistent (i.e. having the finite intersection property) family of definable (over $A$) sets $\{\varphi(M^n) : \varphi\in p\}$, with a preference for the set definition. For $n$ and a structure $\NN=(N,\ldots)$, let $S_n(N)$ denote the set of complete $n$-types over $N$. Unless otherwise specified all types we consider are global (over $M$) and complete. By this identification an $n$-type is an ultrafilter on $\defn$.

A type $p$ with object variable $v$ is ($A$-)definable if, for every partitioned formula $\varphi(u,v)$, the family of sets of the form $\varphi(u,M^{l(v)})$ in $p$, is ($A$-)definable. In other words if, for every $\varphi(u,v)$, the set $\{u : \varphi(u, M^{l(v)})\in p\}$ is ($A$-)definable. 

Note that, if $p$ is an ($A$-)definable type, then the projection $\pi(p)$ is also an ($A$-)definable type.

\subsection{Preorders}\label{remark:preorder}

Recall that a preorder is a reflexive transitive relation. A preordered set $(X,\lleq)$ is a set $X$ together with a preorder $\lleq$ on it. It is definable if the preorder is definable. A preorder $\lleq$ on a set $X$ induces an equivalence relation $\sim$ on $X$ given by $x\sim y$ if and only if $x\lleq y$ and $y\lleq x$. We use notation $x\prec y$ to mean $x\lleq y$ and $x \not\sim y$. Given a set $Y\subseteq X$ let $x\lleq Y$ (respectively $x\prec Y$) mean $x\lleq y$ (respectively $x\prec y$) for every $y\in Y$. For every $x,y\in X$ let $[x,y]_{\lleq}=\{z\in X : x\lleq z\lleq y\}$ and $(x,y)_{\lleq}=\{z\in X : x\prec z\prec y\}$.

Let $p$ be an $n$-type and $\GG$ be the collection of all definable functions $f:M^n\rightarrow \exR$ such that $dom(f)\in p$. Then $p$ induces a total preorder $\lleq$ on $\GG$ given by $f \lleq g$ if and only if $\{u\in M^n : f(u)\leq g(u)\}\in p$. In other words, $f\lleq g$ when $f(\xi)\leq g(\xi)$ for some (every) realization $\xi$ of $p$.


Let $\FF=\{f_x : x\in X\}\subseteq \GG$ be a uniform family. Without loss of clarity we will often abuse notation and refer to $\lleq$ too as the total preorder on the index set $X$ given by $x\lleq y$ if and only if $f_x \lleq f_y$. Note that, if $\FF$ and $p$ are $A$-definable, then $\lleq$ on $X$ is $A$-definable too. 



\subsection{Transversals and intersection properties}

For $\SSS$ a family of sets $\SSS$ and $F$ a set let $\SSS\cap F=\{ S\cap F : S\in\SSS\}$. The set $F$ is a transversal of $\SSS$ if it intersects every set in $\SSS$, i.e. if $\SSS\cap F$ does not contain the empty set. In this paper we are interested in the property that a definable family of sets has a finite transversal, and on a similar property that we introduce in Section~\ref{section:transversals} (Definition~\ref{dfn:ftt}).

We say that a family of sets $\SSS$ is $n$-consistent if any subfamily of cardinality at most $n$ has nonempty intersection. A family is consistent if it is $n$-consistent for any $n$. We say that $\SSS$ is $n$-inconsistent if every subfamily of cardinality $n$ has empty intersection. We say that $\SSS$ has the $(p,q)$-property, for $p\geq q>0$, if every subfamily $\SSS'\subseteq \SSS$ of size $p$ contains a subfamily of size $q$ with nonempty intersection. Note that, for $p$ large with respect to $q$, this is a rather weak intersection property.

Note that $\SSS$ does not have the $(p,q)$-property if and only if there exists a subfamily of $\SSS$ of size $p$ that is $q$-inconsistent. 

As an example of a fact involving transversals we include the following proposition, whose proof, which can be derived from Mirsky's theorem \cite[Theorem 2]{mirsky71}, we leave to the reader.
\begin{proposition}\label{prop:Mirsky}
Let $\SSS\subseteq \mathcal{P}(M)$ be a family of closed intervals. Let $k\geq 1$ be the maximum such that there exists $k$ pairwise disjoint sets in $\SSS$. Then $\SSS$ has a transversal of size $k$.
\end{proposition}

\subsection{$\VC$ theory}

In this subsection we present the basic notions of $\VC$ theory, and state the Alon-Kleitman-Matou\v{s}ek $(p,q)$-theorem, as well as a known result bounding the $\VC$-condensity of formulas in o-minimal structures. These results have applications in Section~\ref{section:transversals}. The following brief introduction serves as context for them. Only the notion of $\VC$-codensity will appear in the proofs in this paper, always in reference to the theorems.  

A pair $(X,\SSS)$, where $X$ is a set and $\SSS$ is a family of subsets of $X$, is called a \emph{set system}. For a set $F\subseteq X$ we say that $\SSS$ \emph{shatters} $F$ if $\SSS\cap F=\mathcal{P}(F)$. The $\VC$-dimension of $\SSS$, denoted by $\VC(\SSS)$, is the maximum cardinality of a finite set $F$ such that $\SSS$ shatters $F$, if one such set exists. Otherwise we write $\VC(\SSS)=\infty$. The shatter function $\pi_\SSS:\omega\rightarrow \omega$ is defined by $\pi_\SSS(n)=\max\{ |\SSS\cap F| : F\subseteq X,\, |F|=n\}$. So the $\VC$-dimension of $\SSS$ is the largest $n$ such that $\pi_{\SSS}(n)=2^n$. 

Suppose that $\SSS$ has $\VC$-dimension $k<\infty$. The Sauer-Shelah Lemma states then that $\pi_{\SSS}= O(n^k)$ (that is, $\pi_{\SSS}(n)/n^k$ is bounded at infinity). The $\VC$-density of $\SSS$, denoted by $\vc(\SSS)$, is the infimum over all real numbers $r\geq 0$ such that $\pi_{\SSS}= O(n^r)$. Hence $\vc(\SSS)\leq \VC(\SSS)$. 

The dual set system of $\SSS$ is the set system $(X^*,\SSS^*)$, where $X^*=\SSS$ and $\SSS^*=\{\SSS_x : x\in X\}$ where $\SSS_x=\{S\in \SSS : x\in S\}$. Hence $\SSS^*$ shatters $\SSS'\subseteq X^*$ if and only if every field in the Venn diagram induced by $\SSS'$ on $X$ contains at least one point. The dual shatter function of $\SSS$ is $\pi^*_\SSS=\pi_{\SSS^*}$. Similarly there is the dual $\VC$-dimension of $\SSS$, $\VC^*(\SSS)=\VC(\SSS^*)$, and $\VC$-codensity, $\vc^*(\SSS)=\vc(\SSS^*)$. These satisfy that $\vc^*(\SSS)\leq \VC^*(\SSS)\leq 2^{\VC(\SSS)+1}$.

A theory is NIP if every uniform family of sets in any model has finite $\VC$-dimension. Every o-minimal theory is NIP.  

For convenience we state Matou\v{s}ek's theorem in terms of $\VC$-codensity. For a finer statement see~\cite{matousek04}.

\begin{theorem}[Alon-Kleitman-Matou\v{s}ek $(p,q)$-theorem]\label{thm:vc_tranversal}
Let $p \geq q >0$ be natural numbers and let $(X,\SSS)$ be a set system such that $\vc^*(\SSS)<q$. Then there is $n$ such that, for every finite subfamily $\FF\subseteq \SSS$, if $\FF$ has the $(p,q)$-property then it has a transversal of size at most $n$. 
\end{theorem}


The following Corollary will be useful in section~\ref{section:transversals}. It is a reformulation of the main result for weakly o-minimal structures (a class that contains o-minimal structures) in~\cite{vc_density} by Aschenbrenner, Dolich, Haskell, Macpherson and Starchenko. It was previously proved for o-minimal structures by Wilkie (unpublished) and Johnson-Laskowski~\cite{john_las_10}, and for o-minimal expansions of the field of reals by Karpinski-Macintyre~\cite{kar_mac_00}.

\begin{theorem}[\cite{vc_density}, Theorem 6.1]\label{cor:vc_density}
Let $\SSS\subseteq \defn$ be a uniform family of sets. Then \mbox{$\vc^*(\SSS)\leq n$}. 
\end{theorem}

Since throughout this paper $p$ and $q$ are employed as terminology for types, in the subsequent sections we address the $(p,q)$-property in terms of $m$ and $n$, e.g. the $(m,n)$-property.

\section{Types}\label{section:types}

This section contains results about types in o-minimal structures that will be used later on. Our main result is Theorem~\ref{thm:types}, which shows that \emph{types with a uniform basis} are dense in a rather strong sense among all types. The content of this section was motivated by results on types and downward directed families of sets in o-minimal expansions of ordered groups and fields in~\cite{atw1}.

\begin{definition}
A family of sets $\SSS$ is \emph{downward directed} if, for every pair $S, S'\in\SSS$, there exists $S''$ such that $S''\subseteq S\cap S'$. For convenience we also ask that $\SSS$ does not contain the empty set. Equivalently $\SSS$ is a \emph{filter basis}.   
\end{definition}

We are interested in the interface between types and uniform downward directed families of sets. 

\begin{definition}
Let $p$ be an $n$-type. A \emph{type basis} for $p$ is a partial type $q\subseteq p$ that is downward directed and generates $p$, i.e. $p=\{X\in\defn : Y\subseteq X \text{ for some } Y\in q\}$. 


A \emph{uniform type basis} for $p$ is a type basis given by a uniform family of sets, e.g. $q=\{\varphi(u, M^n) : u\in \Omega\}$. 
\end{definition}

In general, we only consider type bases that are uniform. In particular we refer to bases as ``definable" in the sense of ``uniform definability" (not just ``type definability"). Note that a type with a uniform basis is $A$-definable if and only if it has an $A$-definable basis.  

In o-minimal and weakly o-minimal structures every $1$-type $p$ has a uniform basis. When $p$ is not realized this basis is given by all the open intervals in $p$. A natural question to ask is whether all types have a uniform basis. We observe that this already fails among definable types in any o-minimal expansion of $(\mathbb{R},+,<)$.

The dimension of a type is the lowest dimension of a set in it. In~\cite[Corollary 32]{atw1} we proved the following two statements. If $\MM$ expands an ordered group, then every definable type with a uniform basis is of dimension at most $2$. And if $\MM$ expands an ordered field then one such type is of dimension at most $1$. The Marker-Steinhorn theorem (Theorem~\ref{thm:marker_steinhorn}) implies that every type in an o-minimal expansion of $(\R,<)$ is definable. So, in any o-minimal expansion of $(\R,<,+)$ or $(\R,<,+,\cdot)$, any type of sufficiently large dimension lacks a uniform basis. Finally recall that, for any $n$, any expansion of the partial type $\{M^n\setminus X : X\in \defn, \, \dim X<n\}$ is $n$-dimensional.  

Nevertheless, we end this section by arguing, using the next two lemmas, that types with
a uniform basis are dense among all types in the context of a topology strictly stronger than the usual Stone type topology, and an analogous density result
among definable types.



\begin{definition}
Let $\SSS,\,\FF \subseteq \defn$ be consistent families of sets. We say that $\SSS$ is \emph{finer than $\FF$} if, for every $F\in\FF$, there is $S\in\SSS$ such that $S\subseteq F$. That is, the filter generated by $\SSS$ is finer than the one generated by $\FF$. 

Let $\SSS\subseteq \defn$ be downward directed and $\XX\subseteq \defn$. We say that $\SSS$ is \emph{complete for $\XX$} if, for every $X\in\XX$, there is $S\in\SSS$ such that either $S\subseteq X$ or $S \cap X=\emptyset$. We say that $\SSS$ is \emph{complete} if it is complete for $\defn$. In other words, if $\SSS$ is a type basis for some type.  

\end{definition}

We may now proceed with the results in this section. Theorem~\ref{thm:types} will follow from Lemmas~\ref{lemma:existence_complete_directed_families} and~\ref{lemma:existence_complete_directed_families_2.1}. The following facts are easy to prove and will be used later on, often without notice.

\begin{fact}\label{fact:downward_directed_family_1}
Let $\SSS$ be a downward directed family and $X$ be a set such that, for every $S\in \SSS$, $S\cap X\neq \emptyset$. Then the family $\SSS\cap X=\{S\cap X : S\in\SSS\}$ is downward directed. 
\end{fact}

\begin{fact}\label{fact:downward_directed_family_2}
Let $\SSS$ be a downward directed family and $\XX$ be a finite covering of a set $X$. If $S\cap X\neq \emptyset$ for every $S\in\SSS$ then there exists some $Y\in\XX$ such that $S\cap Y\neq \emptyset$ for every $S\in\SSS$. 
\end{fact}

\begin{fact}\label{fact:downward_directed_family_3}
Let $\SSS$ be a definable downward directed family that is complete for $\XX$ a finite definable partition of a set $X$. If $S\cap X\neq \emptyset$ for every $S\in\SSS$ then there exists a unique $Y\in\XX$ such that $S\subseteq Y$ for some $S\in\SSS$.
\end{fact}


\begin{lemma}\label{lemma:existence_complete_directed_families}
Let $\SSS$ be a uniform downward directed family. There exist a uniform downward directed family of cells $\CCC$ that is complete and finer than $\SSS$. If $\SSS$ is a definable family then $\CCC$ can be chosen definable. 
\end{lemma}

The work of Thomas, Walsberg and the author in~\cite{atw1} involves showing that, if $\MM$ expands an ordered field, every downward directed definable family $\SSS$ admits a finer (complete by o-minimality) family of the form $\{\gamma[(0,t)] : 0<t\}$ for some definable curve $\gamma:(0,\infty)\rightarrow \cup\SSS$. Moreover, whenever $\MM$ expands an ordered group, there exists one such finer family of the form $\{\gamma[(0,s)\times (t,\infty)] : 0<s<t\}$ for some definable map $\gamma:(0,\infty)\times (0,\infty)\rightarrow \cup\SSS$. As mentioned before, this provides us with information on the maximum dimension of sets in a complete definable downward directed family. In general however the picture is different. In a dense linear order $(M,<)$, for any $n$ the definable nested family of $n$-dimensional sets of the form $\{ \al x_1,\ldots, x_n\ar : a<x_1<\cdots<x_n\}$, for $a\in M$, is complete. 


\begin{proof}[Proof of Lemma~\ref{lemma:existence_complete_directed_families}]
Fix $\SSS\subseteq \mathcal{P}(M^n)$ a uniform downward directed family. We may assume that $\cap \SSS=\emptyset$, since otherwise it suffices to take $\CCC=\{\{x\}\}$ for any $x\in \cap \CCC$. We prove the lemma by induction on $n$. We assume that $\SSS$ is definable. When $\SSS$ is not definable the base case follows from the fact that by o-minimality every $1$-type has a uniform basis of intervals, and the inductive step follows the same proof we present. 

\textbf{Base case: $n=1$.} 

Let \mbox{$H=\{ t\in M : \exists S\in\SSS, S\cap (-\infty,t]= \emptyset \}$}. If $H$ is empty then, by o-minimality, for every $S\in\SSS$ there is $t_S$ such that $(-\infty, t_S)\subseteq S$, in which case we may take $\CCC=\{(-\infty, t) : t\in M\}$. 

Now suppose that $H$ is nonempty. Note that $H$ is an interval in $M$ (possibly right closed) unbounded from below. Let $a=\sup H$.

If $a=\max H<\infty$ then there exists $S_a\in \BB$ such that $S_a\cap (-\infty,t]= \emptyset$. We note that, for any $S\in\SSS$, there exists $t_S>a$ such that $(a,t_S)\subseteq S$, so we may take $\CCC=\{(a,t) : t>a\}$. Otherwise by o-minimality there exists $S\in \SSS$ and $r_S>a$ such that $(a,r_S)\cap S=\emptyset$, but then any set $S'\in \SSS$ with $S'\subseteq S_a \cap S$ satisfies that $(-\infty, r_S)\cap B' =\emptyset$, contradicting that $a=\sup H$.

If $a\notin H$ then, for every $S\in \SSS$, $S \cap (-\infty, a]\neq \emptyset$. We show that it suffices to take $\CCC:=\{(t,a) : t<a\}$. Let $S_a\in\SSS$ be such that $a\notin S_a$. Suppose towards a contradiction that there exists $S\in\SSS$ and $r_S<a$ such that $(r_S,a)\cap S=\emptyset$. Since $r_S\in H$, there is $S'\in \SSS$ such that $S' \cap (-\infty,r_S]=\emptyset$. But then, any $S''\in\BB$ with $S''\subseteq S'\cap S \cap S_a$ satisfies that $(-\infty, a]\cap S'' =\emptyset$, contradiction. This completes the proof of the base case.
 
\textbf{Inductive step: $n>1$.} 

Suppose that there exists a cell $C\subseteq M^n$ with $\dim C=m<n$ such that $S\cap C\neq \emptyset$ for every $S\in \SSS$. Let $\pi_C: C\rightarrow M^m$ be a projection that is a homeomorphism onto an open cell. By Fact~\ref{fact:downward_directed_family_1} the family $\{\pi_C(S\cap C) : S\in\SSS\}$ is a definable downward directed family of subsets of $M^m$. By inductive hypothesis it admits a finer complete definable downward directed family of cells $\CC'$. By passing to a subfamily of $\CC'$ if necessary we may assume that every cell in $\CC'$ is a subset of $\pi_C(C)$. The definable family of cells $\CC=\{ \pi_C^{-1}(C') : C'\in\CC'\}$ is then downward directed, complete, and finer than $\SSS$. 

Hence onwards we assume that, for every cell $C\subseteq M^n$ with $\dim C <n$, there exists some $S\in\SSS$ such that $S\cap C=\emptyset$. By Fact~\ref{fact:downward_directed_family_2} and o-minimal cell decomposition it follows that, for every definable set $X\subseteq M^n$,
\begin{equation}\label{eqn:proof_complete_directed_families_1}\tag{$\dagger$}
\text{ if } \dim X < n \text{ then there exists } S\in \SSS \text{ such that } S\cap X=\emptyset.
\end{equation}
In particular $\dim S=n$ for any $S\in\SSS$.

Consider the definable downward directed family $\pi(\SSS)=\{\pi(S): S\in \SSS\}$. By induction hypothesis there exists a definable downward directed family of cells $\BB$ in $M^{n-1}$ that is complete and finer than $\pi(\SSS)$. Let $\FF=\{ S \cap  (B\times M) : S\in\SSS, B\in\BB\}$. We show that this definable family is downward directed. 

Let $S, S'\in \SSS$ and $B, B'\in\BBB$. Let $S''\in \SSS$ be such that $S''\subseteq S\cap S'$. Since $\BBB$ is finer than $\pi(\SSS)$ and downward directed there exists $B''\in\BBB$ with $B''\subseteq \pi(S'')\cap B\cap B'$. Then $\emptyset\neq S''\cap (B''\times M) \subseteq S \cap (B\times M) \cap S' \cap (B'\times M)$.
 
Clearly $\FF$ is finer than $\SSS$. We proceed by proving the lemma for $\FF$ in place of $\SSS$. 

Let $\FF=\{\varphi(u,M^n) : u\in \Omega\subseteq M^m\}$. Let $\DD$ be a cell decomposition of $\varphi(M^m,M^n)$. Then, by uniform cell decomposition, for every $u\in \Omega$ the family of fibers $\{D_u : D\in\DD,\, D_u\neq \emptyset\}$ is a cell decomposition of $\varphi(u,M^n)$ and the family $\{\pi(D_u) : D\in\DD,\, D_u\neq \emptyset\}$ is a cell decomposition of $\pi(\varphi(u,M^n))$. For each $u\in\Omega$ let $\DD_u=\{D_u : D\in\DD\}$. 

For the rest of the proof we assume that, for every $F\in \FF$, there is a unique $u\in \Omega$ with $F=\varphi(u,M^n)$. Then for clarity we identify $F$ with $u$ by writing $D_F$ in place of $D_u$ and $\DD_F$ in place of $\DD_u$, onwards omitting the subscript $u$ entirely. This assumption is valid if $\MM$ has elimination of imaginaries. Nevertheless it is adopted entirely for clarity and the proof can be written in terms of $u\in \Omega$ instead of $F=\varphi(u,M^n)$. 

We now wish to show that there exists a definable family of open cells $\CC=\{ C_F\in \DD_F : F\in\FF\}$ such that, for every $F\in \FF$, the following two conditions hold. 
\begin{enumerate}[(i)] 
\item \label{itm1:lemma_types_cells_proof} $C_F\cap F'\neq\emptyset$ for every $F'\in\FF$. 
\item \label{itm2:lemma_types_cells_proof} For every $C \in \DD_F$, if $C\cap F'\neq \emptyset$ for every $F'\in \FF$, then $C=C_F$ or $\pi(C)=\pi(C_F)$ and $C_F < C$, i.e. if $C_F=(f,g)$ and $C=(f',g')$ then $g\leq f'$.  
\end{enumerate} 

Note that, for every $F\in \FF$, there can be at most one $D\in \DD$ such that $D_F=C_F$ (i.e. satisfying~\eqref{itm1:lemma_types_cells_proof} and~\eqref{itm2:lemma_types_cells_proof} above).

Observe that, for every $D\in\DD$, the family of $F\in\FF$ such that $D_F=C_F$ is definable. Consider these induced subfamilies of $\FF$, as $D$ ranges in $\DD$. If $\CC$ as described above exists, then these subfamilies cover $\FF$. Additionally, by downward directedness, one of them must be finer than $\FF$. Hence, if $\CC$ exists, it is definable, and in particular, after passing to a subfamily of $\FF$ if necessary, we may in fact assume that $\CC$ is a family of fibers over a single cell in $\DD$. We show that $\CC$ exists.

For any $F\in\FF$ let $\DD'_F:=\{C\in \DD_F : C\cap F'\neq \emptyset\, \forall F'\in \FF\}$. By Fact~\ref{fact:downward_directed_family_2} this family is nonempty. By~\eqref{eqn:proof_complete_directed_families_1} every cell in $\DD'_F$ must be open. Now note that, since $\BB$ is complete and, by construction, $\pi(\FF)$ is finer that $\BB$, $\pi(\FF)$ is complete. By Fact~\ref{fact:downward_directed_family_3} it follows that, for every $F\in\FF$, there exists a unique set in $\pi(\DD_F)$, say $B_F$, such that $\pi(F')\subseteq B_F$ for some $F'\in\FF$.  Consequently, for any $C, C'\in \DD'_F$, it holds that $\pi(C)=\pi(C')=B_F$. Let $C_F=(f,g)$ be the cell in $\DD'_F$ satisfying that, for any other cell $(f',g')\in\DD'_F$, $g \leq f'$.

We have defined $\CC=\{ C_F : F\in\FF\}$. This family is clearly finer than $\FF$. We prove that it is downward directed and complete.  

\begin{claim}\label{claim:existence_complete_direted_families_1}
For every pair $F, F' \in \FF$, $C_F\cap C_{F'} \neq \emptyset$. 
\end{claim}  
Let $C_F=(f,g)$ and $C_{F'}=(f',g')$, and suppose towards a contradiction that $C_F\cap C_{F'}=\emptyset$. Recall that $\pi(\FF)$ is downward directed and complete, and notation $B_F$ for the projection of any (every) cell in $\DD'_F$. By~\eqref{itm1:lemma_types_cells_proof}, there exists some $B\in \pi(\FF)$ such that $B\subseteq B_F\cap B_{F'}$. Without loss of generality suppose that $g'|_B\leq f|_B$.

Let $\DD_F\setminus \DD'_F= \{ C(1), \ldots, C(l)\}$. By definition of $\DD'_F$ for every $1\leq i\leq l$ there exists some $F(i)\in \FF$ such that $C(i)\cap F(i)=\emptyset$. Let $F''\in \FF$ be such that $F''\subseteq F \cap (\cap_{i=1}^l F(i)) \cap (B \times M)$. Then, by~\eqref{itm2:lemma_types_cells_proof}, we have that $F''\subseteq (f|_B,+\infty)$. But then $F''\cap C_{F'}=\emptyset$, contradicting~\eqref{itm1:lemma_types_cells_proof}.


\begin{claim}
The family $\CC_{\FF}$ is downward directed. 
\end{claim}
Let $C, C' \in \CC_{\FF}$. Recall that $\FF$ is finer than $\SSS$. By~\eqref{eqn:proof_complete_directed_families_1}, let $F\in \FF$ be such that $F\cap \partial C=\emptyset$ and $F\cap \partial C'=\emptyset$. By Claim~\ref{claim:existence_complete_direted_families_1} $C_F \cap C \neq \emptyset$. Since every cell is definably connected and every cell in $\CC_{\FF}$ is open, it follows that $C_F\subseteq C$. By the same argument $C_F\subseteq C'$. Hence $C_F \subseteq C \cap C'$. 

Finally we show that $\CC$ is complete. Suppose otherwise, in which case there exists a definable set $X\subseteq M^n$ such that, for every $C \in \CC$, $C\cap X \neq \emptyset$ and $C\setminus X \neq \emptyset$. By~\eqref{eqn:proof_complete_directed_families_1}, $\dim X = n$, and we may find $F\in \FF$ such that $F \cap bd(X)=\emptyset$, having in particular that $C_F\cap bd(X)=\emptyset$. It follows that $C_F$ is not definably connected, contradicting that $C_F$ is a cell. 
\end{proof}

Notice that, in the base case of the proof of Lemma~\ref{lemma:existence_complete_directed_families}, if $\SSS$ is definable then the \textquotedblleft$a$" is definable over the same parameters as $\SSS$. It follows that $\CC$ is definable over the same parameters as $\SSS$. We might ask if this holds in higher dimensions. If $\MM$ has definable choice functions then it does, by the usual argument involving passing to an elementary substructure. In Appendix~\ref{section:appendix B} Will Johnson proves that, in general within o-minimality, the answer to that question is no.

We now present a second lemma connecting types and uniform downward directed families of sets. 

\begin{lemma}\label{lemma:existence_complete_directed_families_2.1}
Let $p$ be an $n$-type and $\SSS$ be a uniform family of subsets of $M^n$. There exists a uniform downward directed family of sets $\FF\subseteq p$ that is complete for $\SSS$. If $p$ is definable then $\FF$ can be chosen definable over the same parameters as $p$.  
\end{lemma}
\begin{proof}
Since otherwise the result is immediate we may assume that $p$ is not realised in $\MM$. We proceed by induction on $n$. For simplicity we prove the case where $p$ is definable. The case where $p$ is not definable follows by the same arguments. The fact that $\FF$ can be chosen definable over the same parameters as $p$ follows by keeping track of parameters. 

We start by reducing the proof to the case where $\SSS$ is a definable family of cells and $\SSS\subseteq p$.

By expanding $\SSS$ if necessary let us assume that it is a family of the form $\SSS=\{\varphi(u,M^n) : u\in M^m\}$. Let $\DD$ be a cell decomposition of $M^{m}\times M^n$ compatible with $\varphi(M^m, M^n)$. Then, for every $u\in M^m$, the family of fibers $\{D_u: D\in \DD\}$ is a partition of $M^n$ compatible with $\varphi(u,M^n)$. This means that there exists some $D\in \DD$ such that $D_u\in p$, and moreover either $D_u\subseteq \varphi(u,M^n)$ or $D_u\cap \varphi(u,M^n)=\emptyset$.

Consider the union $\CC$ of the families $\{D_u : u\in M^m, D_u\in p\}$ for $D\in\DD$. Then $\CC$ is a definable family of cells with $\CC\subseteq p$ and, for every $S\in\SSS$, there is some $C\in \CC$ such that either $C\subseteq S$ or $C\cap S=\emptyset$. Clearly it suffices to prove the lemma for $\CC$ in place of $\SSS$.

Hence onwards we assume that $\SSS$ is a definable family of cells and $\SSS\subseteq p$. This means in particular that we must prove the existence of a downward directed definable family $\FF\subseteq p$ that is finer than $\SSS$. 



\textbf{Base case: $n=1$.} 

In this case the result follows from the fact that every non-realized $1$-type has a uniform basis of intervals. In particular when the type is definable this basis is of the form $\{(a,t) , t>a\}$ for some $a\in M\cup \{-\infty\}$, or $\{(t,a) : t<a\}$ for some $a\in M\cup\{+\infty\}$. These are complete definable families of nested sets.  


\textbf{Inductive step: $n>1$.} 

Suppose there exists $S^*\in \SSS$ that is the graph of some definable continuous function $\pi(S^*)\rightarrow M$. Note that the type $\pi(p)=\{ \pi(X) : X\in p\}$ is definable. Note that the definable family $\SSS\cap \SSS=\{S\cap S' : S,S'\in \SSS\}$ is a subset of $p$. We apply the induction hypothesis to the type $\pi(p)$ and family $\pi(\SSS\cap \SSS)$ and reach a downward directed definable family $\GG\subseteq \pi(p)$ that is finer than $\pi(\SSS\cap \SSS)$. Let $\FF=\{ S\cap (G\times M) : S\in \SSS, G\in \GG\}$. Then $\FF$ is clearly definable and finer than $\SSS$. We show that it is a subset of $p$ and downward directed. Firstly, since $\GG\subseteq \pi(p)$, for every $G\in \GG$ it holds that $G \times M\in p$, and thus $S\cap (G\times M)\in p$ for every $S\in \SSS$. Secondly, since $\GG$ is downward directed and finer that $\pi(\SSS\cap \SSS)$, for every $S, S'\in \SSS$ there exists $G\in \GG$ such that $G\subseteq \pi(S^*\cap S) \cap \pi(S^*\cap S')$. Hence in particular $S^*\cap (G\times M) \subseteq S\cap S'$. It follows that, for any $G', G''\in \GG$, if $G'''\subseteq G\cap G'\cap G''$, then $S^*\cap (G'''\times M) \subseteq S\cap (G'\times M) \cap S' \cap (G''\times M)$.

We now prove the case where every cell $S\in\SSS$ is of the form $(f_S,g_S)$ for definable continuous functions $f_S,g_S:\pi(S)\rightarrow \exR$. Consider the definable families of sets $\SSS_0=\{ (-\infty, g_S) : S\in \SSS\}$ and $\SSS_1=\{(f_S,+\infty) : S\in\SSS\}$. We observe that it suffices to find definable downward directed families $\FF_0$ and $\FF_1$ contained in $p$ such that $\FF_0$ is complete for $\SSS_0$ and $\FF_1$ is complete for $\SSS_1$. In that case let $\FF=\{F\cap F' : F\in\FF_0, \, F'\in \FF_1\}$. Then $\FF$ is clearly contained in $p$, in particular it does not contain the empty set, and is downward directed. Moreover, for any $S=(f,g)\in \SSS$, if $F_0\in\FF_0$ is such that $F_0\subseteq (-\infty, g)$ and $F_1\in\FF_1$ is such that $F_1\subseteq (f,+\infty)$, then $F_0\cap F_1 \subseteq S$, so $\FF$ is finer than $\SSS$. We prove the existence of $\FF_0$ and $\FF_1$. 

For any two $S,S'\in\SSS$ let $X(S,S')=\{x\in \pi(S)\cap \pi(S') : g_S(x)\leq g_{S'}(x)\}$. These sets are definable uniformly over $S, S'\in \SSS$ (formally over $\Omega^2$, if $\Omega$ is the index set of $\SSS$). Let $\XX=\{ X(S,S') : S,S'\in \SSS\}$. By induction hypothesis there exists a definable downward directed family $\BB\subseteq \pi(p)$ that is complete for $\XX$. Let $\FF_0=\{(B\times M) \cap (-\infty, g_S) : B\in \BB,\,S\in \SSS\}$. We claim that this family is downward directed, contained in $p$ and finer than $\SSS_0$. The last property is obvious. Moreover, for every $B\in \BB$, since $B\in \pi(p)$ then $B\times M \in p$ and so, for every $S\in \SSS$, $(B\times M)\cap (-\infty, g_S) \in p$. It remains to prove that $\FF_0$ is downward directed. 


Let us fix $B', B''\in \BB$ and $S, S'\in \SSS$. Note that $\pi(S)\cap \pi(S')$ is covered by $X(S,S')$ and $X(S',S)$. Since $\BB$ is a subset of $\pi(p)$ and complete for $\XX$ it follows that there exists $B'''\in \BB$ satisfying either $B'''\subseteq X(S,S')$ or $B'''\subseteq X(S',S)$ (fact~\ref{fact:downward_directed_family_3}). Suppose without loss of generality that it is the former. Then $(B'''\times M)\cap (-\infty, g_{S})\subseteq (-\infty, g_{S'})$. Let $B \in \BB$ be such that $B\subseteq B'\cap B'' \cap B'''$, then $(B\times M)\cap (-\infty, g_S) \subseteq (B'\times M)\cap (-\infty, g_S) \cap (B''\times M)\cap (-\infty, g_{S'})$. We have shown that $\SSS_0$ is downward directed. 

The construction of $\SSS_1$ is analogous. This completes the proof of the lemma.     
\end{proof}

\begin{remark}
Let $p$ be a definable type and $\SSS\subseteq p$ be an $A$-definable family of sets. Is there always an $A$-definable downward directed family $\FF\subseteq p$ finer than $\SSS$? 

The negative answer is given by the following counterexample. Consider $\SSS=\{ M\setminus \{t\} : t\in M\}$, a family which is $\emptyset$-definable, and let $p$ be a non-realised type with basis $\{(s_0,t) : t>s_0\}$ for some fixed $s_0$ with $s_0\notin\text{dcl}(\emptyset)$. Suppose there exists $\FF$ as above and $\emptyset$-definable. Consider $B=\{ s\in M : \forall F\in \FF \,\exists t>s\, (s,t)\subseteq F\}$. Clearly $B$ is $\emptyset$-definable. Since $\FF\subseteq p$ the point $s_0$ is in $B$. If $B$ is finite then $s_0\in dcl(\emptyset)$, contradiction. Suppose that $B$ contains an interval $I$. By definition of $B$ and o-minimality this means that every set $F\in \FF$ must satisfy that $F\cap I$ is cofinite in $I$. By uniform finiteness there is an $m$ such that $|I\setminus F|<m$ for every $F\in \FF$. This however contradicts that $\FF$ is downward directed and finer than $\SSS$. 

On the other hand, in the next section we prove Proposition~\ref{prop:param_version_having_FTT}, which implies that any $A$-definable family $\SSS$ that extends to a definable type $p$ also extends to an $A$-definable type. By virtue of Lemma~\ref{lemma:existence_complete_directed_families_2.1}, it follows that there exists an $A$-definable downward directed family $\FF$ finer than $\SSS$ (however $\FF$ might not be a subset of $p$). For example, given $\SSS=\{ M\setminus \{t\} : t\in M\}$ as in the above counterxample one may choose $\FF=\{ (t,\infty) : t\in M\}$, i.e. the basis of a type at infinity. 
\end{remark}


We may now prove the main result of this section. 


\begin{theorem}\label{thm:types}
Let $\SSS\subseteq \mathcal{P}(M^n)$ be a uniform family of sets. The following are equivalent (with and without the definability condition in parentheses).
\begin{enumerate}[(1)]
\item \label{itm1: thm_types}$\SSS$ extends to a (definable) $n$-type.
\item \label{itm2: thm_types} There exists a uniform (definable) downward directed family finer than $\SSS$. 
\item \label{itm3: thm_types}$\SSS$ extends to an $n$-type with a uniform (definable) basis of cells.  
\end{enumerate}
\end{theorem}
\begin{proof}
$\eqref{itm1: thm_types}\Rightarrow\eqref{itm2: thm_types}$ is given by Lemma~\ref{lemma:existence_complete_directed_families_2.1}, $\eqref{itm2: thm_types}\Rightarrow \eqref{itm3: thm_types}$ by Lemma~\ref{lemma:existence_complete_directed_families}, $\eqref{itm3: thm_types}\Rightarrow\eqref{itm1: thm_types}$ is trivial. 
\end{proof}


\begin{remark}\label{remark:type_topology}
Theorem~\ref{thm:types} highlights an interesting property of types with uniform bases in o-minimal structures. Consider the \textquotedblleft uniform topology" on $S_n(M)$ where basic open sets $A_\SSS$ are indexed by uniform families $\SSS\subseteq \defn$, where $A_\SSS=\{p\in S_n(M) : \SSS\subseteq p\}$. This topology is clearly finer than the usual Stone topology among types. Note that, in this topology, the types with a uniform basis are precisely the isolated types. In particular, by o-minimality, all $1$-types are isolated (so the space $S_1(M)$ is not compact). Moreover, for every $A\subseteq M$, the set of all $A$-definable types is closed. In the context of this topology, Theorem~\ref{thm:types} states that isolated types are dense in $S_n(M)$ and isolated definable types are dense in the closed subspace of all definable types.   
\end{remark}


\begin{remark} \label{remark:parameter_thm_type}
Recall that Appendix~\ref{section:appendix B} shows that a parameter version of Lemma~\ref{lemma:existence_complete_directed_families} is not possible. That is, it is not true in general that every $A$-definable downward directed family of sets extends to an $A$-definable type with a uniform basis. We end this section however by proving the next best thing (Proposition~\ref{prop:parameter_version_existence_downward_directed_families}), that every $A$-definable downward directed family extends to an $A$-definable type. This result may be interpreted as a weak parameter version of Lemma~\ref{lemma:existence_complete_directed_families}. Added to the parameter information in the statement of Lemma~\ref{lemma:existence_complete_directed_families_2.1}, it
allows us to conclude that Theorem~\ref{thm:types} still holds if we substitute the word ``definable" by ``$A$-definable", with the exception of implications \eqref{itm1: thm_types}$\Rightarrow$\eqref{itm3: thm_types} and \eqref{itm2: thm_types}$\Rightarrow$\eqref{itm3: thm_types}. 
\end{remark}

The following method to construct types will be useful in Proposition~\ref{prop:parameter_version_existence_downward_directed_families} and further on. 

\begin{definition}\label{definition:construction_types}
Let $p$ be a $n$-type and $f:M^n\rightarrow M$ be a definable function with $dom(f)\in p$. We define three $n+1$-types.

Let $f|_p$ be the type of sets $S$ such that $\pi(S\cap graph(f))\in p$.

Let $f^{+}|_p$ be the type of sets $S$ such that $\{x\in dom(f) : \{x\}\times (f(x),s) \subseteq S \text{ for some } s>f(x)\}\in p$.

Let $f^{-}|_p$ be the type of sets $S$ such that $\{x\in dom(f) : \{x\}\times (s,f(x)) \subseteq S \text{ for some } s<f(x)\}\in p$.

The definitions of $f^{+}|_p$ and $f^{-}|_p$ hold too if $f$ is, respectively, the constant function $-\infty$ and $+\infty$.
\end{definition}

The next lemma without proof follows by routine application of o-minimality and knowledge of the Fiber Lemma for o-minimal dimension \cite[Chapter 4, Proposition 1.5 and Corollary 1.6]{dries98}.

\begin{lemma}
Given an $n$-type $p$ and definable function $f:M^n\rightarrow M$ with $dom(f)\in p$, the types $f|_p$, $f^+|_p$ and $f^-|_p$ are well defined. If $p$ is $A$-definable and $f$ $B$-definable then these types are $AB$-definable. If $p$ is $d$-dimensional then $f|_p$ is $d$-dimensional and $f^{+}|_p$ and $f^{-}|_p$ are $d+1$-dimensional. 
\end{lemma}


\begin{proposition}\label{prop:parameter_version_existence_downward_directed_families}
Let $A\subseteq M$ and let $\SSS\subseteq \defn$ be an $A$-definable downward directed family of sets. Then $\SSS$ extends to an $A$-definable type. 
\end{proposition}
\begin{proof}
We proceed by induction on $n$. If $n=1$ then this is shown by the base case in the proof of Lemma~\ref{lemma:existence_complete_directed_families}. Suppose that $n>1$. 

By induction hypothesis let $p$ be an $A$-definable type extending $\pi(\SSS)$. For every $S\in\SSS$ we define $f_S:\pi(S)\rightarrow \exR$ by $f_S(x)=\sup S_x$, where $S_x$ is the corresponding fiber $\{t : \al x,t\ar\in S\}$. Let $\FF=\{f_S : S\in\SSS\}$. This family is clearly $A$-definable. Let $\lleq$ be the $A$-definable total preorder induced by $p$ on $\FF$ described in Subsection~\ref{remark:preorder}. We consider two cases.

\textbf{Case $1$:} $\FF$ has a minimum with respect to $\lleq$. Let $\SSS_{\text{min}}\subseteq \SSS$ be the family of $S\in\SSS$ such that $f_S$ is a minimum. Let $\FF_{\text{min}}=\{ f_S : S\in \SSS_{\text{min}}\}$. We prove that, for any given $S\in \SSS_{\text{min}}$, $\SSS$ extends to either $f_S|_p$ or $f_S^-|_p$, and show that these types are $A$-definable. We make use of a claim and a fact.

\begin{claim}\label{claim:param_version_existence_downward_directed_families}
For every $g,h\in \FF_{\text{min}}$, $g|_p=h|_p$ and $g^-|_p=h^-|_p$. We denote these types $\FF_\text{min}|_p$ and $\FF^-_\text{min}|_p$ respectively. 
\end{claim} 

Since the family $\FF_{\text{min}}$ is $A$-definable, it follows from Claim~\ref{claim:param_version_existence_downward_directed_families} and definition~\ref{definition:construction_types} that $\FF_\text{min}|_p$ and $\FF^-_\text{min}|_p$ are $A$-definable. We prove the claim.    

Let $g,h\in \FF_{\text{min}}$. Obeserve that, because $g\lleq h$ and $h\lleq g$, we have that $\{ x : g(x)=h(x)\} \in p$.
Let $Z\in \defn$ be such that $Z\in g|_p$. Then $\{ x : \al x,h(x) \ar \in Z\} \supseteq \{ x: g(x)=h(x),\, \al x,g(x) \ar \in Z\}\in p$. So $Z\in h|_p$. Analogously one proves that $g^-|_p=h^-|_p$. 

The following fact which is easy to check.

\begin{fact}\label{fact:param_version_existence_downward_directed_families}
For every $S\in \SSS$, exactly one of the sets $\{ x : f_S(x)\in S_x\}$ and $\{x : f_S(x)\notin S_x\}$ must belong in $p$. If it is the former then $S\in f_S|_p$. If it is the latter then, by definition of $f_S$ and o-minimality, we have that $S\in f^-_S|_p$. 
\end{fact}  

We assume that $\SSS$ does not extend to $\FF_\text{min}|_p$ and prove that it extends to $\FF^-_\text{min}|_p$. 
Let $S^*\in \SSS$ be such that $S^*\notin \FF_\text{min}|_p$. Let $S'\in \SSS_{\text{min}}$. For any $S\in\SSS$, by downward directedness there exists $S''\in \SSS$ with $S''\subseteq S\cap S'\cap S^*$. Since $S''\subseteq S'$ we have that $S''\in\SSS_{\text{min}}$. Since $S''\subseteq S^*$, $S''\notin \FF_\text{min}|_p= f_{S''}|_p$. By Fact~\ref{fact:param_version_existence_downward_directed_families}, it follows that $S''\in f^-_{S''}|_p$. We conclude that $S\in f^-_{S''}|_p= \FF^-_\text{min}|_p$. 

\textbf{Case $2$:} $\FF$ does not have a minimum with respect to $\lleq$. In this case consider the $n$-type $q$ given by all sets $F\in \defn$ such that $F\in h|_p$ for every $h$ in some subfamily $\FF'\subseteq \FF$ that is unbounded from below with respect to $\lleq$. Note that, since $\FF$, $\lleq$ and $p$ are $A$-definable, $q$ is $A$-definable too. Clearly, for every $F\in\defn$, either $F\in q$ or $M^n\setminus F \in q$. We prove that $q$ is consistent, which allows us to conclude that $q$ is well defined, and then show that $\SSS\subseteq q$. We use the following two short claims.    

\begin{claim}\label{claim:param_version_existence_downward_directed_families_1}
Let $C=(f,g)$ be a cell in $M^n$ with $\pi(C)\in p$ and let $h\in \FF$. If $f \prec h \prec g$ then $C \in h|_p \cap h^-|_p \cap h^+|_p$. 
\end{claim}

By definition of $\lleq$, we have that $\{ x : f(x)<h(x)<g(x)\}\in p$. So $\{ x : \al x, h(x)\ar \in C\} \in p$, meaning that $C\in h|_p$. Similarly one shows what $C\in h^-|_p$ and $C\in h^+|_p$.

\begin{claim}\label{claim:param_version_existence_downward_directed_families_2}
Let $F\in \defn$. If $F\in q$ then it contains a cell $C=(f,g)$ with $\pi(C)\in p$ such that $f\prec \FF$ and $h\prec g$ for some $h\in\FF$. 
\end{claim}

Note that, from Claims~\ref{claim:param_version_existence_downward_directed_families_2} and~\ref{claim:param_version_existence_downward_directed_families_1} it follows that, if $F$ is a set in $q$, then there is some $h\in \FF$ such that, for every $h'\in \FF$ with $h'\lleq h$, it holds that $F\in h'|_p$. Since $\FF$ does not have a minimum, it follows easily that $q$ is consistent. We prove Claim~\ref{claim:param_version_existence_downward_directed_families_2}.


Let $F\in \defn$. Note that, by cell decomposition, either $F$ or $M^n\setminus F$ must contain a cell $C=(f,g)$ with $\pi(C)\in p$ and satisfying that $f\prec \FF$ and $h\prec g$ for some $h\in \FF$. We assume that $F$ does not contain such a cell $C$ and hence fix $C=(f,g)$ with $C\subseteq M^n\setminus F$. We also fix $h\in \FF$ such that $C=(f,g)$ with $f\prec \FF$ and $h\prec g$. By Claim~\ref{claim:param_version_existence_downward_directed_families_1} we have that, for every $h'\in \FF$ with $h'\lleq h$, it holds that $C\in h'|_p$, and in particular $F\notin h'|_p$. It follows that $F\notin q$. This proves Claim~\ref{claim:param_version_existence_downward_directed_families_2}.

Finally, we show that $\SSS\subseteq q$. Towards a contradiction suppose that there exists some $S^*\in \SSS$ with $S^*\notin q$. By Claim~\ref{claim:param_version_existence_downward_directed_families_2} applied to the set $M^n\setminus S^*$ let us fix a cell $C=(f,g)$ with $C\cap S^*=\emptyset$ and $\pi(C)\in p$ and some $S'\in \SSS$ such that $f\prec \FF$ and $f_{S'}\prec g$. By downward directedness let $S\in \SSS$ be such that $S\subseteq S^*\cap S'$. Since $S\subseteq S^*$ we have that $C\cap S=\emptyset$. Since $S\subseteq S'$ it holds that $f_S \lleq f_{S'}\prec g$. But then, by Claim~\ref{claim:param_version_existence_downward_directed_families_1}, $C\in f_S|_p \cap f^-_S|_p$. However, by Fact~\ref{fact:param_version_existence_downward_directed_families}, $S\in f_S|_p \cup f^-_S|_p$, contradicting that $C\cap S=\emptyset$.  
\end{proof}

\begin{remark}
In Appendix~\ref{section:appendix B} Will Johnson constructs, in a structure with the same definable sets as $(\mathbb{R},\leq,\cdot)$ (but a different language), a $\emptyset$-definable downward directed family of sets $\SSS$ that does not extend to an $\emptyset$-definable type with a uniform basis. By Proposition~\ref{prop:parameter_version_existence_downward_directed_families}, such a family $\SSS$ does extend to a $\emptyset$-definable type. 
In the context of the uniform type topology described in Remark~\ref{remark:type_topology}, this means that, in general, the set of isolated $\emptyset$-definable types is not dense in the space of all $\emptyset$-definable types. In contrast, in Remark~\ref{remark:type_topology} we note that, by Theorem~\ref{thm:types} (implication \eqref{itm1: thm_types}$\Rightarrow$\eqref{itm3: thm_types}), this density result does hold for the space of ($M$-)definable types. By the usual argument involving passing to an elementary substructure, the result also holds for $A$-definable types, for any $A\subseteq M$, whenever $\MM$ has definable choice functions. 
\end{remark}

\begin{remark}
In~\cite{dolich04} Dolich gives a definition for a formula being ``good over $A$" that is equivalent to the property that the formula extends to an $A$-definable type which is given inductively by one of the three constructions described in Definition~\ref{definition:construction_types}, where all the functions involved are $A$-definable\footnote{In conversation with Sergei Starchenko he informally referred to these types are ``strongly $A$-definable".}. It follows from Dolich's work that, in o-minimal expansions of ordered fields, any $A$-definable downward directed family of sets extends in fact to an $A$-definable type of this kind. It seems plausible that this is true in all o-minimal structures (perhaps with the assumption of having elimination of imaginaries). This would imply in particular that the type construction in Case $2$ of the proof of Proposition~\ref{prop:parameter_version_existence_downward_directed_families}, which will be used again in the proof of Lemma~\ref{lemma:FIP_transversal}, is not necessary. Using Lemma~\ref{lemma:existence_complete_directed_families_2.1}, it also yields a strong density result for $A$-definable types of this kind. 
\end{remark}

The main result of this section, Theorem~\ref{thm:types}, uses both Lemmas~\ref{lemma:existence_complete_directed_families} and~\ref{lemma:existence_complete_directed_families_2.1} to make an observation about density of types with a uniform basis. On the other hand, in terms of their applicability in further sections, Lemmas~\ref{lemma:existence_complete_directed_families} and~\ref{lemma:existence_complete_directed_families_2.1} are quite distinct. Lemma~\ref{lemma:existence_complete_directed_families}, and specifically the fact that any definable downward directed family extends to a definable type, will be used in proving Lemma~\ref{lemma:order_transversals} and Proposition~\ref{prop:curve_compact_iff_directed_compact}. Meanwhile Lemma~\ref{lemma:existence_complete_directed_families_2.1}, and specifically the fact that any definable family that extends to a definable type admits a finer definable downward directed family, is used to prove Proposition~\ref{prop:compactnes_transversals} (in particular the fact that definable compactness implies type-compactness) and, together with Proposition~\ref{prop:parameter_version_existence_downward_directed_families}, to prove Proposition~\ref{prop:param_version_having_FTT}. Because of their different applications, we will often reference the lemmas directly instead of the theorem. 

\section{Transversals}\label{section:transversals}
Throughout this section we fix a monster model $\U$. For a definable set $X$, let $X(\U)$ denote the interpretation on $X$ in $\U$. Our main result is Theorem~\ref{thm:finite_tame_transversals}, stating that, if a definable family has the $(m,n)$-property for some large enough $n$ and $m\geq n$, then it partitions into finitely many subfamilies, each of which extends to a definable type. We begin by introducing some terminology to refer to the property of a family having such a partition. 



\begin{definition}
Let $\SSS\subseteq \defn$. A transversal $F$ of $\{S(\U) : S\in \SSS\}$ is called \emph{tame} if it belongs in some tame extension of $\MM$. We abuse terminology by calling $F$ a \emph{tame tranversal of $\SSS$}.  
\end{definition}
Let $\SSS$ be a uniform family. We are interested in the property that $\SSS$ admits a tame transversal that is finite. We give an equivalent definition of this property.  
\begin{definition}\label{dfn:ftt}
 Let $\SSS$ be a uniform family and $A\subseteq M$. We say that $\SSS$ has a finite tame transversal (FTT) if there are finitely many subfamilies $\SSS_1,\ldots,\SSS_k \subseteq \SSS$ with $\cup_i \SSS_i = \SSS$ such that, for every $1\leq i\leq k$, the family $\SSS_i$ extends to a definable type.

\end{definition}
By one direction of the Marker-Steinhorn Theorem \cite[Theorem 2.1]{mark_stein_94}, every definable type is realized in some tame elementary extension of $\MM$. So both definitions of having FTT are equivalent. 

Before beginning the work towards Theorem~\ref{thm:finite_tame_transversals}, we show that the types witnessing that a family $\SSS$ has FTT can always be chosen definable over the same parameters as $\SSS$. This corresponds to Proposition~\ref{prop:param_version_having_FTT}. We prove it using Lemma~\ref{lemma:existence_complete_directed_families_2.1} and Proposition~\ref{prop:parameter_version_existence_downward_directed_families} from last section on the relationship between types and downward directed families of sets. First we need a short technical lemma. 

\begin{lemma}\label{lemma:remark_definable_types}
Let $A\subseteq M$. Let $p\in S_n(U)$ be a type definable over a tuple $b\in U^{m}$. Suppose the $\tp(b/M)$ is $A$-definable. Then the restriction of $p$ to $M$ (the partial type of $M$-definable sets in $p$) is $A$-definable. 
\end{lemma}
\begin{proof}
Let $\varphi(x,y)$ be a partitioned formula with $l(x)=n$ and $\psi(y,b)$ be the definition of $p|_\varphi$.  Since the type of $b$ over $M$ is $A$-definable the set $\psi(M^n, b)$, which corresponds to all $y\in M^{l(y)}$ such that $\varphi(x,y)\in p$, is $A$-definable. 
\end{proof}

\begin{proposition}[Parameter chatacterization of having FTT]\label{prop:param_version_having_FTT}
Let $A\subseteq M$ and let $\SSS$ be an $A$-definable family of sets with a finite tame transversal of size $n$. Then $\SSS$ can be paritioned into $n$ subfamilies $\SSS_1,\ldots, \SSS_n$ such that, for every $1\leq i \leq n$, the family $\SSS_i$ extends to an $A$-definable type. 
\end{proposition}
\begin{proof}
Let $\SSS\subseteq Def(M^m)$ be an $A$-definable family of sets having a finite tame transversal of size $n$. By Lemma~\ref{lemma:existence_complete_directed_families_2.1} there exists $n$ formulas $\varphi_1(x,y_1,a_1), \ldots, \varphi_n(x,y_n,a_n)$, with parameters $a_1,\ldots a_n$ respectively, such that, for every $i$, the family $\{\varphi(M^m, u_i, a_i) : u_i\in M^{l(y_i)}\}$ is a downward directed and, moreover, for every $S\in \SSS$, there is some $i$ and $u_i\in M^{l(y_i)}$ such that $\varphi_i(M^m,u_i,a_i)\subseteq S$.       

Let $X$ denote the set of points $\al b_1 \ldots, b_n\ar$ in $\prod_i M^{l(a_i)}$ satisfying that the formulas $\varphi_1(x,y_1,b_1), \ldots, \varphi_n(x,y_n,b_n)$ have the properties described above for $a_i=b_i$. That is, each $\varphi_i(x,y_i,b_i)$ describes a downward directed family and, for every $S\in \SSS$, there is one such family that is finer than $\{S\}$. This set is nonempty, since $\al a_1, \ldots, a_n\ar \in X$. Note that $X$ is $A$-definable.

By Proposition~\ref{prop:parameter_version_existence_downward_directed_families} (applied to a family with only one set) let $p$ be an $A$-definable type extending $X$ and $\al\mathbf{b_1},\ldots, \mathbf{b_n}\ar$ a realization of $p$ in $\prod_i U^{l(a_i)}$. Note in particular that, for every $i$, $\tp(\mathbf{b_i}/M)$ is $A$-definable. By Proposition~\ref{prop:parameter_version_existence_downward_directed_families}, let $p_i\in S_{m}(U)$ be a $\mathbf{c_i}$-definable type extending $\{\varphi_i(x,u_i, \mathbf{c_i}) : u_i\in U^{l(y_i)} \}$. Note that, by definition of $X$, every set $S\in \SSS$ belongs in a type $p_i$ for some $i$. 

Finally for any $i$ let $q_i$ be the restriction of $p_i$ to $M$. By Lemma~\ref{lemma:remark_definable_types} these restrictions are $A$-definable. Moreover every $S\in \SSS$ belongs in $q_i$ for some $i$. 
\end{proof}

The following is the main result of this section. 

\begin{theorem}\label{thm:finite_tame_transversals}
Let $\SSS\subseteq \mathcal{P}(M^k)$ be a definable family of nonempty sets. Let $n>\dim \cup\SSS$. The following are equivalent.   
\begin{enumerate}[(1)]
\item \label{itm:thm_transversals_1} There exists some $m\geq n$ such that $\SSS$ has the $(m,n)$-property.
\item \label{itm:thm_transversals_2} $\SSS$ has a finite tame transversal.   
\end{enumerate}  
\end{theorem}

We precede the proof of Theorem~\ref{thm:finite_tame_transversals} with two lemmas and a short proposition. First however we must introduce some notation on cuts. 

Given a definable totally preordered set $(X, \lleq)$, a cut $(P,Q)$ of $(X, \lleq)$ is a partition of $X$ into two sets $P$ and $Q$ such that $x\prec y$ for every $x\in P$ and $y\in Q$. Such a cut $(P,Q)$ is definable if $P$ (equivalently $Q$) is definable. All the cuts we consider satisfy that either $P$ does not have a maximum or $Q$ does not have a minimum with respect to $\lleq$. That is, for every $x\in P$ and $y\in Q$ the interval $(x,y)_\lleq$ is nonempty. Consequently we may always associate to a cut $(P,Q)$ the partial type $\{ (x,y)_\lleq : x\in P, y\in Q\}$, and by means of this association we often refer to cuts as types. In particular we say that $z\in X(\U)$ realizes $(P,Q)$ if $(x,y)_{\lleq}\in \tp(z/M)$ for every $x\in P$ and $y\in Q$. 


Without context a cut is a cut in $(M,\leq)$. In this setting non-definable and definable cuts are denoted, respectively, irrational and rational cuts in~\cite{pillay_stein_87}, and simply cuts and noncuts in other sources such as~\cite{mark_stein_94} and~\cite{dolich04}.
By o-minimality, any non-isolated $1$-type over $M$ is uniquely characterised by the unique cut it realizes in $(M,\leq)$.
 
It is worth noting that, in the statement of the Lemma~\ref{lemma:order_transversals} below, the existence of the uniform bound $l$ is redundant, since it follows from o-minimality.

\begin{lemma}\label{lemma:order_transversals}
Let $(X, \lleq)$, $X\subseteq M^k$, be a definable totally preordered set, and $l<\omega$. Let $\SSS$ be a definable family of nonempty subsets of $X$, all of which are finite union of at most $l$ intervals with respect to $\lleq$ (including degenerate intervals of the form $[x,x]_\lleq$) with endpoints in $X \cup\{-\infty, +\infty\}$ (where $-\infty$ and $+\infty$ have the natural interpretation with respect to $\lleq$). Then exactly one holds:
\begin{enumerate}[(1)]
\item $\SSS$ has a finite tame transversal. \label{itm:order_transversals_1}
\item There exists $\SSS'\subseteq \SSS$ an infinite subfamily (not necessarily definable) of pairwise disjoint sets. \label{itm:order_transversals_2}
\end{enumerate}
\end{lemma}

Lemma~\ref{lemma:order_transversals} stands on its own right as a powerful observation regarding all definable total preorders in an o-minimal structure. Note that, by letting $(X,\lleq)=(M, \leq)$, it implies a strengthening of the (base) case $\SSS\subseteq \mathcal{P}(M)$ of Theorem~\ref{thm:finite_tame_transversals}, since it shows that, if $\SSS$ is a definable family of nonempty sets with the $(\omega,2)$-property, then it has a finite tame transversal. We do not know whether or not this strengthening holds in higher dimensions (see Question~\ref{question:generalizing_FTT_theorem}). 

The fact that Lemma~\ref{lemma:order_transversals} applies to any definable totally preordered set $(X,\lleq)$, and in particular to preorders like those described in Subsection~\ref{remark:preorder}, will be useful later on in constructing an inductive proof in Lemma~\ref{lemma:FIP_transversal}.

\begin{proof}[Proof of Lemma~\ref{lemma:order_transversals}]
Clearly~\eqref{itm:order_transversals_1} and~\eqref{itm:order_transversals_2} are mutually exclusive. We assume the negation of \eqref{itm:order_transversals_1} and prove \eqref{itm:order_transversals_2}. Throughout the proof let interval refer to an interval in $(X, \lleq)$ with endpoints in $X \cup\{-\infty, +\infty\}$. 

Hence let us assume that $\SSS$ does not admit a finite tame transversal. Let $\SSS=\{\varphi(u,M^k) : u\in \Omega\}$. Then by first order logic compactness there exists $c\in \Omega(\Mon)$, $\mathbb{S}=\varphi(c,U^k)$, with the property that, for every definable $\Omega'\subseteq \Omega$, if $\{ \varphi(u,M^k) : u\in\Omega'\}$ extends to a definable type, then $c\notin \Omega'(\Mon)$. In particular this holds for $\Omega'=\Omega_x=\{u\in \Omega : \MM\models\varphi(u,x)\}$ for every $x\in X$, and so $\mathbb{S}\cap X=\emptyset$. Moreover note that, if $\sim$ denotes the equivalence relation induced by the preorder $\lleq$ on $X$, then by definition every $S\in\SSS$ is compatible with the equivalence classes of $\sim$, i.e. for every $x, y\in X$ if $x\sim y$ and $x\in S$ then $y\in S$ too. This must also be satisfied by $\badS$ with respect to the interpretation of $\lleq$ on $X(\Mon)$ (which we also denote $\lleq$). So each equivalence class in $X(\Mon)$ which includes points from $X$ must be disjoint from $\badS$. Hence each point $x\in\mathbb{S}$ induces a cut $(P,Q)$ in $X$ by $P=\{y\in X : y\prec x\}$ and $Q=X\setminus P$. 


Let $(P,Q)$ and $(P',Q')$ be different cuts in $M$, with, say, $P\subsetneq P'$, and let $x_{(P,P')}\in P'\setminus P$. Since $x_{(P,P')}\in X$ we have that $x_{(P',P)}\notin \mathbb{S}$ and so, if $x$ and $y$ are points in $\mathbb{S}$ that realize $(P,Q)$ and $(P',Q')$ respectively, then these points do not belong in the same subinterval of $\mathbb{S}$. In other words, each subinterval of $\badS$ contains only realizations from a single cut in $(X,\lleq)$. Since every set in $\mathbb{S}$ is union of at most $l$ intervals, there are at most $l$ many distinct cuts $(P_i,Q_i)$, $0\leq i \leq m$, such that every point in $\mathbb{S}$ realizes one of these cuts. 

For simplicity we assume that that all the cuts $(P_i,Q_i)$ are such that $P_i\neq \emptyset$ and $Q_i\neq\emptyset$. The proof readily adapts to the other case. Let the indexing be such that $(P_i,Q_i)$, for $0\leq i \leq n$, are the cuts that are non-definable, and $(P_i,Q_i)$, for $n<i \leq m$, are definable.  


For every $0\leq i \leq m$ let us fix $a_i\in P_i$ and $b_i\in Q_i$ such that $\{a_i, b_i\}\subseteq P_j$ or $\{a_i,b_i\} \subseteq Q_j$ for every $j\neq i$. We show that there exists $S\in\SSS$ with $S\subseteq \cup_i (a_i, b_i)_\lleq$ such that, for every $0\leq i\leq m$, there are $a'_i\in P_i$, $b'_i\in Q_i$ with $a_i\lleq a'_i \prec b'_i \lleq b_i$ such that $(a'_i,b'_i)_\lleq\cap S=\emptyset$ (in particular $S\cap \mathbb{S}=\emptyset$). We may then apply again the result with parameters $a'_i$, $b'_i$ in place of $a_i$ and $b_i$ and reach a second set $S'\in \SSS$ that will be disjoint from $S$. Repeating this process yields a countably infinite pairwise disjoint subfamily of $\SSS$ as desired. 

Let $(P,Q)$ be a cut. We say that a definable set $S$ is disjoint from $(P,Q)$ if there exists any $a\in P$ and $b\in Q$ such that $(a, b)_\lleq\cap S=\emptyset$. Hence we must find $S\in \SSS$ contained in $\cup_i (a_i, b_i)_\lleq$ that is disjoint from every cut $(P_i,Q_i)$, $0\leq i \leq m$.  

We focus first on the definable cuts.  
\begin{claim}\label{claim:order_transversals2}
Let $(P,Q)$ be a definable cut in $X$, then $\badS$ is disjoint from the cut $(P(\Mon),Q(\Mon))$ in $X(\Mon)$.
\end{claim}
For any $a\in P(\Mon)$ and $b\in Q(\Mon)$ let $(a,b)_{\lleq(\Mon)}=\{c'\in X(\Mon) : a \lleq c' \lleq b\}$.
We show that there is no $a\in P(\Mon)$ and $b\in Q(\Mon)$ with $(a,b)_{\lleq(\Mon)} \subseteq \badS$.
Since $\badS$ is union of finitely many intervals with respect to $\lleq$ it follows that there exists $a\in P(\Mon)$ and $b\in Q(\Mon)$ with $(a,b)_{\lleq(\Mon)}\cap \badS=\emptyset$. 

Suppose that there is $a\in P(\Mon)$ and $b\in Q(\Mon)$ with $(a,b)_{\lleq(\Mon)} \subseteq \badS$. Let $\Sigma\subseteq \Omega$ be denote the set of index elements $u$ such that $(a,b)_\lleq \subseteq \varphi(u,M^k)$ for some $a\in P$, $b\in Q$. Since the cut is definable then so is $\Sigma$. Moreover $c$, where $\badS=\varphi(c,U^k)$, belongs in $\Sigma(\Mon)$. On the other hand the definable family $\{(a,b)_{\lleq} : a\in P, b\in Q\}$ is downward directed and so, by Lemma~\ref{lemma:existence_complete_directed_families}, extends to a definable type. This type must include $\{\varphi(u,M^k) :u\in \Sigma\}$. But then by construction $c\notin \Sigma(\Mon)$, contradiction. This proves Claim~\ref{claim:order_transversals2}.  

We will need the following fact.
\begin{fact}\label{fact:order_transversals}
For all choice of $x_i\in P_i$, $y_i\in Q_i$, $0\leq i \leq m$, there exists $S\in \SSS$ such that $S\subseteq \cup_i (x_i,y_i)_\lleq$ and is disjoint from every cut $(P_j, Q_j)$, for $n<j\leq m$. 
\end{fact}
By Claim~\ref{claim:order_transversals2} this is witnessed by $\badS$ in $\Mon$, so it also holds in $\MM$. 

Finally we require a claim regarding non-definable cuts. 
\begin{claim}\label{claim:order_transversals}
Let $(P,Q)$ be a non-definable cut. If $F\subseteq X$ is a definable set such that, for every $x\in P$, there is $y\in F$ with $x\lleq y$, then $F\cap Q\neq \emptyset$. Similarly if, for every $x\in Q$, there is $y\in F$ with $y \lleq x$, then $P\cap F\neq \emptyset$. 
\end{claim}
To prove Claim~\ref{claim:order_transversals} note that, since $F\subseteq X$ is definable, then its downward closure with respect to $\lleq$, namely $F'=\{x\in X :x\lleq y \text{ for some } y\in F\}$, is definable too. If $F\cap P$ is cofinal in $P$ then $P \subseteq F'$, so if $P$ is not definable then it must be that $F'\cap Q\neq \emptyset$, hence $F\cap Q\neq \emptyset$. The case where $P\cap Q$ is unbounded from below in $Q$ is analogous. This proves Claim~\ref{claim:order_transversals}.

For any $x_0, y_0, \ldots, x_n, y_n\in X$ let $\rel(x_0, y_0, \ldots, x_n, y_n)$ denote the formula, including parameters $a_i, b_i$ for $i>n$, asserting that there exists $S\in\SSS$ with $S\subseteq \bigcup_{0\leq i\leq n} (x_i, y_i)_\lleq \cup \bigcup_{j>n} (a_j, b_j)_\lleq$ such that $S$ is disjoint from the cuts $(P_j, Q_j)$ for $n<j$. Note that Fact~\ref{fact:order_transversals} implies that $\rel(x_0, y_0, \ldots, x_n, y_n)$ holds for all choice of $x_i\in P_i$ and $y_i\in Q_i$, $0< i\leq n$. We define by reverse recursion a family of definable sets $F_i, G_i \subseteq X$, for $0< i \leq n$, as follows. Let
\begin{align*}
F_n(x_0, y_0, \ldots, x_{n-1}, y_{n-1})&=\{ y \succcurlyeq a_n : \rel(x_0, y_0, \ldots, x_{n-1}, y_{n-1}, a_n, y)\} \text{ and} \\
G_n(x_0, y_0, \ldots, x_{n-1}, y_{n-1})&=\{ x \lleq b_n :  \rel(x_0, y_0, \ldots, x_{n-1}, y_{n-1}, x, b_n)\}.
\end{align*}

For $0< i <n$ let 
\begin{align*}
F_i(x_0, y_0, \ldots , x_{i-1}, y_{i-1})=\{& y \succcurlyeq a_{i} : \exists x'\in F_{i+1}(x_0, y_0, \ldots, x_{i-1}, y_{i-1}, a_{i}, y) 
\\
&\exists y'\in G_{i+1}(x_0, x_1, \ldots, x_{i-1}, y_{i-1}, a_{i}, y),\, x'\lleq y'\}\\
\omit\text{and} 
\\ 
G_i(x_0, y_0, \ldots , x_{i-1}, y_{i-1})=\{& x \lleq b_{i} : \exists x'\in F_{i+1}(x_0, y_0, \ldots, x_{i-1}, y_{i-1}, x, b_i) 
\\
&\exists y' \in G_{i+1}(x_0, y_0, \ldots, x_{i-1}, y_{i-1}, x, b_i),\, x'\lleq y'\}.
\end{align*}

Clearly these sets are definable. For every $0< i \leq n$ we prove the following.  
\begin{enumerate}
\item[$(\text{I}_{i})$] \label{itm:proof_order_trans_1} For all choice of $x_j\in P_j$, $y_j\in Q_j$, $0 \leq j <i$, it holds that there exists $x\in F_i(x_0, y_0, \ldots, x_{i-1}, y_{i-1})$ and $y\in G_i(x_0, y_0, \ldots, x_{i-1}, y_{i-1})$ with $x\lleq y$. 
\item[$(\text{II}_{i})$] \label{itm:proof_order_trans_2}  For all choice of $a_j\lleq x_j \prec y_j \lleq b_j$, $0\leq j < i$, if there exists $x\in F_i(x_0, y_0, \ldots, x_{i-1}, y_{i-1})$ and $y\in G_i(x_0, y_0, \ldots, x_{i-1}, y_{i-1})$ with $x\lleq y$, then there exists $S\in \SSS$ with $S\subseteq \bigcup_{j < i} (x_j, y_j)_\lleq \cup \bigcup_{j\geq i} (a_j, b_j)_\lleq$ that is disjoint from any cut $(P_j, Q_j)$ for $j \geq i$. 
\end{enumerate}

We proceed by reverse induction on $i$, where in the inductive step $(\text{I}_{i})$ follows from $(\text{I}_{i+1})$, and $(\text{II}_{i})$ from $(\text{II}_{i+1})$. We will then derive the lemma from $(\text{I}_{1})$ and $(\text{II}_{1})$. 
\sloppy
Let $i=n$. Fix any $x_j\in P_j$ and $y_j\in Q_j$ for every $j < n$. By Fact~\ref{fact:order_transversals} and Claim~\ref{claim:order_transversals} $F_n(x_0, y_0, \ldots, x_{n-1}, y_{n-1})\cap P_n\neq \emptyset$ and $G_n(x_0, y_0, \ldots, x_{n-1}, y_{n-1})\cap Q_n\neq \emptyset$, so in particular there exists $x\in F_n(x_0, y_0, \ldots, x_{n-1}, y_{n-1})$ and $y\in G_n(x_0, y_0, \ldots, x_{n-1}, y_{n-1})$ with $x\lleq y$. This proves $(\text{I}_{n})$. To prove $(\text{II}_{n})$ let $a_j\lleq x_j \prec y_j \lleq b_j$ for $j < n$. If there exists $x\in F_n(x_0, y_0, \ldots, x_{n-1}, y_{n-1})$ and $y\in G_n(x_0, y_0, \ldots, x_{n-1}, y_{n-1})$ with $x\lleq y$ then it cannot be that $F_n(x_0, y_0, \ldots, x_{n-1}, y_{n-1})\subseteq Q_n$ and $G_n(x_0, y_0, \ldots, x_{n-1}, y_{n-1})\subseteq P_n$. Suppose that $F_n(x_0, y_0, \ldots, x_{n-1}, y_{n-1})\cap P_n\neq \emptyset$, the other case being similar. Then by definition of $F_n$ there exists $y\in P_n$ and some $S\in\SSS$ with $S\subseteq \cup_{0\leq j<n} (s_j, t_j)_\lleq \cup (a_n,y)_\lleq \cup \cup_{j>n} (a_j,b_j)_\lleq$ such that $F$ is disjoint from every cut $(P_j, Q_j)$, $j>n$. By construction $S$ is disjoint from $(P_n, Q_n)$ too, proving $(\text{II}_{n})$.

Suppose now that $i<n$. Fix $x_j\in P_j$ and $y_j\in Q_j$ for $j < i$. By $(\text{I}_{i+1})$ and Claim~\ref{claim:order_transversals} it holds that there exists $x\in F_i(x_0, y_0, \ldots, x_{i-1}, y_{i-1})\cap P_{i}$ and $y\in G_i(x_0, y_0, \ldots, x_{i-1}, y_{i-1})\cap Q_{i}$, satisfying in particular $x\lleq y$. This proves $(\text{I}_{i})$.  

Now fix $a_j\lleq x_j \prec y_j \lleq b_j$ for $j < i$. Once again if there are $x\in F_i(x_0, y_0, \ldots, x_{i-1}, y_{i-1})$ and $y\in G_i(x_0, y_0, \ldots, x_{i-1}, y_{i-1})$ with $x\lleq y$ then it must be that either $F_i(x_0, y_0, \ldots, x_{i-1}, y_{i-1})\cap P_i\neq \emptyset$ or $G_i(x_0, y_0, \ldots, x_{i-1}, y_{i-1})\cap Q_i\neq \emptyset$. We assume the former, being the proof given the latter analogous. Hence let $y\in F_i(x_0, y_0, \ldots, x_{i-1}, y_{i-1})\cap P_i$. By definition of $F_i$ there exists some $x'\in F_{i-1}(x_0, y_0, \ldots, x_{i-1}, y_{i-1}, a_i, y)$ and $y'\in G_{i-1}(x_0, y_0, \ldots, x_{i-1}, y_{i-1}, a_i, y)$ with $x'\lleq y'$. But then by $(\text{II}_{i+1})$ there exists $S\in \SSS$ such that $S\subseteq \bigcup_{0\leq j<i} (x_j, y_j)_\lleq \cup (a_{i}, y)_\lleq \cup \bigcup_{j>i} (a_j, b_j)_\lleq$, and $S$ is disjoint from any cut $(P_j, Q_j)$ for $j>i$. However note that, by construction, $S$ is also disjoint from the cut $(P_{i}, Q_{i})$, which proves $(\text{II}_{i})$. 
 
Finally we derive the lemma from $(\text{I}_{1})$, $(\text{II}_{1})$ and Claim~\ref{claim:order_transversals} by repeating the arguments in the inductive step. That is, by $(\text{I}_{1})$ and Claim~\ref{claim:order_transversals} there exists some $y \in P_0$ such that $x'\lleq y'$ for some $x'\in F_1(a_0, y)$, $y'\in G_1(a_0,y)$. Then by $(\text{II}_{1})$ there is some $S\in \SSS$ with $S\subseteq (a_0, y)_\lleq \cup \bigcup_{j>0} (a_j, b_j)_\lleq$ that is disjoint from any cut $(P_j, Q_j)$ for $j>0$. By construction $S$ is also disjoint from $(P_0, Q_0)$. 
\end{proof}   

Lemma~\ref{lemma:order_transversals} can be generalized to higher dimensions as follows.

\begin{proposition}
Let $\SSS\subseteq \mathcal{P}(M^n)$ be a definable family of nonempty sets. At least one of the following holds. 
\begin{enumerate}[(1)]
\item \label{itm1:small_transversal} There exists an infinite subfamily of $\SSS$ of pairwise disjoint sets.
\item \label{itm2:small_transversal} There exists a tame transversal $T$ for $\SSS$ with $\dim T <n$. 
\end{enumerate}
\end{proposition}
\begin{proof}
We assume the negation of \eqref{itm1:small_transversal} and prove \eqref{itm2:small_transversal}. We work by induction on $n$. The case $n=1$ is given by Lemma~\ref{lemma:order_transversals}. Suppose that $n>1$. Let $\pi_1$ denote the projection onto the first coordinate and consider the definable family $\pi_1(\SSS)=\{\pi_1(S) : S\in\SSS\}$. If $\SSS$ does not contain an infinite subfamily of pairwise disjoint sets then the same holds for $\pi_1(\SSS)$. By the case $n=1$ there exists a finite set $\{\xi_0, \ldots, \xi_m\}$ that is a tame transversal for $\pi_1(\SSS)$. Let $\NN=(N,<,\ldots)$ be a tame extension of $\MM$ that contains $\{\xi_0, \ldots, \xi_m\}$. Then $\cup_{0\leq i \leq m}\{\xi_i\}\times N^{n-1}$ is a tame tranversal for $\SSS$ of dimension $n-1$.  
\end{proof}

In order to prove Theorem~\ref{thm:finite_tame_transversals} we first show that any consistent definable family has a finite tame transversal. We in fact prove the following more precise statement. Recall that a family of sets $\SSS$ is $n$-consistent if any subfamily of at most $n$ sets has nonempty intersection. 

\begin{lemma}\label{lemma:FIP_transversal}
Let $\SSS\subseteq \mathcal{P}(M^n)$ be a definable family of sets. If $\SSS$ is $2^{n}$-consistent then it has a finite tame transversal. 
\end{lemma}
\begin{proof}
We proceed by induction on $n$. The case $n=1$ is given by Lemma~\ref{lemma:order_transversals}. Suppose that $n>1$. Let $\SSS=\{\varphi(u,M^n) : u\in \Omega\}$. 


Consider the family $\{\pi(S\cap S') : S, S'\in \SSS\}$. This family is definable and $2^{n-1}$-consistent. Hence by induction hypothesis there exists a finite collection of definable $(n-1)$-types $p_0, \ldots, p_m$ such that, for every $S, S'\in \SSS$, $\pi(S\cap S')$ belongs in one of them. We construct a finite tame transversal for $\SSS$ given by types whose projection is one of $p_0, \ldots p_m$. The approach is to use cell decomposition to witness the \textquotedblleft fiber over the $p_i$" of each set in $\SSS$ as a finite union of intervals on some definable preordered set, and then apply Lemma~\ref{lemma:order_transversals}.

Let $\DD$ be a uniform cell decomposition of $\SSS$, i.e. a cell decomposition of $\varphi(\Omega,M^n)$. For any $D\in \DD$ and $u\in \Omega$ with $D_u\neq \emptyset$ let $f_{D_u}$ and $g_{D_u}$ denote the functions such that $D_u=(f_{D_u}, g_{D_u})$ or $D_u=graph(f_{D_u})=graph(g_{D_u})$. 
Let $\mathcal{H}_i=\{ f_{D_u}, g_{D_u} : D\in \DD, u\in \Omega, \pi(D_u)\in p_i\}$. We briefly observe that this family is definable. 

For every $i$ and $D\in \DD$ let $\Omega(i, D)=\{u\in \Omega : \pi(D_u)\in p_i\}$. By definability of $p_i$ these sets are definable. In particular, the families $\{f_{D_u}, g_{D_u} : u\in \Omega(i,D)\}$ are definable. Note that $\mathcal{H}_i=\{f_{D_u}, g_{D_u} : D\in\DD, u\in \Omega(i,D)\}$. By re-indexing $\mathcal{H}_i$ in terms of a finite disjoint union $X_i$ of sets $\Omega(i,D)$ for $D$ ranging in $\DD$, we may take $\mathcal{H}_i$ to be a definable family $\{ h_x : x\in X_i\}$. 

Let $\lleq_i$ be the total preorder induced by $p_i$, in particular on $\mathcal{H}_i$, described in Subsection~\ref{remark:preorder}. As usual we abuse notation and also let $\lleq_i$ denote the induced definable preorder on the index set $X_i$. Let $(X,\lleq)$ be the definable totally preordered set that extends the disjoint union of $(X_i,\lleq_i)$, $0\leq i \leq m$, by letting $x\lle y$ for any $x\in X_i$, $y\in X_j$ with $i < j$. 


For every $u\in\Omega$, let $B_u$ denote the union in $(X,\lleq)$ of intervals $[y,z]_{\lleq}$ where $y,z\in X_i$ for some $0\leq i\leq m$ and satisfy that there exists $D\in \DD$ with $D_u=(h_y,h_z)$ or $D_u=graph(h_y)=graph(h_z)$. So every set $B_u$ is the union of at most $(m+1)|\DD|$ closed intervals in $(X,\lleq)$. By definability of $\SSS$, $\lleq$ and each $X_i$, the family $\BB_S=\{B_u : u\in\Omega\}$ is definable. Onwards for clarity we write $B_S$ in place of $B_u$ where $S=\varphi(u,M^n)$. This is valid because $B_u=B_v$ whenever $S=\varphi(u,M^n)=\varphi(v,M^n)$. However this is done only for clarity, since the proof can also be completed in terms of the subscript $u$. 

Since every two sets $S, S'\in \SSS$ satisfy that $\pi(S\cap S')\in p_i$ for some $0\leq i \leq m$, the family $\BB_{\SSS}$ is $2$-consistent. Consequently, by Lemma~\ref{lemma:order_transversals}, there exists finitely many definable types $q_0, \ldots, q_l$ such that, for every $S\in\SSS$, the set $B_S$ belongs in one of them. We complete the proof by fixing an arbitrary $0\leq j\leq l$ and showing that the subfamily $\SSS_j=\{S\in \SSS : B_S\in q_j\}$ admits a finite tame transversal. We will use the construction of types appearing in Definition~\ref{definition:construction_types}.

Let $i$ be such that $X_i\in q_j$. We define $Q=\{x\in X_i : (-\infty, x]_{\lleq}\in q_j\}$. Note that, by definability of $q_j$, this set is definable. We complete the proof by considering two cases. First we assume that $Q$ has a minimum $\hat{x}$ with respect to $\lleq$, and show that the family $\SSS_j$ has a tame transversal of size at most three, given by the types $h_{\hat{x}}|_{p_i}$, $h^+_{\hat{x}}|_{p_i}$ and $h^-_{\hat{x}}|_{p_i}$. Then we assume that $Q$ does not have a minimum, and show that $\SSS_j$ extends to a definable type (i.e. it has a tame transversal of size one).
We make use of the following short claim. 

\begin{claim}\label{claim:FIP_transversal}
Let $S\in\SSS_j$. Then there exists $y,z\in X_i$ with $[y,z]_{\lleq}\in q_i$, where $h_y$ and $h_z$ define a cell inside $S$, i.e. $dom(h_y)=dom(h_z)$ and $(h_y,h_z)\subseteq S$ or $graph(h_y)=graph(h_z)\subseteq S$. Moreover $y\lleq Q$ and $z\in Q$. 
\end{claim}

Let $S\in\SSS_j$. Recall that, definition of $\SSS_j$, we have that $B_S\in q_j$. By definition of $B_S$ there must exists $y,z\in X_i$ with $[y,z]_{\lleq}\in q_i$ such that $h_y$ and $h_z$ define a cell inside $S$. Clearly $z\in Q$. Moreover, for every $x\in Q$, we have that $[y,z]_{\lleq}\cap (-\infty,x]_{\lleq} \in q_i$, and so $y\lleq x$. Hence $y\lleq Q$. This completes the prove of the claim. We proceed with the proof by cases.

\textbf{Case $1$:} The set $Q$ has a minimum $\hat{x}$. Let $y,z\in X_i$ be as in Claim~\ref{claim:FIP_transversal}. Note that $y\lleq \hat{x}\lleq z$. By definition of $\lleq$ we have that 
\[
\{ a : h_y(a) \leq h_{\hat{x}}(a) \leq h_z(a)\}\in p_i.
\] 
It follows that at least one of the following three sets must belong in $p_i$
\begin{align*}
&\{ a : h_y(a) = h_{\hat{x}}(a) = h_z(a)\} \\
&\{ a : h_y(a) \leq h_{\hat{x}}(a) < h_z(a)\} \\
&\{ a : h_y(a) < h_{\hat{x}}(a) \leq h_z(a)\}.
\end{align*}
Since $h_y$ and $h_z$ define a cell inside $S$, we conclude that $S \in h_{\hat{x}}|_{p_i} \cup h^+_{\hat{x}}|_{p_i} \cup h^-_{\hat{x}}|_{p_i}$. 

\textbf{Case $2$:} $Q$ does not have a minimum. Let $q$ be the type of all sets $S\in\defn$ such that $S\in h_x|_{p_i}$ for all $x$ in some subset of $Q$ that is unbounded from below in $Q$ with respect to $\lleq$. We show that this type is consistent and definable, and that $\SSS_j \subseteq q$. This follows the arguments in Case $2$ in the proof of Proposition~\ref{prop:parameter_version_existence_downward_directed_families}, and so we are concise in the presentation.

To prove that $\SSS_j \subseteq q$ let us fix $S\in \SSS_j$. Let $y,z\in X_i$ be as in Claim~\ref{claim:FIP_transversal}. Since $Q$ does not have a minimum we have that $y\prec Q$ and moreover there exists $x\in Q$ with $x\prec z$. For every $x\in Q$ with $x\prec z$ we have that $y\prec x \prec z$. By definition of $\lleq$, from the inequality $y\prec x \prec z$ it follows that $S\in h_x|_{p_i}$. We conclude that $\SSS_j \subseteq q$.

Note that, by definability of $p_i$, $\mathcal{H}_i$ and $Q$, the type $q$ is definable. It remains to show consistency of $q$. By cell decomposition it suffices to do it for cells. 
Let $C$ be a cell in $q$. Since $Q$ does not have a minimum there must exist $x, x'\in Q$ with $x \prec x'$ such that $C\in h_{x}|_{p_i} \cap h_{x'}|_{p_i}$. But then we have that
\[
\{a :  h_{x}(a)< h_{x'}(a),\, \al a,h_x(a) \ar \in C, \, \al a,h_{x'}(a) \ar \in C\}\in p_i
\]
In particular we derive that $C$ must be a cell of the form $(f,g)$ for two functions $f$ and $g$ with domain in $p_i$. Moreover observe that, for any $x\in Q$, it holds that  $C\in h_x|_{p_i}$ if and only if 
\[
\{a : \al a, h_{x}(a)\ar \in C \} = \{a : f(a)<h_x(a)<g(a)\}\in p_i.
\]
In other words, $C\in h_x|_{p_i}$ if and only if $f\prec_i h_x \prec_i g$. Now let $C_1=(f_1,g_1), C_2=(f_2,g_2)$ be two cells in $q$. For any $k\in\{1,2\}$ observe that, since $C_k\in q$, then it must be that $f_k\prec_i \{h_x : x\in Q\}$ and $h_{x_k}\prec_i g_k$ for some $x_k\in Q$. Let $x'=\min_{\lleq}\{x_1, x_2\}$.
For every $x\in Q$ with $x\lleq x'$ and any $k\in\{1,2\}$ we have that $f_k \prec_i h_{x} \prec_i g_k$, and so $C_k \in h_{x}|_{p_i}$. It follows that $C_1\cap C_2\in q$. This completes the proof of the lemma. 
\end{proof}

We may now prove Theorem~\ref{thm:finite_tame_transversals}. Recall that a family of sets is $n$-inconsistent if every subfamily of size $n$ has empty intersection. Recall also that, for a definable set $X$, we denote its boundary (in the o-minimal euclidean topology) by $bd(X)$. If $X$ is contained in a set $Y$ then let $bd_Y(X)=bd(X)\cap Y$.

\begin{proof}[Proof of Theorem~\ref{thm:finite_tame_transversals}]
To prove that~\eqref{itm:thm_transversals_2} implies~\eqref{itm:thm_transversals_1} note that, if $\SSS$ admits a tame transversal of cardinality $l$, then it has a covering of $l$ consistent subfamilies (those given by the sets containing a given element in the transversal). It follows that, for every $n>0$, $\SSS$ has the $(nl+1,n+1)$-property. We assume the negation of~\eqref{itm:thm_transversals_2} and derive the negation of~\eqref{itm:thm_transversals_1}. Note that, if $\SSS$ does not have a finite tame transversal, then $\SSS$ must be infinite (otherwise it suffices to take a point from every set in $\SSS$ to build a finite transversal in $M^k$). It follows that we may assume throughout that $\dim \cup \SSS\geq 1$, and in particular that $n>1$. Observe that the negation of~\eqref{itm:thm_transversals_1} is the statement that, for every $m\geq n$, there exists a subfamily of $\SSS$ of cardinality $m$ that is $n$-inconsistent. Instead of this we will show the following apparently stronger statement.

Let $X_1 \subseteq \cdots \subseteq X_{n}$ be a collection of nested definable subsets of $M^k$ satisfying that $\dim X_i < i$ for every $1\leq i < n$ and $\cup\SSS\subseteq X_{n}$. For every $m\geq n$, there exists $\FF\subseteq \SSS$ with $|\FF|=m$ such that, for every $1\leq i \leq n$ and subfamily $\FF'\subseteq \FF$ with $|\FF'|=i$, it holds that $\cap \FF' \cap X_i=\emptyset$. That is, for every $1\leq i \leq n$, the family $\FF\cap X_i$ is $i$-inconsistent. In particular, since $\cup \SSS\subseteq X_{n}$, the family $\FF$ is $n$-inconsistent.


Let $\prop{m}{n}$ denote the statement that the above holds for some fixed $m \geq n > 1$ and any given definable family $\SSS \subseteq \mathcal{P}(M^k)$ of nonempty sets with $n > \dim \cup \SSS$ and without a finite tame transversal. We prove $\prop{m}{n}$ by induction. We let $\prop{2}{2}$ be the base case. Suppose that $m>2$. If $m>n$ then we use $\prop{m-1}{n}$ to derive $\prop{m}{n}$. Otherwise we use $\prop{m-1}{n-1}$ to derive $\prop{m}{n}$. We prove all cases simultaneously. We divide the argument into two parts. 

Let us fix $m \geq n > 1$ and $\SSS\subseteq \mathcal{P}(M^k)$ a definable family of nonempty sets with $n > \dim \cup\SSS$ that does not admit a finite tame transversal; and $X_1\subseteq \cdots \subseteq X_{n}$ a family of definable nested subsets of $M^k$ with $\dim X_i < i$ for $1\leq i < n$ and $\cup\SSS\subseteq X_{n}$. The first part of the argument consists of the following claim. 

\begin{claim}\label{claim:trans_theorem}
Let $\SSS'\subseteq \SSS$ be a finite subfamily and, for every $1\leq i < n$, let $Y_{i}=\cup_{S\in \SSS'} bd_{X_{i+1}} (S \cap X_{i+1})$. There exists a subfamily $\FF\subseteq \SSS$ of cardinality $m-1$ that is $n$-inconsistent and such that $\FF \cap (X_i \cup Y_{i})$ is $i$-inconsistent for $1 \leq i < n$. 
\end{claim}

Note that, by o-minimality, $\dim Y_i < i$ for every $i$. So $\dim (X_i \cup Y_{i})<i$ for every $1\leq i < n$. Note that, if $m=n$, then any subfamily $\FF\subseteq \SSS$ of cardinality $m-1=n-1$ is $n$-inconsistent vacuously. It follows that, if $m=n=2$, the claim follows easily from the fact that $\SSS$ does not admit a finite transversal in $M^k$.

Suppose that $m>2$. If $m>n$ then the claim is given by $\prop{m-1}{n}$ applied to $\SSS$ and the collection of nested sets $\{X_1\cup Y_{1},\ldots, X_{n-1}\cup Y_{n-1}, X_{n}\}$. Suppose that $m=n>2$. We may apply $\prop{m-1}{n-1}$ to the definable family $\SSS \cap (X_{n-1}\cup Y_{n-1})$ and collection of nested sets $\{X_1\cup Y_1, \ldots, X_{n-1}\cup Y_{n-1}\}$. It follows that there exists $\FF\subseteq \SSS$ with cardinality $m-1=n-1$ such that $\FF \cap (X_{n-1}\cup Y_{n-1})\cap (X_i \cup Y_{i})=\FF\cap (X_i \cup Y_{i})$ is $i$-inconsistent for $1\leq i < n$. Since $|\FF|=n-1$ the family $\FF$ is $n$-inconsistent vacuously. This completes the proof of the claim. 

\sloppy
We continue with the second part of the proof. For any $S\in \SSS$ consider the family of all ordered tuples $\al S_1, \ldots, S_{m-1}\ar \in \prod_{1\leq i <m} \SSS$ such that $\{ S_1,\ldots, S_{m-1}\}$ is $n$-inconsistent and $\{S_1,\ldots, S_{m-1}\} \cap (X_{i} \cup bd_{X_{i+1}}(S\cap X_{i+1}))$ is $i$-inconsistent for every $1\leq i < n$. We denote this family $\SSS^{m}(S)$. 
Note that, by Claim~\ref{claim:trans_theorem}, the family $\{ \SSS^{m}(S) : S\in \SSS\}$ is consistent.

We now explicitly identify the family $\{ \SSS^{m}(S) : S\in \SSS\}$ with a family of definable sets as follows. Let $\varphi(u,v)$ be such that $\SSS=\{\varphi(u,M^k) : u\in \Omega\}$. For any $u\in \Omega$, $S=\varphi(u,M^k)$, let $\Omega^m(u)$ denote the set of all index tuples $\al u_1,\ldots, u_{m-1}\ar$ such that $\al\varphi(u_1,M^k), \ldots, \varphi(u_{m-1}, M^k)\ar \in \SSS^{m}(S)$. Clearly the family $\{ \Omega^m(u) : u\in \Omega\}$ is definable. Like $\{ \SSS^{m}(S) : S\in \SSS\}$, it is consistent.

By Lemma~\ref{lemma:FIP_transversal} it follows that $\{ \Omega^m(u) : u\in \Omega\}$ admits a finite tame transversal. That is, we may partition $\Omega$ into finitely many definable subfamilies $\Omega_1,\ldots, \Omega_s$ satisfying that, for each $1\leq i \leq s$, the family $\{\Omega^m(u) : u\in \Omega_i\}$ extends to a definable type. Note that, since by assumption $\SSS$ does not have a finite tame transversal, there must be some $1\leq i\leq s$ such that $\SSS_i=\{\varphi(u,M^k) : u\in \Omega_i\}$ does not have a finite tame transversal either. Hence, by passing if necessary to a subfamily $\SSS_i$ of $\SSS$, we may assume that $\{ \Omega^m(u) : u\in \Omega\}$ extends to a definable type $p$.

Let $\NN=(N,\ldots)$ be a tame extension of $\MM$ that realizes $p$. Such an extension exists by one direction of the Marker-Steinhorn theorem \cite[Theorem 2.1]{mark_stein_94}. Onwards let $S^*=S(\NN)$ for every $S\in\SSS$ and $X_i^*=X_i(\NN)$ for every $1\leq i \leq n$. That is, we use an asterisk to denote the interpretation in $\NN$ of a definable set in $\MM$. Let $c=\al c_1,\ldots, c_{m-1}\ar$ be a realization in $\NN$ of $p$ and, for $1\leq i <m$, let $S_i=\varphi(c_i, N^k)$. Then we have that $\{ S_1,\ldots, S_{m-1}\}$ is $n$-inconsistent and $\{ S_1,\ldots, S_{m-1}\} \cap (X^*_{i} \cup bd_{X^*_{i+1}}(S^* \cap X^*_{i+1}))$ is $i$-inconsistent for every $S\in \SSS$ and $1\leq i < n$. We prove that there exists $S\in \SSS$ such that the family $\{ S_1,\ldots, S_{m-1}, S^*\}$ witnesses $\prop{m}{n}$ in $\NN$, by satisfying that $\{ S_1,\ldots, S_{m-1}, S^*\}\cap X^*_i$ is $i$-inconsistent for every $1\leq i \leq n$. Since $\NN$ is an elementary extension of $\MM$, then an analogous family witnessing $\prop{m}{n}$ in $\MM$ must exist too.

Let $\CC$ be a cell decomposition of $N^k$ compatible with each intersection of sets in $\{S_1,\ldots, S_{m-1}, X^*_1, \ldots, X^*_{n}\}$. For each cell $C\in \CC$ let us choose a point $\xi_C\in C$. Since $\SSS$ does not have a finite tame transversal, there exists $S\in\SSS$ such that $S^*\cap \{\xi_C :C\in \CC\}=\emptyset$. Since $X_1$ is finite (because $\dim X_1 < 1$), we also choose $S$ so that $S\cap X_1=\emptyset$. Let $\FF=\{ S_1,\ldots, S_{m-1}, S^*\}$. We show that $\FF \cap X^*_{i}$ is $i$-inconsistent for every $1\leq i \leq n$.

Towards a contradiction suppose that there exists some $1 \leq i \leq n$ and a subset $\FF'\subseteq \FF$ of cardinality $i$ such that $\cap \FF' \cap X^*_{i} \neq \emptyset$. Note that, by construction of $\{ S_1, \ldots, S_{m-1} \}$, it must be that $S^*\in\FF'$. If $i=1$ then $\FF'=\{S^*\}$ and $S^* \cap X^*_{1} \neq \emptyset$ contradicts that $S\cap X_1=\emptyset$. Suppose that $i>1$. Let $\FF''=\FF'\setminus \{S^*\}$ and let us fix a cell $C\in \CC$ such that $C\cap (\cap\FF') \cap X^*_{i} \neq \emptyset$. Since $\CC$ is compatible with $\cap \FF''\cap X^*_{i}$ it must be that $C\subseteq \cap\FF'' \cap X^*_{i}$. Moreover recall that, by construction of $\{ S_1, \ldots, S_{m-1} \}$, the intersection $\cap \FF''$ (and in particular $C$) is disjoint from $bd_{X^*_i}(S^*\cap X^*_i)$. By definable connectedness from $C\cap S^*\neq \emptyset$ it follows that $C\subseteq S^*$. This however contradicts that $\xi_C \notin S^*$. This completes the proof of the theorem.
\end{proof}



Note that, if we add to the statement of Theorem~\ref{thm:finite_tame_transversals} the condition that all the sets in $\SSS$ are finite, then the transversal can always be assumed to be in $M^k$.

With the use of the Alon-Kleitman-Matou\v{s}ek $(p,q)$-theorem (Theorem~\ref{thm:vc_tranversal}), Lemma~\ref{lemma:FIP_transversal} can be used to prove a different version of Theorem~\ref{thm:finite_tame_transversals}, where the $n$ in the statement is substituted by any interger greater than the $\VC$-codensity of $\SSS$. 

\begin{corollary}\label{prop:vc_tame_transversals}
Let $\SSS$ be a definable family of sets with the $(m,n)$-property, where $m\geq n$ and $n$ is greater than the $\VC$-codensity of $\SSS$. Then $\SSS$ admits a finite tame transversal. 
\end{corollary}
\begin{proof}
For simplicity we assume that $\SSS \subseteq \mathcal{P}(M)$, the proof otherwise being analogous. Applying Theorem~\ref{thm:vc_tranversal}, there exists a natural number $l>0$ such that every finite $\SSS'\subseteq \SSS$ has a transversal of cardinality at most $l$. 

Let $\FF=\{ F(S)\subseteq M^{l} : S\in \SSS\}$ be the definable family of sets given by $F(S)=(S\times M^{l-1})\cup (M\times S\times M^{l-2})\cup\cdots\cup (M^{l-1} \times S)$. We observe that $\FF$ is consistent. For this it suffices to note that, if $\SSS'$ is a finite subfamily of $\SSS$, and $\{x_1,\ldots, x_l\}\subseteq M$ is a transversal for $\SSS'$, then $\al x_1,\ldots, x_l \ar\in M^{l}$ is in every $F(S)$ for $S\in\SSS'$. 

So, applying Lemma~\ref{lemma:FIP_transversal}, $\FF$ has a finite tame transversal $T\subseteq U^{l}$. Let $H=\cup_{1\leq i \leq l} \pi_i(T)$, where $\pi_i$ denotes the projection to the $i$-th coordinate. We claim that $H$ is a tame tranversal for $\SSS$. This follows from the observation that, for every $S\in \SSS$ and $\al x_1,\ldots, x_l\ar \in F(S)$, there is some $i$ such that $x_i\in S$.  
\end{proof}

Recall that Theorem~\ref{cor:vc_density} states that any uniform family $\SSS\subseteq Def(M^k)$ has $\VC$-codensity at most $k$. One may use this to recover Theorem~\ref{thm:finite_tame_transversals} from Corollary~\ref{prop:vc_tame_transversals}.


\begin{remark}\label{remark:existence_ddf_from_FTT}
Recall that Lemma~\ref{lemma:existence_complete_directed_families}, and specifically the fact that every definable downward directed family extends to a definable type, is used in proving Theorem~\ref{thm:finite_tame_transversals}, by playing a part in the proof of Lemma~\ref{lemma:order_transversals}. It is worth noting that the fact that every definable downward directed family extends to a definable type can be recovered from Theorem~\ref{thm:finite_tame_transversals} (or from Lemma~\ref{lemma:FIP_transversal}). This follows from the observation that, if a downward directed family $\SSS$ has a finite tame transversal, given by types $p_1, \ldots, p_n$, then there must exist $1\leq i \leq n$ such that $\SSS\subseteq p_i$. Otherwise, for every $i$ there is $S(i)\in\SSS$ with $S(i)\notin p_i$. But then it follows that any $S\subseteq S(1)\cap \cdots \cap S(n)$ is not contained in $p_i$ for any $i$, contradicting that the types $p_i$ witness a finite tame transversal of $\SSS$. 
\end{remark}

The following example shows that the bounds in Theorem~\ref{thm:finite_tame_transversals} and Corollary~\ref{prop:vc_tame_transversals} cannot be improved. 
\begin{example}
For any $u,v\in M$ consider the \textquotedblleft cross" $S(u,v)=(\{u\}\times M)\cup (M\times \{v\})$. The family of crosses $\SSS=\{S(u,v) : u,v\in M\}$ is a definable family of subsets of $M^2$ with $\vc$-codensity $2$. To prove the latter we leave it to the reader to check that ${n\choose 2} + n \leq \pi_{\SSS}^*(n)\leq (n+1)^2$. Moreover $\SSS$ is $2$-consistent, since $\{\al u_1, v_2\ar, \al u_2, v_1\ar\}\subseteq S(u_1,v_1)\cap S(u_2,v_2)$ for every $u_1,v_1,u_2,v_2\in M$. We observe that $\SSS$ does not have a finite transversal in any elementary extension.  

Let $\MM\lleq\NN=(N,\ldots)$ and $X\subseteq N^2$ be a finite set. Let $X'\subseteq M$ be the set of coordinates of points in $X$. Pick any $u, v \notin X'\cap M$. Then clearly the interpretation of $S(u,v)$ in $\NN$ is disjoint from $X$. 
\end{example}

Recall that Lemma~\ref{lemma:order_transversals} shows that a definable family of subsets of $M$ has a finite tame transversal if and only if it has the $(\omega,2)$-property. We ask whether Theorem~\ref{thm:finite_tame_transversals} can be improved to generalize this to higher dimensions. 
\begin{question}\label{question:generalizing_FTT_theorem}
Let $\SSS$ be a definable family of nonempty sets and let $n>\dim\cup\SSS$. If $\SSS$ has the $(\omega, n)$-property then does it have a finite tame transversal?
\end{question}
Note that the answer to the above question is positive whenever $\MM$ is $\omega_1$-saturated, since any definable family of sets with the $(\omega, n)$-property is going to have that $(m,n)$-property for some large enough $m\geq n$. 

We now present a proposition that follows from the ideas in the proof of Lemma~\ref{lemma:order_transversals}. 

\begin{proposition}
Let $\SSS\subseteq \mathcal{P}(M)$ and $k$ be such that every $S\in \SSS$ is union of at most $k$ intervals and points. Suppose that $\SSS$ is $(k+1)$-consistent. Then at least one holds:
\begin{enumerate}[(i)]
\item $\SSS$ has a finite transversal in $M$,
\item $\SSS$ extends to a definable type. 
\end{enumerate}
\end{proposition}
\begin{proof}
Suppose that $\SSS$ as in the proposition does not have a finite transversal in $M$ and also does not extend to a definable type. We prove that $\SSS$ is not $(k+1)$-consistent. Since the construction and arguments that follow are similar to those in the proof of Lemma~\ref{lemma:order_transversals} we are concise in the presentation.

If $\SSS$ does not have a finite transversal in $M$, then for every finite $X\subseteq M$, there is $S\in\SSS$ with $S\cap X=\emptyset$. By first order logic compactness let $\badS \in \SSS(\Mon)$ be such that $\badS\cap M=\emptyset$. Since $\badS$ is union of at most $k$ points and intervals, there are at most $k$ cuts $(P_1,Q_1), \ldots, (P_m,Q_m)$ in $M$ such that every point in $\badS$ realizes one of these cuts. We show that, for every $1\leq i \leq m$, there exists some $S_i\in \SSS$ that is disjoint from $(P_i,Q_i)$. It follows that the family $\badS, S_1(\Mon),\ldots, S_m(\Mon)$, which has size at most $k+1$, has empty intersection. Hence $\SSS$ is not $(k+1)$-consistent.  

Let us fix $1\leq i \leq m$. If $(P_i,Q_i)$ is definable then, since by assumption $\SSS$ does not extend to a definable type, there exists some $S_i\in\SSS$ that is disjoint from the cut. Suppose that $(P_i,Q_i)$ is not definable. Let us fix $a\in P_i$ and $b\in Q_i$ with $\{a,b\}\subset P_j$ or $\{a,b\}\subset Q_j$ for every $j\neq i$. Note that, since $(P_i,Q_i)$ is not definable, we have that $(a,+\infty)\cap P_i\neq \emptyset$. Consider the definable set $F=\{ t\in M : t>a \text{ and there exists }  S\in \SSS \text{ with } (a,t)\cap S=\emptyset\}$. Then $\badS$ witnesses the fact that $(a,+\infty)\cap P_i \subseteq F$. Since $(P_i,Q_i)$ is not definable we have that $F\cap Q_i\neq \emptyset$ (see Claim~\ref{claim:order_transversals} in the proof of Lemma~\ref{lemma:order_transversals}). By definition of $F$ we conclude that there exists some $S_i\in \SSS$ disjoint from $(P_i,Q_i)$.    
\end{proof}

It is possible that the above proposition can be generalized to higher dimensions, albeit that author has not been able to prove a precise statement for said generalization.

We end with a proposition adapting a result on transversals from finite and compact combinatorics (Proposition~\ref{prop:Mirsky}) to definable families in o-minimal structures. Since the result is not central to this paper, we will, like in the previous proposition, be concise in the proof. 

\begin{proposition}
Let $\SSS\subseteq \mathcal{P}(M)$ be a definable family of intervals. Let $k\geq 1$ be the maximum such that there exists $k$ pairwise disjoint sets in $\SSS$. Then $\SSS$ has a tame transversal of size $k$.
\end{proposition}
\begin{proof}
We proceed by induction on $k$. 

For the base case $k=1$ let $H$ be the definable set of all left endpoints of intervals in $\SSS$. Let $a=\sup H$. Since $k=1$, meaning that $\SSS$ is $2$-consistent, note that, by definition of $H$, every right endpoint of an interval in $\SSS$ must be greater or equal to $a$. Now suppose that $\SSS$ does not extend to the definable type with basis $\{(t,a) : t<a\}$. Then there must exist $S\in \SSS$ with $S\subseteq [a,+\infty)$. Additionally  suppose that $\SSS$ does not extend to the definable type with basis $\{(a,t) : a<t\}$. Then there must exists $S'\in \SSS$ with $S'\subseteq (-\infty,a]$. Finally, suppose that $\SSS$ does also not extend to $\tp(a/M)$, and let $S''\in \SSS$ be such that $a\notin S''$. If $S''\subseteq (-\infty,a)$ then $S''\cap S=\emptyset$ and if $S''\subseteq (a,+\infty)$ then $S''\cap S'=\emptyset$, contradicting that $k=1$.  

Now let $k>1$. Let $\SSS_{0}\subseteq \SSS$ be the subfamily of leftmost intervals in any pairwise disjoint subfamily of $\SSS$ of size $k$. Note that, but induction hypothesis, the definable family $\SSS\setminus \SSS_0$ has a tame transversal of size $k-1$. We complete the proof by noting that $\SSS_0$ is $2$-consistent, and so, by the base case, extends to a definable type. 

Suppose towards a contradiction that there are $S, S'\in \SSS_0$ with $S\cap S'=\emptyset$. Without loss of generality let $S<S'$. Let $\FF\subseteq \SSS$ be a subfamily of size $k$ of pairwise disjoint sets that includes $S'$ as leftmost interval. Then the family $\{S\}\cup\FF$ is a family of size $k+1$ of pairwise disjoint sets. Contradiction.   
\end{proof}

\subsection{Forking, dividing and definable types}\label{section:forking_dividing}


In this subsection we reformulate Theorem~\ref{thm:finite_tame_transversals} as a statement about the relationship between dividing and definable types known in o-minimal and some more general NIP theories. This is the subject of ongoing research among NIP theories~\cite{simon15}. Until the end of the subsection we drop the assumption that $\MM$ and $\Mon$ are o-minimal. Unless stated otherwise types are global over $U$.  

Recall that formula $\varphi(x,b)$ \emph{$n$-divides over $A$} (a small set), for some $n\geq 1$, if there exists a sequence of elements $(b_i)_{i<\omega}$ in $U^{l(b)}$, with $\tp(b_i/A)=\tp(b/A)$, such that $\{ \varphi(x,b_i) : i<\omega\}$ is $n$-inconsistent. Equivalently, if the family of sets of the form $\varphi(U^{l(x)}, b')$, where $\tp(b'/A)=\tp(b/A)$, does not have the $(\omega, n)$-property. A formula $\varphi(x,b)$ \emph{divides} if it $n$-divides for some $n$. Conversely, note that a formula $\varphi(x,b)$ does not divide over $A$ if the family $\{ \varphi(U^{l(x)},b') : \tp(b'/A)=\tp(b/A)\}$ has the $(\omega, n)$-property for every $n$. Hence, not dividing is an intersection property. 

A formula \emph{forks over $A$} if it implies a finite disjunction of formulas that divide each over $A$. In $NTP_2$ theories (a class which includes NIP and simple theories) forking and dividing over a model are equivalent notions (see Theorem $1.1$ in~\cite{cher_kap_12}).

The next equivalence was proved first for o-minimal expansions of ordered fields\footnote{Dolich actually works with \textquotedblleft nice" o-minimal theories, a class of structures which includes o-minimal expansions of ordered fields.} by Dolich~\cite{dolich04} (where he considers forking over any set and not just models) and for unpackable $\VC$-minimal theories, a class which includes o-minimal theories, by Cotter and Starchenko~\cite{cotter_star_12}. The following is the best generalization up to date, due to Simon and Starchenko~\cite{simon_star_14}. 

\begin{theorem}\label{thm:forking_definable_types}[\cite{simon_star_14}, Theorem 5] Let $T$ be a dp-minimal $\LL$-theory with monster model $\Mon$ satisfying that, for every $A\subseteq U$, every nonempty unary $A$-definable set extends to an $A$-definable type in $S_1(A)$. Let $M\models T$ and $\varphi(x,b)\in\LL(U)$. The following are equivalent. 
\begin{enumerate}[(i)]
\item\label{itm:thm_forking_definable_types_1} $\varphi(x,b)$ does not fork over $M$. 
\item\label{itm:thm_forking_definable_types_2} $\varphi(x,b)$ extends to an $M$-definable type. 
\end{enumerate}
\end{theorem}

For precise definitions of unpackable $\VC$-minimal and dp-minimal theory see~\cite{cotter_star_12} and~\cite{simon_star_14} respectively. Note that, by the equivalence between forking and dividing in $NTP_2$ theories, \eqref{itm:thm_forking_definable_types_1} and~\eqref{itm:thm_forking_definable_types_2} above are also equivalent to $\varphi(x,b)$ not dividing over $M$. 

We describe how our Theorem~\ref{thm:finite_tame_transversals} implies an improvement of the o-minimal part of Theorem~\ref{thm:forking_definable_types}. In particular, Theorem~\ref{thm:finite_tame_transversals} can be restated as the equivalence between a formula not $n$-dividing over a small set $A$ for some large enough $n$ and said formula extending to an $A$-definable type. 


Consider a weaker form of Theorem~\ref{thm:finite_tame_transversals}, given by substituting the lower bound $\dim \cup\SSS$ in the statement of the theorem by $k$ where $\SSS \subseteq \mathcal{P}(M^k)$. It is easy to see that this result is equivalent to Theorem~\ref{thm:FTT_forking} below. The key tool is noticing that, by a first-order logic compactness argument, a formula $\varphi(x,b)$ does not $n$-divide over a small set $A$ if and only if there exists $\Omega\in\tp(b/A)$ and some $m\geq n$ such that $\{ \varphi(U^{l(x)},u) : u\in \Omega\}$ has the $(m,n)$-property. 

\begin{theorem}\label{thm:FTT_forking}
Let $\Mon$ be o-minimal. Let $\varphi(x,b)$ be a formula, $n>l(x)$, and let $A\subseteq U$ be a small set. The following are equivalent.   
\begin{enumerate}[(1)]
\item \label{itm:thm_transversals_1} $\varphi(x,b)$ does not $n$-divide over $A$. 
\item \label{itm:thm_transversals_2} $\varphi(x,b)$ extends to an $A$-definable type. 
\end{enumerate}  
\end{theorem}

Note that ``not $n$-dividing" for a given $n$ is in general a strictly weaker property that ``not dividing". On the other hand, if a formula extends to an $A$-definable type then clearly it does not divide over $A$. Hence~\eqref{itm:thm_transversals_1} and~\eqref{itm:thm_transversals_2} above are also equivalent to $\varphi(x,b)$ not dividing over $A$. Moreover Cotter-Starchenko~\cite[Corollary 5.6]{cotter_star_12} showed that, in unpackable $\VC$-minimal theories (in particular o-minimal theories), dividing and forking over sets (not just models) are equivalent notions, so~\eqref{itm:thm_transversals_1} and~\eqref{itm:thm_transversals_2} above are also equivalent to $\varphi(x,b)$ not forking over $A$. Hence Theorem~\ref{thm:FTT_forking} implies a parameter version of the o-minimal part of Theorem~\ref{thm:forking_definable_types}, and moreover improves the result further by considering the condition of ``not $n$-dividing", for a given large enough $n$, in place of just ``not dividing" (equivalently ``not forking").  

We end the subsection by pointing out that the equivalence between $n$-dividing and dividing that follows from Theorem~\ref{thm:FTT_forking} can also be derived using VC theory. That is, the Alon-Kleitman-Matou\v{s}ek $(p,q)$-theorem (Theorem~\ref{thm:vc_tranversal}) implies the equivalence between a formula $\varphi(x,b)$ dividing over $A$ and $n$-dividing over $A$, where $n$ is any integer greater than the $\VC$-codensity of $\{\varphi(U^{l(x)}:b') : \tp(b'/A)=\tp(b/A)\}$. By Theorem~\ref{cor:vc_density}, in the o-minimal setting these two properties are also equivalent to having that $\varphi(x,b)$ $n$-divides over $A$ for any $n>l(x)$.

\section{Definable compactness}\label{section:compactness}

In this section we use results from Sections~\ref{section:types} and~\ref{section:transversals} to characterize a notion of definable compactness among definable topologies in o-minimal structures. In our main result (Theorem~\ref{thm:compactness}) we prove the equivalence of various known definitions.  

\begin{definition}
A \emph{definable topological space} $(X,\tau)$, $X\subseteq M^n$,  is a topological space such that there exists a basis for $\tau$ that is definable. 

A \emph{definable curve in $X$} is a definable map $\gamma:(a,b)\rightarrow X$. We say that it \emph{$\tau$-converges} to $x\in X$ as $t\rightarrow a$, and write $\tau$-$\lim_{t\rightarrow a} \gamma=x$, if, for every $\tau$-neighborhood $A$ of $x$, there exists $a<t(A)<b$ such that $\gamma(s)\in A$ whenever $a<s \leq t(A)$. $\tau$-Convergence as $t\rightarrow b$ is defined analogously. We say that $\gamma$ is \emph{$\tau$-completable} if both $\tau$-$\lim_{t\rightarrow a} \gamma$ and $\tau$-$\lim_{t\rightarrow b} \gamma$ exist.
\end{definition}  

For examples of definable topologies in o-minimal structures besides the canonical \textquotedblleft euclidean" topology, by which we mean the order topology on $M$ and induced product topology on $M^n$, see the \emph{definable (manifold) spaces} studied by van den Dries~\cite[Chapter 10]{dries98} and the definable metric spaces of Walsberg~\cite{walsberg15}. For a treatment in general model theory see the work of Pillay~\cite{pillay87}.

\begin{definition}
Let $(X,\tau)$ be a definable topological space. Then $(X,\tau)$ is: 
\begin{enumerate}[(1)]
\item \emph{definably compact} if every definable downward directed family of $\tau$-closed subsets of $X$ has nonempty intersection,
\item \emph{curve-compact} if every definable curve in $X$ is $\tau$-completable. 
\end{enumerate}
\end{definition}  

The above notion of curve-compactnes is adapted from a definition for definable compactness introduced by Peterzil and Steinhorn~\cite{pet_stein_99} for definable manifold spaces, and used in more general settings by Thomas~\cite{thomas12} and Walsberg~\cite{walsberg15}. In establishing a parallelism with general topology, we rename this property curve-compactness, and reserve the former name for a property in terms of closed sets that has been explored in recent years by Johnson~\cite{johnson14}, Fornasiero~\cite{fornasiero}, and Thomas, Walsberg and the author~\cite{atw1}.


\begin{definition}

Let $(X,\tau)$ be a definable topological space and $p$ be a (possibly partial) type such that $X\in p$. We say that $x\in X$ is a \emph{limit\footnote{Fornasiero~\cite{fornasiero}, as well as Thomas, Walsberg and the Author~\cite{atw1}, use the word ``specialization" (borrowed from real algebraic geometry) to refer to limits of types. Here we use instead the terminology from Hrushovski and Loeser~\cite[Chapter 4]{hruloeser16}.} of $p$} (with respect to $\tau$) if $x$ is contained in every $\tau$-closed set in $p$. 

We say that $(X,\tau)$ is \emph{type-compact} if, for every definable type $p\in S_X(M)$, there exists $x\in X$ that is a limit of $p$.
\end{definition}

More appropriate names for the notions of ``curve-compactness" and ``type-compactness" would be ``definable curve-compactness" and ``definable type-compactness" respectively. Nevertheless, we decided to avoid here the use of the word ``definable", which is already gruelingly prevalent in this paper.

The equivalence \eqref{itm:specialization_compactness_i}$\Leftrightarrow$\eqref{itm:specialization_compactness_ii} in Proposition~\ref{prop:pre_specialization_compactness} below gives a characterization of type-compactness that holds in any model theoretic structure. In particular it does not require o-minimality.  

\begin{proposition}\label{prop:pre_specialization_compactness}
Let $(X,\tau)$ be a definable topological space. The following are equivalent. 
\begin{enumerate}[(1)]
\item \label{itm:specialization_compactness_i} $(X,\tau)$ is type-compact.
\item \label{itm:specialization_compactness_ii} Any definable family of $\tau$-closed sets that extends to a definable type has nonempty intersection.
\setcounter{specialization_compactness}{\value{enumi}}
\end{enumerate}
If $\MM$ expands $(\mathbb{R},<)$, then \eqref{itm:specialization_compactness_i} and \eqref{itm:specialization_compactness_ii} are also equivalent to:
\begin{enumerate}[(1)]
\setcounter{enumi}{\value{specialization_compactness}}
\item $(X,\tau)$ is compact. 
\end{enumerate}
\end{proposition}
\begin{proof}
To prove \eqref{itm:specialization_compactness_i}$\Rightarrow$\eqref{itm:specialization_compactness_ii}, suppose that $(X,\tau)$ is type-compact and let $\CC$ be a definable family of $\tau$-closed sets that extends to a definable type $p$. Let $x\in X$ be a limit of $p$. Then clearly $x\in \cap\CC$. 

The key element to the rest of the proof is noticing that any closed set in a topological space is an arbitrary intersection of basic closed sets.  

Suppose that \eqref{itm:specialization_compactness_ii} holds. To prove that $(X,\tau)$ is type-compact let $p\in S_X(M)$ be a definable type. Let $\BB$ denote a definable basis for the topology $\tau$. Now let $\CC$ denote the definable family of basic $\tau$-closed sets in $p$, i.e. the family of sets $C\in p$ of the form $X\setminus A$ for some $A\in \BB$. By \eqref{itm:specialization_compactness_ii}, let $x\in X$ be such that $x\in \cap \CC$. Then $x$ is a limit of $p$.  

Now suppose that $\MM$ expands $(\mathbb{R},<)$. Clearly if $(X,\tau)$ is compact then it is type-compact. Conversely, suppose that $(X,\tau)$ is type-compact and let $\CC$ be a consistent family of $\tau$-closed sets. The intersection $\cap\CC$ can be re-written as an intersection of basic closed sets. In particular we may assume that $\CC$ contains only definable sets. Now, by the Marker-Steinhorn Theorem (Theorem~\ref{thm:marker_steinhorn}), every type over $M$ is definable. Consequently $\CC$ extends to a definable type $p$. Let $x$ be a limit of $p$, then $x\in \cap \CC$. So $(X,\tau)$ is compact. 
\end{proof}

\begin{remark}
Note that the fact that type-compactness implies classical topological compactness shown in Proposition~\ref{prop:pre_specialization_compactness} holds in any structure satisfying that all types are definable. In fact, if $\varphi(u,v)$ defines a basis $\BB$ for the topology, i.e. $\BB=\{\varphi(u, M^{l(v)}) : u\in M^{l(u)}\}$, then it is enough to have that every $\varphi$-type is definable, where by $\varphi$-type (with object variable $v$) we mean any maximal consistent family of boolean combinations of sets in $\BB$. In particular this holds whenever $\varphi(u,v)$ is stable. On the other hand in~\cite[Proposition 1.2]{pillay87} it was noted that, in a stable structure, any infinite $T_1$ definable topological space must be discrete, and consequently the only type-compact $T_1$ spaces are the finite ones.   
\end{remark}

We now present the main result of this section. 

\begin{theorem} \label{thm:compactness} Let $(X,\tau)$ be a definable topological space. The following are equivalent. 
\begin{enumerate}[(1)]
\item \label{itm:compactness_1} $(X,\tau)$ is definably compact. 
\item \label{itm:compactness_2} $(X,\tau)$ is type-compact.
\item \label{itm:compactness_2.5} Any definable family of $\tau$-closed sets that extends to a definable type has nonempty intersection.
\item \label{itm:compactness_3} Any consistent definable family of $\tau$-closed sets has a finite transversal in $X$. 
\item \label{itm:compactness_5} Any definable family $\CC$ of nonempty $\tau$-closed sets with the $(m,n)$-property, where $m\geq n >\dim \cup\CC$, has a finite transversal in $X$.
\item \label{itm:compactness_4} Any definable family $\CC$ of nonempty $\tau$-closed sets with the $(m,n)$-property, where $m\geq n$ and $n$ is greater than the $\VC$-codensity of $\CC$, has a finite transversal in $X$. 
\setcounter{thm_compactness}{\value{enumi}}
\end{enumerate}
Moreover all the above imply and, if $\tau$ is Hausdorff or $\MM$ has definable choice, are equivalent to:
\begin{enumerate}[(1)]
\setcounter{enumi}{\value{thm_compactness}}
\item \label{itm:compactness_6} $(X,\tau)$ is curve-compact. 
\end{enumerate}
\end{theorem}

\begin{remark}\label{remark:peterzil_pillay_aschenbrenner}
~
\begin{enumerate}[(i)]
\item \label{itm:remark:peterzil_pillay_aschenbrenner_1} In \cite[Theorem 2.1]{pet_pillay_07} Peterzil and Pillay extracted from~\cite{dolich04} the following. Suppose that $\MM$ has definable choice (e.g. expands an ordered group). Let $\Mon$ be a monster model and $\phi(x,b)$ be a formula in $\LL(U)$ such that $\phi(U^{l(x)}, b)$ is closed and bounded (in the euclidean topology). If the family $\{ \phi(U^{l(x)},b') : \tp(b'/M)=\tp(b/M)\}$ is consistent, then $\phi(U^{l(x)},b)$ has a point in $M^{l(x)}$. They derive from this that any definable family of closed and bounded (equivalently definably compact with respect to the euclidean topology) sets that is consistent has a finite transversal \cite[Corollary 2.2 $(i)$]{pet_pillay_07}. Our work generalizes these results in a number of ways: we drop the assumption of having definable choice in $\MM$, and consider any $M$-definable topology. We also weaken the intersection assumption (in their work they actually observe that it suffices to have $k$-consistency for some $k$ in terms of $l(x)$ and $l(b)$) without the use of $\VC$ or forking literature. 


\item \label{itm:remark:peterzil_pillay_aschenbrenner_2} $(p,q)$-theorems are closely related to so-called Fractional Helly theorems (see~\cite{matousek04}), which branched from the classical Helly theorem. In its infinite version, this theorem states that any family of closed and bounded convex subsets of $\mathbb{R}^n$ that is $(n+1)$-consistent has nonempty intersection. With an eye towards definably extending Lipschitz maps, Aschenbrenner and Fischer proved (\cite{aschen_fischer_11}, Theorem B) a definable version of Helly's Theorem (i.e. for definable families) in definably complete expansions of real closed fields. 

Our Theorem~\ref{thm:compactness} and the arguments in Section $3.2$ in~\cite{aschen_fischer_11} allow a generalization of the o-minimal part of Aschenbrenner's and Fischer's definable Helly Theorem, by asking that the sets be definably compact and closed in any definable topology, instead of closed and bounded in the euclidean sense. Moreover, by appropriately adapting Corollary $2.6$ in~\cite{aschen_fischer_11}, one may show that every definable family of convex subsets of $M^n$ that is $(n+1)$-consistent extends to a definable type. 

\end{enumerate}
\end{remark}

\begin{remark} 
In~\cite[Section 6]{atw1}, Thomas, Walsberg and the author introduce the notion of \emph{definable net} $\gamma:(\Omega,\lleq)\rightarrow (X,\tau)$ to mean a definable map from a definable (downwards) directed set $(\Omega, \lleq)$ into a definable topological space $(X,\tau)$. A \emph{subnet} (a \emph{Kelley subnet}) of $\gamma$ is a net of the form $\gamma'= \gamma \circ \mu$ where $\mu:(\Omega', \lleq')\rightarrow (\Omega, \lleq)$ is a downward cofinal map on some directed set $(\Omega',\lleq')$. We say that such a net $\gamma'$ is definable if all of $(\Omega',\lleq')$, $\mu$ and $\gamma$ are definable.

Classically, a topological space is compact if and only if every net in it has a converging subnet. Following the classical proof of this result one may show that, in any model theoretic structure (regardless of the axiom of o-minimality), definable compactness implies that every definable net has a definable converging subnet. The reverse implication follows whenever the structure has definable choice. In our o-minimal $\MM$ we have that, by the relationship between \eqref{itm:compactness_1} and \eqref{itm:compactness_6} in Theorem~\ref{thm:compactness} (and more specifically by Lemma~\ref{lemma:choice_compact_spaces} below), this reverse implication also holds whenever the topology is Hausdorff. See~\cite[Corollary 44]{atw1} for a proof of the equivalence between definable compactness, curve-compactness, and the property that every definable net has a convergent definable subnet, in o-minimal expansions of ordered groups.

\end{remark}

We divide the proof of Theorem~\ref{thm:compactness} into a number of propositions. Note that the equivalence \eqref{itm:compactness_2}$\Leftrightarrow$\eqref{itm:compactness_2.5} is given by Proposition~\ref{prop:pre_specialization_compactness}. In Proposition~\ref{prop:compactnes_transversals} we prove, using results from previous sections, the equivalence between~\eqref{itm:compactness_1}, \eqref{itm:compactness_2.5}, \eqref{itm:compactness_3}, \eqref{itm:compactness_5} and~\eqref{itm:compactness_4}. In Proposition~\ref{prop:curve_compact_iff_directed_compact} we prove the implication~\eqref{itm:compactness_1}$\Rightarrow$\eqref{itm:compactness_6}, and reverse implication when $\tau$ is Hausdorff or $\MM$ has definable choice. 




\begin{proposition}\label{prop:compactnes_transversals}
Let $(X,\tau)$ be a definable topological space. The following are equivalent. 
\begin{enumerate}[(1)]
\item \label{itm:compactness_transversals_1} The topology $\tau$ is definably compact. 
\item \label{itm:compactness_transversals_1.5} Any definable family of $\tau$-closed sets that extends to a definable type has nonempty intersection.
\item \label{itm:compactness_transversals_2} Any consistent definable family of $\tau$-closed sets has a finite transversal in $X$.
\item \label{itm:compactness_transversals_3} Any definable family $\CC$ of nonempty $\tau$-closed sets with the $(m,n)$-property, where $m\geq n >\dim \cup\CC$, has a finite transversal in $X$.
\item \label{itm:compactness_transversals_4} Any definable family $\CC$ of nonempty $\tau$-closed sets with the $(m,n)$-property, where $m\geq n$ and $n$ is greater than the $\VC$-codensity of $\CC$, has a finite transversal in $X$.  
\end{enumerate}
\end{proposition}
\begin{proof}
We show that $\eqref{itm:compactness_transversals_1} \Rightarrow \eqref{itm:compactness_transversals_1.5}$, $\eqref{itm:compactness_transversals_1.5} \Rightarrow \eqref{itm:compactness_transversals_2} \wedge \eqref{itm:compactness_transversals_3} \wedge\eqref{itm:compactness_transversals_4}$, and $\eqref{itm:compactness_transversals_2} \vee \eqref{itm:compactness_transversals_3} \vee \eqref{itm:compactness_transversals_4} \Rightarrow \eqref{itm:compactness_transversals_1}$. 

Note that, if a downward directed family of sets has a finite transversal, then, by Fact~\ref{fact:downward_directed_family_2}, it has nonempty intersection. Hence \eqref{itm:compactness_transversals_2}, \eqref{itm:compactness_transversals_3} and~\eqref{itm:compactness_transversals_4} each imply~\eqref{itm:compactness_transversals_1}. The implications~\eqref{itm:compactness_transversals_1.5}$\Rightarrow$\eqref{itm:compactness_transversals_3} and~\eqref{itm:compactness_transversals_1.5}$\Rightarrow$\eqref{itm:compactness_transversals_4} follow from Theorem~\ref{thm:finite_tame_transversals} and Corollary~\ref{prop:vc_tame_transversals} respectively. The implication \eqref{itm:compactness_transversals_1.5}$\Rightarrow$\eqref{itm:compactness_transversals_2} follows because~\eqref{itm:compactness_transversals_3}$\Rightarrow$\eqref{itm:compactness_transversals_2} is obvious. It remains thus to show that~\eqref{itm:compactness_transversals_1}$\Rightarrow$\eqref{itm:compactness_transversals_1.5}. 




Suppose that $(X,\tau)$ is definably compact and let $\CC$ be a definable family of $\tau$-closed sets that extends to a definable type. By Lemma~\ref{lemma:existence_complete_directed_families_2.1}, let $\FF\subseteq p$ be a definable downward directed family of sets that is finer than $\CC$. Note that the family $\{cl_\tau(F) : F\in\FF\}$ is definable, downward directed and finer that $\CC$. By definable compactness there exists $x\in \bigcap\{cl_\tau(F) : F\in\FF\}$. In particular $x\in \cap\CC$.  
\end{proof}

\begin{remark}
It is worth noting which results from Sections~\ref{section:types} and~\ref{section:transversals} are used in proving the non-trivial implications in Proposition~\ref{prop:compactnes_transversals}. By virtue of Proposition~\ref{prop:pre_specialization_compactness}, let us refer to condition \eqref{itm:compactness_transversals_1.5} in Proposition~\ref{prop:compactnes_transversals} as ``type-compactness". 

We prove that definable compactness implies type-compatness (\eqref{itm:compactness_transversals_1}$\Rightarrow$\eqref{itm:compactness_transversals_1.5}) using Lemma~\ref{lemma:existence_complete_directed_families_2.1}. It is easy to see that the converse implication can be proved using Lemma~\ref{lemma:existence_complete_directed_families}, and in particular the fact that every definable downward directed family extends to a definable type, which, as noted in Remark~\ref{remark:existence_ddf_from_FTT}, is a condition weaker than Theorem~\ref{thm:finite_tame_transversals}.

Furthermore, we prove that type-compactness implies~\eqref{itm:compactness_transversals_2},~\eqref{itm:compactness_transversals_3} and~\eqref{itm:compactness_transversals_4} using Theorem~\ref{thm:finite_tame_transversals} and Corollary~\ref{prop:vc_tame_transversals}. Alternatively this can be done, and generalized to a broader dp-minimal setting, using Theorem~\ref{thm:forking_definable_types} due to Simon and Starchenko, as well as $\VC$-theory, in particular the Alon-Kleitman-Matou\v{s}ek $(p,q)$-theorem (Theorem~\ref{thm:vc_tranversal}) and known bounds on $\VC$-density in NIP theories (see~\cite{vc_density}) such as Theorem~\ref{cor:vc_density}.
\end{remark}

The following Corollary is a direct consequence of Propositions~\ref{prop:pre_specialization_compactness} and~\ref{prop:compactnes_transversals}.

\begin{corollary}\label{cor:definably_compact_iff_compact}
Suppose that $\MM$ extends $(\mathbb{R},<)$, the linear order of reals, and let $(X,\tau)$ be a definable topological space. Then $(X,\tau)$ is definably compact if and only if it is compact. 
\end{corollary}

One may use Corollary~\ref{cor:definably_compact_iff_compact} and  to show that, whenever $\MM$ expands the real line, the chatacterization of definable compactness given by Theorem~\ref{thm:compactness}~\eqref{itm:compactness_4} follows directly from the Alon-Kleitman-Matou\v{s}ek $(p,q)$-theorem (Theorem~\ref{thm:vc_tranversal}).

We now prove the connection between definable compactness and curve-compactness stated in Theorem~\ref{thm:compactness}. That is, that definable compactness implies curve-compactness in general, and that both notions are equivalent when the topology is Hausdorff or when the underlying structure $\MM$ has definable choice. This is Proposition~\ref{prop:curve_compact_iff_directed_compact}. We follow the proposition with an example of a non-Hausdorff topological space definable in the trivial structure $(M,<)$ that is curve-compact but not definably compact.  

In the case where $\MM$ extends an ordered field, the equivalence between curve- and type-compactness, as well as their equivalence with classical compactness whenever $\MM$ expands the real line, can be derived using the theory of tame pairs. This was done for definable metric spaces in~\cite[Proposition 6.6]{walsberg15} and in general in~\cite[Corollaries 47 and 48]{atw1}.

The next lemma allows us to apply definable choice in certain instances even when the underlying structure $\MM$ does not have the property.

\begin{lemma}[Definable choice in compact Hausdorff spaces]\label{lemma:choice_compact_spaces}
Let $C\subseteq M^m$ be a definable nonempty $\tau$-closed set in a curve-compact Hausdorff definable topological space $(X,\tau)$. Suppose that $\tau$ and $C$ are $A$-definable. Then there exists a point $x\in C \cap dcl(A)$, where $dcl(A)$ denotes de definable closure of $A$. 

Consequently, for any $A$-definable family $\{\varphi(u,M^m) : u\in\Omega\}$ of nonempty subsets of $X$, there exists an $A$-definable selection function $s:\Omega\rightarrow X$ such that $s(u)\in cl_\tau(\varphi(u,M^m))$ for every $u\in \Omega$. 
\end{lemma}
\begin{proof}
We prove the first paragraph of the lemma. The uniform result is derived in the usual way by the use of first-order logic compactness.

For this proof we adopt the convention of the one point euclidean space $M^0=\{0\}$. In particular any projection $M^{k}\rightarrow M^0$ is simply the constant function $0$ and any relation $E\subseteq M^0\times M^{k}$ is definable if and only if its restriction to $M^{k}$ is.  

Let $C$ and $\tau$ be as in the Lemma. Let $0\leq n\leq m$ be such that there exists an $A$-definable function $f:D\subseteq M^n\rightarrow C$, for $D$ a nonempty set. If $n$ can be chosen to be zero then the lemma follows. We prove that this is the case by backwards induction on $n$. 

Note that $n$ can always be chosen equal to $m$ by letting $f$ be the identity on $C$. Suppose that $0<n\leq m$. For any $x\in M^{n-1}$ let $D_x$ denote the fiber $\{t\in M : \al x,t\ar \in D\}$. For each $x\in\pi(D)$ let $s_x=\sup D_x$, and consider the $A$-definable set $H=\{x\in \pi(D): s_x\in D_x\}$. 

If $H\neq \emptyset$ then let $g$ be the map $x\mapsto f(s_x):H\rightarrow C$. 
If $H=\emptyset$ then let $g$ be the map $x\mapsto \tau\text{-}\lim_{t\rightarrow s_x^-} f(x,t): \pi(D)\rightarrow C$ which, by curve-compactness and Hausdorffness, is well defined. In both cases $g$ is an $A$-definable nonempty partial function $M^{n-1}\rightarrow C$. 
\end{proof}
Lemma~\ref{lemma:choice_compact_spaces} implies that, even when $\MM$ does not have definable choice, if an $A$-definable family of nonempty closed sets $\CC$ in a curve-compact Hausdorff $A$-definable topology has a finite transversal, then it also has one in $dcl(A)$. To prove this it suffices to note that, for any $k\geq 1$, the set of $k$-tuples of points corresponding to a transversal of $\CC$ is $A$-definable and closed in the product topology. 

\begin{proposition}\label{prop:curve_compact_iff_directed_compact}
Let $(X,\tau)$ be a definable topological space. If $(X,\tau)$ is definably compact then it is curve-compact. 

Suppose that either $\tau$ is Hausdorff or $\MM$ has definable choice. Then $(X,\tau)$ is definably compact if and only if it is curve-compact. 
\end{proposition}

We prove the left to right direction through a short lemma. 
\begin{lemma}\label{lemma:curve_compactness_implies_compactness}
Let $(X,\tau)$ be a definably compact definable topological space. Then $(X,\tau)$ is curve-compact. 
\end{lemma}
\begin{proof}
Let $\gamma:(a,b)\rightarrow X$ be a definable curve in $X$. Consider the definable family of $\tau$-closed nested sets $\CC_\gamma=\{ cl_\tau(\gamma[(a,t)]) : a<t<b\}$. By definable compactness there exists $x\in\cap\CC_\gamma$. Clearly $x\in \tau\text{-}\lim_{t\rightarrow a^-} \gamma(t)$. Similarly one shows that $\gamma$ $\tau$-converges as $t\rightarrow b$.    
\end{proof}

We now prove a simpler case of the left to right implication.  

\begin{lemma}\label{lemma:case_nested_family}
Let $(X,\tau)$ be a definable topological space. Suppose that either $\tau$ is Hausdorff or $\MM$ has definable choice. 
Let $\CC$ be a nested definable family of $\tau$-closed nonempty subsets of $X$. If $(X,\tau)$ is curve-compact then $\cap\CC\neq \emptyset$. 
\end{lemma}
\begin{proof}
Let $\CC=\{\varphi(u,M^m) : u\in\Omega\}$, $\Omega\subseteq M^n$, and $(X,\tau)$ be as in the lemma. We proceed by induction on $n$. 
 
\textbf{Base case: $n=1$.}

Consider the definable total preorder $\lleq$ in $\Omega$ given by $u\lleq v$ if and only if $\varphi(u,M^m) \subseteq \varphi(v,M^m)$. If $\Omega$ has a minimum $u$ with respect to $\lleq$ then $\varphi(u,M^n)\subseteq C$ for every $C\in\CC$ and the result is obvious. We suppose that $(\Omega,\lleq)$ does not have a minimum and consider the definable nested family of infinite sets $\{ (-\infty, u)_{\lleq} : u\in \Omega\}$. By the base case in the proof of Lemma~\ref{lemma:existence_complete_directed_families}, this family extends to a definable type with a basis of sets $\{(a,t): t>a\}$ for some $a\in M\cup\{-\infty\}$ or $\{(t,a) : t<a\}$ for some $a\in M\cup \{+\infty\}$. We consider the former case, being the proof for the latter analogous. This means that, for every $u\in \Omega$, there exists $t(u)>a$ such that $v\prec u$ for every $a<v<t(u)$.

If $\MM$ has definable choice or if $\tau$ is Hausdorff (Lemma~\ref{lemma:choice_compact_spaces}) there exists a definable function $f:\Omega\rightarrow \cup \CC$ satisfying that $f(u)\in \varphi(u,M^m)$ for every $u\in \Omega$. Recall that, for every $u,v\in \Omega$, if $v\lleq u$ then $\varphi(v,M^m)\subseteq \varphi(u,M^m)$, and in particular $f(v)\in \varphi(u,M^n)$. It follows that, for every $u\in \Omega$, there exists $t(u)>a$ such that $f(v)\in \varphi(u,M^n)$ for every $a<v<t(u)$. Let $b>a$ be such that $(a,b)\subseteq \Omega$ and $\gamma$ be the restriction of $f$ to $(a,b)$. We conclude that, for every $C\in \CC$, $\tau\text{-}\lim_{t\rightarrow a} \gamma(t) \in C$. 

\textbf{Inductive step: $n>1$.}

For any $x\in\pi(\Omega)$ let $\CC(x)$ denote the family $\{\varphi(x,t,M^m) : t\in \Omega_x\}$, and set $C(x):=\cap \CC(x)$. By the base case the definable family of $\tau$-closed sets $\DD=\{C(x) : x\in\pi(\Omega)\}$ does not contain the empty set. Clearly $\cap \DD=\cap \CC$. We observe that the family $\DD$ is nested and the result follows from the induction hypothesis. 

Given $x ,y\in \pi(\Omega)$, if for every $C\in \CC(x)$ there is $C'\in\CC(y)$ with $C'\subseteq C$ then $\cap \CC(y)  \subseteq \cap \CC(x)$. Otherwise there is $C\in\CC(x)$ such that, for every $C'\in\CC(y)$, it holds that $C\subseteq C'$, in which case $\cap \CC(x)\subseteq C \subseteq \cap \CC(y)$.     
\end{proof}

We may now prove Proposition~\ref{prop:curve_compact_iff_directed_compact}.

\begin{proof}[Proof of Proposition~\ref{prop:curve_compact_iff_directed_compact}]

By Lemma~\ref{lemma:curve_compactness_implies_compactness} we must only prove the right to left implication.
Let $(X,\tau)$, $X\subseteq M^m$, be a curve-compact Hausdorff definable topological space. Let $\CC\subseteq \mathcal{P}(X)$ be a definable downward directed family of subsets of $(X,\tau)$, not necessarily $\tau$-closed. We show that $\bigcap\{cl_\tau(C) : C\in\CC\}\neq \emptyset$. 

We proceed by induction on $n=\min\{\dim C : C\in\CC\}$. By Lemma~\ref{lemma:existence_complete_directed_families}, after passing to a finer family if necessary we may assume that $\CC$ is a complete family of cells. For the length of this proof, given two partial $\exR$-valued functions $f$ and $g$ let $(f,g)=\{\al x, t\ar: x\in dom(f)\cap dom(g) \text{ and } f(x)<t<g(x)\}$, relaxing thus the classical notation by allowing that $f$ and $g$ have different domains.  

If $n=0$ there exists a finite set in $\CC$ and so (Fact~\ref{fact:downward_directed_family_2}) $\bigcap\CC\neq \emptyset$. We now prove the case $n=m>0$. Hence suppose that every $C\in \CC$ is an open cell $C=(f_C,g_C)$, for definable continuous functions $f_C ,g_C:\pi(C)\rightarrow \exR$ with $f_C < g_C$. For any $C\in\CC$ let $D(C)=\bigcap \{cl_\tau(f_C,g_{C'}) : C'\in \CC\}$. We first show that these sets are nonempty.

Let us fix $C=(f,g)$ and, for any $x\in\pi(C)$, let $C^0(x)=\tau\text{-}\lim_{t\rightarrow f(x)^+} \al x,t\ar$. If $\tau$ is Hausdorff then this point is unique, otherwise we use definable choice to pick one such point definably on $x$. The definable set $C^0=\{C^0(x) : x\in\pi(C)\}$ has dimension $\dim(C)-1=n-1$. For any $C'=(f',g')\in\CC$, the definable set $\{x\in \pi(C)\cap \pi(C'): f(x)<g'(x)\}$ is nonempty, since otherwise we would have $C\cap C' =\emptyset$. It follows that $C^0 \cap cl_\tau(f,g')\neq \emptyset$. Note that the definable family $\{C^0 \cap cl_\tau(f, g_{C'}) : C'\in\CC\}$ is downward directed. By inductive hypothesis there is a point that belongs in the $\tau$-closure of $C^0 \cap cl_\tau(f, g_{C'})$, and in particular in $cl_\tau(f, g_{C'})$, for all $C'\in\CC$. So $D(C)\neq \emptyset$.

Note that, for every $C\in\CC$, $D(C)\subseteq cl_\tau(C)$. We now note that the definable family of nonempty $\tau$-closed sets $\{D(C) : C\in\CC\}$ is nested. Then, by Lemma~\ref{lemma:case_nested_family}, $\bigcap \{ D(C) : C\in\CC\}\neq \emptyset$, and thus $\bigcap\{ cl_\tau(C) : C\in\CC\}\neq \emptyset$. 

Let us fix $C_1=(f_1,g_1), C_2=(f_2,g_2)\in\CC$. We may partition $B=\pi(C_1)\cap \pi(C_2)$ by the definable sets 
\[
B_1=\{x\in B : f_1(x)\leq f_2(x)\} \text{ and } B_2=\{x\in B : f_1(x)>f_2(x)\}. 
\]
Since $\CC$ is complete there exists some $i\in\{1,2\}$ and $C\in \CC$ such that $\pi(C)\subseteq B_i$. Without loss of generality suppose that $i=1$, and let us fix $C_3\in\CC$ with $\pi(C_3)\subseteq B_1$. For any $C=(f,g)\in \CC$ let $C'=(f',g')\in\CC$ be such that $C'\subseteq C \cap C_3$. Clearly $(f_2,g')\subseteq (f_1,g')\subseteq (f_1,g)$. It follows that $D(C_2)\subseteq D(C_1)$. This completes the proof of the case $n=m>0$. 

Finally we prove the case $0<n<m$ by adapting the arguments above. Fix $\hat{C}\in\CC$ with $\dim \hat{C}=n<m$ and a projection $\hat{\pi}:M^m\rightarrow M^n$ such that $\hat{\pi}|_{\hat{C}}:\hat{C}\rightarrow \hat{\pi}(\hat{C})$ is a bijection. By passing to a finner subfamily if necessary we may assume that every set in $\CC$ is contained in $\hat{C}$. It follows that the definable downward directed family $\hat{\pi}(\CC)=\{\hat{\pi}(C) : C\in\CC\}$ contains only open cells in $M^n$. Note moreover that this family is complete. 

Set $h:=(\hat{\pi}|_{\hat{C}})^{-1}$. Then the proof in the case $n=m$ can be applied to $\hat{\pi}(\CC)$ instead of $\CC$, making sure to write throughout $cl_\tau h(\cdot)$ in place of $cl_\tau(\cdot)$, and letting $\hat{\pi}(C)^0=\{\tau\text{-}\lim_{t\rightarrow f(x)^+} h(x,t) : x\in\pi(\hat{\pi}(C))\}$ for any $\hat{\pi}(C)=(f,g)\in\hat{\pi}(\CC)$. It follows that $\bigcap\{ cl_\tau(C) : C\in\CC\}\neq \emptyset$. This completes the proof of the proposition.  
\end{proof}

The following is an example of a non-Hausdorff definable topological space in $(M,<)$ that is curve-compact but not definably compact.

\begin{example}
Let $X=\{\al x,y\ar\in M^2 : y<x\}$. 
Consider the family $\BB$ of sets 
\begin{align*}
A(x',x'',x''',y',y'',y''')=&\{\al x,y\ar : y<y', y<x\} \cup \\
&\{\al x,y\ar : y''<y<y''' \wedge (y<x<y''' \vee x'<x<x'' \vee x'''<x)\}
\end{align*}
definable uniformly over $y'<y''<y'''<x'<x''<x'''$.

\begin{figure}[h!]
\centering
\includegraphics[scale=0.6]{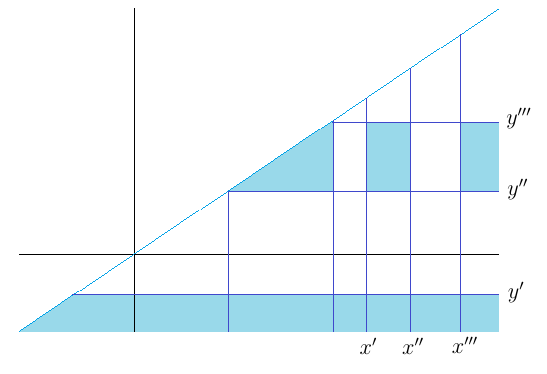}
\caption{\footnotesize Depicting (in blue) the set $A(x',x'',x''',y',y'',y''')$. \normalsize}
\end{figure}

Given any $A_0=A(x_0',x_0'',x_0''',y_0',y_0'',y_0''')$ and $A_1=A(x_1',x_1'',x_1''',y_1',y_1'',y_1''')$ in $\BB$ and any $\al x,y\ar\in A_1\cap A_2$, since every set in $\BB$ is open in the euclidean topology, we may find $y''<y<y'''<x'<x<x''$ such that the box $(x',x'')\times (y'',y''')$ is a subset of $A_1\cap A_2$. Let $y'<\min\{y'',y'_0,y'_1\}$ and $x'''>\max\{x'',x'''_0,x'''_1\}$. Then $\al x,y\ar\in A(x',x'',x''',y',y'',y''')\subseteq A_1 \cap A_2$. Hence the family $\BB$ is a definable basis for a topology $\mytau$. 

This topology is $T_1$, i.e. every singleton is closed. For every $y\in M$, $\mytau\text{-}\lim_{t\rightarrow y^+} \al t, y\ar = \mytau\text{-}\lim_{t\rightarrow +\infty} \al t, y\ar = (M\times\{y\})\cap X$, and, for every $x\in M$, $\mytau\text{-}\lim_{t\rightarrow x^-} \al x, t\ar=(M\times\{x\})\cap X$ and $\mytau\text{-}\lim_{t\rightarrow -\infty} \al x, t\ar=X$. In particular $\mytau$ is not Hausdorff. 

Now suppose that $\MM=(M,<)$. By quantifier elimination we know that in this structure any definable partial map $M\rightarrow M$ is piecewise either constant or the identity. Let $\gamma=(\gamma_0,\gamma_1):(a,b)\rightarrow X$ be an injective definable curve in $X$. Let $I=(a,c)\subseteq (a,b)$ be an interval where $\gamma_0$ and $\gamma_1$ are either constant or the identity. Since the graph of the identity is disjoint from $X$ and $\gamma$ is injective it must be that $\gamma_i$ is constant and $\gamma_{1-i}$ is the identity on $I$ for some $i\in\{0,1\}$. 

Suppose that $\gamma_1|_I$ has a constant value $y$. Then, by the observations made above about the topology, $\gamma$ $\mytau$-converges as $t\rightarrow a$ to either $\al a , y \ar$ (if $y<a$) or $(M\times\{y\}) \cap X$ (if $a=y$). On the other hand, if $\gamma_0|_I$ has a constant value $x$, then $\gamma$ $\mytau$-converges as $t\rightarrow a$ to either $\al x , a \ar$ (if $a>-\infty$) or the whole space $X$ (if $a=-\infty$). Treating the limit as $t\rightarrow b$ similarly allows us to conclude that $\gamma$ is $\mytau$-completable. Hence the space $(X,\mytau)$ is curve-compact.  

Meanwhile, the definable downward directed family of $\mytau$-closed sets $\{ X \cap (R\times [u,+\infty)) : u\in M\}$ has empty intersection. So $(X,\mytau)$ is not definably compact. 
\end{example}

We end the paper with a question. Notice that all the characterizations of definable compactness in Theorem~\ref{thm:compactness} are upfront expressible with infinitely many sentences (possibly with parameters) in the language of $\MM$ (you have to check all relevant definable families of closed sets or all definable curves). In~\cite{pet_stein_99} it is shown that, in the euclidean topology, a set is curve-compact if and only if it is closed and bounded. Being closed and bounded is expressible by a single formula. Moreover, given a definable family of sets with the euclidean topology, the subfamily of those that are closed and bounded is definable. We ask the following.   
\begin{question}
\begin{enumerate}[(i)]
\item Can definable compactness of an arbitrary definable topology be expressed with a single formula in the language of $\MM$? 

\item More generally, given a definable family of definable topological spaces, is the subfamily of all which are definably compact definable?
\end{enumerate} 
\end{question}



\newpage
\appendix

\section{Appendix: Proof of The Marker-Steinhorn Theorem}\label{section:appendix A}

In this appendix we use the approach of various proofs in this paper, based on the the preorder described in Subsection~\ref{remark:preorder}, to extract from~\cite{mark_stein_94} a short proof of the Marker-Steinhorn Theorem (that avoids a treatment by cases). It is worth noting that a similar approach via definable linear orders was used by Walsberg~\cite{walsberg19} to prove the theorem in o-minimal expansions of ordered groups.  

Recall that $\MM$ denotes an o-minimal structure and that, unless stated otherwise, definable means in $\MM$ over $M$. Recall that an elementary extension $\NN$ of $\MM$ is \emph{tame} if, for every $s\in N$, the set $\{t\in M : t<s\}$ has a supremum in $M\cup\{+\infty\}$. This is also referred to as $\MM$ being \emph{Dedekind complete} in $\NN$~\cite{mark_stein_94}.

In full generality the Marker-Steinhorn theorem shows that a type $p\in S_n(M)$ is definable if and only if it is realized in some tame extension. We state and prove the \textquotedblleft if" direction, and direct the reader to~\cite[Theorem 2.1]{mark_stein_94} for the full result.

The approach of the proof below involves using cell decomposition to reduce the problem to showing that any cut in a definable preordered set that is realized in a tame extension is definable. This of course is obvious for the preordered set $(M,\leq)$, by definition of tameness. 

\begin{theorem}[Marker-Steinhorn Theorem]\label{thm:marker_steinhorn}
Let $\NN=(N,\ldots)$ be a tame extension of $\MM$. For every $a\in N^n$, the type $\tp(a/M)$ is definable. 
\end{theorem}
\begin{proof}
We must prove that, for every $a\in N^n$ and formula $\varphi(u,a)$, $l(u)=m$, the set $\varphi(M^m,a)=\{ u\in M^{m} : \NN\models \varphi(u,a)\}$ is definable. We do this by induction on $n$ and $m$, where in the inductive step we assume that it holds for any $\al  n',m'\ar$ smaller than $\al n,m\ar$ in the lexicographic order. We may assume that $a\notin M^n$. 

The case $n=1$ (for any $m$) follows easily from o-minimality and tameness. In particular, let $s_a$ be the supremum in $M\cup\{\pm\infty\}$ of $(-\infty,a)\cap M$. If $s_a<a$ then $\tp(a/M)$ has a definable basis of the form $\{(s_a,t) : s_a<t\in M\}$, otherwise it has a definable basis of the form $\{(t,s_a) : s_a>t\in M\}$. 

Suppose that $n>1$ and let $a=\al c,d\ar\in N^{n-1}\times N$. Let $\psi_i(N^m,N^n)$, $0\leq i \leq l$, be a ($0$-definable) cell decomposition of $\varphi(N^m, N^n)$. For every $u\in M^m$, the set $\varphi(u,N^n)$ is partitioned by the cells $\psi_i(u,N^n)$, $0\leq i \leq l$. In particular $\NN\models \varphi(u,a)$ if and only if $\NN\models \psi_i(u,a)$ for some $i$. So to prove the theorem it suffices to pass to an arbitrary $0\leq i \leq l$ and show that $\psi_i(M^m,a)$ is definable. Hence, without loss of generality we assume that all the sets $\varphi(u,N^n)$ are cells.  

By induction hypothesis $\tp(c/M)$ is definable.
Suppose that there exists an $M$-definable partial function $f:N^{n-1}\rightarrow N$ with $f(c)=d$. Then, for every $u\in N^m$, we have that $\NN\models \varphi(u,a)$ if and only if $\NN\models \exists t\, (\varphi(u,c,t)\wedge (f(c)=t))$. And so the result follows. 


Hence onwards we assume that there exists no $M$-definable partial function $f$ with $f(c)=d$. In particular we have that every $u\in M^m$ with $\NN\models \varphi(u,a)$ satisfies that $\varphi(u,N^n)$ is a cell of the form $(f_u,g_u)$ for two definable continuous functions $f_u$ and $g_u$. Let 
\[
\Omega=\{u\in M^{m} : \varphi(u,N^n) \text{ is of the form } (f_u,g_u) \text{ and } \NN\models \exists t\, \varphi(u,c,t)\}.
\]
Clearly $\varphi(M^m,a)\subseteq \Omega$. Since $\tp(c/M)$ is definable the set $\Omega$ is definable. We prove that the sets $\{ u\in \Omega: f_u(c)<d\}$ and $\{ u\in \Omega: g_u(c)>d\}$ are definable. Then their intersection equals $\varphi(M^m,a)$. The proof of definability is the same for both sets, so we show it only for the former. 

Let $\lleq$ be the total preorder on $\Omega$ induced by $\{ f_u(c) : u\in \Omega\}$ described in Subsection~\ref{remark:preorder}, i.e. for $u,w\in\Omega$, $u\lleq w$ if and only if $f_u(c)\leq f_w(c)$. By definability of $\tp(c/M)$ the preordered set $(\Omega,\lleq)$ is definable. Let $P=\{ u\in \Omega: f_u(c)<d\}$ and $Q=\{ u\in \Omega: f_u(c)>d\}=\Omega\setminus P$. We must prove that $P$ (equivalently $Q$) is definable. If $P$ has a supremum in $\Omega\cup\{\pm\infty\}$ with respect to $\lleq$ then the result is immediate, so we assume otherwise. In particular we have that, for every $u\in P$ and $v\in Q$, $\dim (u,v)_{\lleq}>0$. 

Note that, to prove the definability of $P$, it suffices to show that there exists a definable set $P'\subseteq P$ that is cofinal in $P$, since then $P=\{u\in \Omega : u\lleq v \text{ for some } v\in P'\}$. Similarly it is enough to show the existence of a defiable $Q'\subseteq Q$ coinitial in $Q$. So we may always pass to a subset $\Omega'\subseteq \Omega$ such that either $\Omega'\cap P$ is cofinal in $P$ or $\Omega'\cap Q$ is coinitial in $Q$, and then prove definability of $\Omega'\cap P$ or $\Omega'\cap Q$. Hence, by passing to a subset if necessary, we may assume that, for any $u\in P$ and $v\in Q$,
\begin{equation}\label{eqn:MS_1}
\dim (u,v)_{\lleq} = \dim \Omega. \tag{$\star$} 
\end{equation}
Moreover, note that we may also pass to a set in any given definable finite partition of $\Omega$. In particular, by cell decomposition, we assume that $\Omega$ is a cell.  

Suppose that $m=1$. For each $u\in \Omega(\NN)$ with $c\in dom(f_u)$ (note that this includes every $u$ in $\Omega$) let $\hat{f}(u)=f_u(c)$. Then $\hat{f}$ is definable in $\NN$ over $M\cup\{c\}$. By o-minimality there exists a partition (definable over $M\cup \{c\}$) of the domain of $\hat{f}$ into points and intervals such that, on each interval, $\hat{f}$ is continuous and either constant or strictly monotonic. Since $\tp(c/M)$ is definable the intersections of these cells with $\Omega$ are definable. Note that, on any such intersection, the restriction of the preorder $\lleq$ is either $\leq$, $\geq$, or $\leq \cup \geq$ (the trivial relation where any two points are indistinguishable), depending respectively on whether $\hat{f}$ is strictly increasing, decreasing or constant. We fix one such interval $I$ and show that $I\cap P$ is definable. 

If $I\cap P=I\cap \Omega$ or $I\cap Q=I\cap \Omega$ then the result is immediate. Otherwise there exist $u,v\in I \cap \Omega$ such that $\hat{f}(u)<d$ and $\hat{f}(v)>d$. By continuity there must exist $r$ in the subinterval of $I$ with endpoints $u$ and $v$ with $\hat{f}(r)=d$. By tameness $J=(-\infty,d) \cap M$ is definable. Finally note that $\hat{f}|_I$ is not constant, and thus it is strictly monotonic. If it is increasing then it must be that $J \cap I= P \cap I$ and otherwise $J\cap I=Q\cap I$. 

Now suppose that $m>1$. For every $x$ in the projection $\pi(\Omega)$ to the first $m-1$ coordinates, let $\lleq_x$ be the definable preorder on the fiber $\Omega_x$ given by $s \lleq_x t$ if and only if $\al x,s\ar \lleq \al x,t\ar$. Note that, following the arguments in the case $m=1$, the sets $P_x$ and $Q_x$ are definable, and moreover $\Omega_x$ can be partitioned into finitely many points and intervals where the restriction of $\lleq_x$ is either $\leq$, $\geq$, or $\leq \cup \geq$. 

If there exists $x\in \pi(\Omega)$ such that $\{x\}\times P_x$ is cofinal in $P$ or $\{x\} \times Q_x$ is coinitial in $Q$, then we are done. Suppose otherwise. We complete the proof by partitioning $\Omega$ into finitely many definable sets with the following property. For each set $\Sigma$ in the partition and $x\in \pi(\Sigma)$, either $\Sigma_x \subseteq P_x$ or $\Sigma_x\subseteq Q_x$. Observe that the set $\Theta$ of all $x\in \pi(\Sigma)$ such that $\Sigma_x \subseteq P_x$ is  described by
\[
x\in \pi(\Sigma),\, \forall t\in \Sigma_x,\, f_{(x,t)}(c) < d,
\]
so, by induction hypothesis (applied in the case $\al n, m-1\ar$), this set is definable. It follows that $\Sigma\cap P = \bigcup_{x\in\Theta} (\{x\}\times \Sigma_x)$ is definable, and we may conclude that $P$ is definable. 

Recall that $\Omega$ is a cell. If it is defined as the graph of a function then, by taking $\Sigma$ in the above paragraph to be $\Omega$, we are done, so we assume otherwise. 
Let $x\in \pi(\Omega)$. By~\eqref{eqn:MS_1}, since $\{x\}\times P_x$ is not cofinal in $P$ and $\{x\} \times Q_x$ is not coinitial in $Q$, the dimension of 
\[
\{u \in \Omega : \{x\}\times P_x \prec u \prec \{x\} \times Q_x\}
\] 
equals $\dim \Omega$. Consider $r\in \Omega_x$ to be the right endpoint of some maximal subinterval $I'$ of $P_x$. Then in particular $r$ is the left endpoint of a subinterval $I''$ of $Q_x$. If $r\notin P_x$ then the set $\{u\in \Omega : \{x\}\times I' \prec u \prec \al x,r\ar\}$ has dimension $\dim(\Omega)$. If however $r\in P_x$, then the set $\{u\in \Omega : \al x,r\ar \prec u \prec \{x\}\times I''\}$ has dimension $\dim(\Omega)$. If $r$ is the right endpoint in $\Omega_x$ of a a maximal subinterval of $Q_x$ then the analogous holds. Note that, if there exists $s, t\in \Omega_x$ with $s\in P_x$ and $t\in Q_x$, there will always be some $r$ in the closed interval between $s$ and $t$ with the described properties. 

For any $u=\al x,t\ar\in \Omega$, let $L(x,t)$ be the set of $v\in \Omega$ such that either $\{x\}\times I' \prec v \prec u$ or $u \prec v \prec \{x\}\times I'$ for some interval $I'\subseteq \Omega_x$ with right endpoint $t$. Similarly let $U(x,t)$ be the set of $v\in \Omega$ such that either $\{x\}\times I'' \prec v \prec u$ or $u \prec v \prec \{x\}\times I''$ for some interval $I''\subseteq \Omega_x$ with left endpoint $t$. These sets are definable uniformly on $u\in \Omega$. Let $\Lambda$ be the definable set of all $u\in\Omega$ such that $U(u)\cup L(u)$ has dimension $\dim(\Omega)$. By the above paragraph, for every $x\in \pi(\Omega)$ and $s,t\in \Omega_x$, if $s\in P_x$ and $t\in Q_x$, then there is some $r$ in the closed interval between $s$ and $t$ such that $\al x, r\ar \in \Lambda$.  

We now show that, for every $x\in \pi(\Omega)$, the fiber $\Lambda_x$ is finite. Then the proof is completed by taking any finite cell partition of $\Omega$ compatible with $\Lambda$, since, for any $\Sigma$ in said partition and $x\in \pi(\Sigma)$, the fiber $\Sigma_x$ is going to be either a point or an interval contained in $P_x$ or $Q_x$.

Towards a contradiction suppose that $\Lambda_x$ is infinite for some $x\in \Omega$. Let $J'$ be a subinterval of $\Lambda_x$ where $\lleq_x$ is either $\leq$, $\geq$ or $\leq \cup \geq$. If $\lleq_x$ equals $\leq \cup \geq$ then, for every $t\in J'$, the sets $U(x,t)$ and $L(x,t)$ are empty, contradicting that $\dim (L(x,t) \cup U(x,t))=\dim \Omega$. Suppose that $\lleq_x$ is either $\leq$ or $\geq$. Note that, for any distinct $s,t\in J'$, the sets $L(x,s)$, $U(x,s)$, $L(x,t)$ and $U(x,t)$ are pairwise disjoint. Using the fact that $\dim (L(x,t) \cup U(x,t))=\dim \Omega$ for every $t\in J'$, and applying the Fiber Lemma for o-minimal dimension \cite[Chapter 4, Proposition 1.5 and Corollary 1.6]{dries98}, we derive that  
\[
\dim \left(\bigcup_{t\in J'} (L(x,t) \cup U(x,t))\right) > \dim(\Omega),
\]
contradiction.  
\end{proof} 

\section{Appendix (by Will Johnson): $A$-definability in Lemma~\ref{lemma:existence_complete_directed_families}}
\label{section:appendix B}

Let $M$ be an o-minimal structure\footnote{In this appendix we resort to the convention of using the same terminology for a structure and its universe.}.
Let $\mathcal{S} \subseteq \mathcal{P}(M^n)$ be an $A$-definable
downward directed family.
\begin{itemize}
\item By Theorem~\ref{thm:types}, $\mathcal{S}$ can be extended to a
  definable type $p$ with a definable basis.
\item By Proposition~\ref{prop:parameter_version_existence_downward_directed_families}, $\mathcal{S}$ can be extended to an $A$-definable type.
\end{itemize}
It is natural to ask whether we can do both at the same time.  Can
$\mathcal{S}$ be extended to an $A$-definable type with a definable
basis?  As we will see below, the answer is \textbf{no} in general.

This shows that a parameter version of Lemma~\ref{lemma:existence_complete_directed_families} is not possible: in
our counterexample, there is no \emph{complete} $A$-definable downward
directed family $\mathcal{S}'$ finer than $S$.  Otherwise,
$\mathcal{S}'$ would be a definable basis for an $A$-definable type
$p$ extending $S$.

\subsection{Types with definable bases}

Let $M$ be any structure, not necessarily o-minimal.

Recall that, for any $M'\succeq M$, a type $p\in S_x(M')$ is an \emph{heir} of $p|_M$ if, for every formula $\varphi(x,a,c)\in p$, where $a$ and $c$ are tuples of parameters from $M'$ and $M$ respectively, there exists $b\in M^{l(a)}$ such that $\varphi(x,b,c)\in p$. It is easy to see that the unique heir of a definable type is its canonical extension using the same defining scheme.

\begin{lemma}\label{sdef}
  Let $p(x)$ be an $n$-type over $M$ with a definable basis.
  \begin{enumerate}
  \item \label{sdef1} $p(x)$ is a definable type.
  \item \label{sdef2} If $M' \succeq M$ and $p^{M'}$ is the heir of
    $p$, then $p^{M'}$ has a definable basis.
  \item \label{sdef3} Suppose $M$ is $\aleph_1$-saturated.  Let
    $\mathcal{S}_1, \mathcal{S}_2, \ldots$ be a sequence of definable downward
    directed families of subsets of $M^n$. Suppose that
    $\mathcal{S}_{i+1}$ is finer than $\mathcal{S}_i$ for all $i$, and
    suppose $\bigcup_{i = 1}^\infty \mathcal{S}_i$ generates
    $p(x)$.  Then there is some $i_0 \in \mathbb{N}$ such that
    $\mathcal{S}_{i_0}$ is finer than $\mathcal{S}_i$ for all $i$. In particular $\mathcal{S}_{i_0}$ generates $p(x)$.
  \end{enumerate}
\end{lemma}
\begin{proof}
  ~
  \begin{enumerate}
  \item Easy. 
  \item For any definable family $\mathcal{F}$ in $M$, let $\FF(M')$ denote its interpretation in $M'$.  Let
    $\mathcal{D}$ be the definable downward directed family generating
    $p(x)$.  For any definable family $\mathcal{F}$ in $M$, we
    have
    \begin{equation*}
      \forall X \in \mathcal{F} ~ \exists Y \in \mathcal{D} : (Y
      \subseteq X~\vee~Y \subseteq M^n \setminus X).
    \end{equation*}
    This is expressed by a first-order formula, so it remains true in
    $M'$.  Therefore $\mathcal{D}(M')$ generates an $n$-type
    over $M'$, which is clearly $p^{M'}$.
  \item The families $\mathcal{D}$ and $\bigcup_{i = 1}^\infty
    \mathcal{S}_i$ generate the same type, so $\mathcal{D}$
    is finer than $\bigcup_{i = 1}^\infty \mathcal{S}_i$.  Therefore
    \begin{equation}
      \forall i \in \mathbb{N}: (\mathcal{D} \text{ is finer than }
      \mathcal{S}_i). \label{need-later}
    \end{equation}
    Additionally, the union $\bigcup_{i = 1}^\infty \mathcal{S}_i$
    is finer than $\mathcal{D}$.  Therefore
    \begin{equation*}
      \forall X \in \mathcal{D} ~ \exists i \in \mathbb{N} ~ \exists Y \in
      \mathcal{S}_i : Y \subseteq X.
    \end{equation*}
    By $\aleph_1$-saturation, there is some finite $N$ such that
    \begin{equation*}
      \forall X \in \mathcal{D} ~ \exists i \le N ~ \exists Y \in
      \mathcal{S}_i : Y \subseteq X.
    \end{equation*}
    Thus $\bigcup_{i = 1}^N \mathcal{S}_i$ is finer than $\mathcal{D}$.  Now
    $\mathcal{S}_N$ is finer than $\mathcal{S}_i$ for all $i \le N$, so
    $\mathcal{S}_N$ is finer than $\bigcup_{i = 1}^N \mathcal{S}_i$.  Then
    \begin{equation*}
      \mathcal{S}_N \text{ is finer than } \left(\bigcup_{i = 1}^N
      \mathcal{S}_i\right) \text{ is finer than } \mathcal{D} \text{ is finer than
      } \left(\bigcup_{i = 1}^\infty \mathcal{S}_i \right),
    \end{equation*}
    by (\ref{need-later}).  So $\mathcal{S}_N$ is finer than $\bigcup_{i =
      1}^\infty \mathcal{S}_i$, which implies that $\mathcal{S}_N$
    is finer than $\mathcal{S}_i$ for all $i$. \qedhere
  \end{enumerate}
\end{proof}

\subsection{Divisible ordered abelian groups}

For any linearly ordered set $(S,\le)$, let $\Cut(S)$ be the set of
cuts in $S$, i.e., downward-closed subsets $\xi \subseteq S$.  Let
$\overline{S}$ denote the disjoint union $S \sqcup \Cut(S)$.  The
linear order on $S$ extends to a linear order on $\overline{S}$ in which
\begin{align*}
  x < \xi & \iff x \in \xi &&\text{ for } x \in S, ~ \xi \in \Cut(S) \\
  \xi_1 \le \xi_2 & \iff \xi_1 \subseteq \xi_2 &&\text{ for } \xi_1, \xi_2 \in \Cut(S).
\end{align*}
We let $-\infty$ and $+\infty$ denote the least and greatest elements
of $\Cut(S)$, i.e., $\emptyset$ and $S$.  If $x \in S$, we let $x^+$
and $x^-$ denote the cuts immediately above and below $x$:
\begin{align*}
  x^+ &:= \{y \in S : y \le x\} \\
  x^- &:= \{y \in S : y < x\}.
\end{align*}
When $S = \mathbb{Q}$, we identify an irrational number $\alpha \in \mathbb{R}
\setminus \mathbb{Q}$ with the corresponding cut.  Thus
\begin{equation*}
  \overline{\mathbb{Q}} = \mathbb{R} \cup \{\pm \infty\} \cup \{q^+, ~ q^- : q \in \mathbb{Q}\}.
\end{equation*}
Fix a divisible ordered abelian group $(M,+,\le)$.  Embed $M$ into a
monster model $\mathbb{U}$.
\begin{lemma} \label{p-type}
  For any $\xi \in \Cut(\mathbb{Q})$, there is a 2-type $p_\xi(x,y)$
  over $M$ generated by the following formulas:
  \begin{itemize}
  \item $x > t$, for each $t \in M$.
  \item $y > qx + s$, for each $q \in \mathbb{Q}$ with $q < \xi$ and $s \in M$.
  \item $y < qx+s$, for each $q \in \mathbb{Q}$ with $q > \xi$ and $s \in M$.
  \end{itemize}
\end{lemma}
\begin{proof}
  Take $a \in \mathbb{U}$ with $a > M$.  The prime model $M[a]$ is then $\mathbb{Q}
  \cdot a \oplus M$, which is isomorphic to $\mathbb{Q} \times M$ with the
  lexicographic order.  Take $b \in \mathbb{U}$ with
  \begin{align*}
    b &> qa + M \text{ for } q < \xi \\
    b &< qa + M \text{ for } q > \xi.
  \end{align*}
  Then $(a,b)$ realizes $p_\xi(x,y)$, so $p_\xi(x,y)$ is a consistent partial type.
  For completeness, let $(a',b')$ be another pair realizing
  $p_\xi(x,y)$.  Then $a' > M$, because of the formulas $x > t$.  By
  o-minimality, $\tp(a'/M) = \tp(a/M)$.  Moving $(a',b')$ by an
  element of $\Aut(\mathbb{U}/M)$, we may assume $a' = a$.  Then for any $q
  \in \mathbb{Q}$ and $s \in M$, we have
  \begin{align*}
    b > qa + s & \iff \xi > q \iff b' > qa + s \\
    b < qa + s & \iff \xi < q \iff b' < qa + s.
  \end{align*}
  Therefore $b$ and $b'$ determine the same cut in $M[a]$.  By
  o-minimality, $\tp(b/Ma) = \tp(b'/Ma)$.  Thus $\tp(a,b/M) =
  \tp(a',b'/M)$, as desired.  So $p_\xi$ is a (complete) type.
\end{proof}
\begin{lemma}\label{r-type}
  For any $q \in \mathbb{Q}$ and $\xi \in \overline{M}$, there is a
  2-type $r_{q,\xi}(x,y)$ over $M$ generated by the following
  formulas:
  \begin{itemize}
  \item $x > t$, for each $t \in M$.
  \item $y \le qx + s$ for $s \in M$ with $s \ge \xi$.
  \item $y \ge qx + s$ for $s \in M$ with $s \le \xi$.
  \end{itemize}
\end{lemma}
\begin{proof}
  Similar to Lemma~\ref{p-type}.
\end{proof}
\begin{lemma} \label{classify}
  Let $p(x,y)$ be a 2-type over $M$ containing the formulas
  $x > t$ for all $t \in M$.  Then $p(x,y)$ is one of the types
  $p_\xi$ or $r_{q,\xi}$ from Lemma~\ref{p-type},\ref{r-type}.
\end{lemma}
\begin{proof}
  Take $(a,b) \in \mathbb{U}^2$ realizing $p(x,y)$.  Then $a > M$, so the
  prime model $M[a] = \mathbb{Q} \cdot a \oplus M$ is isomorphic to $\mathbb{Q}
  \times M$ with the lexicographic order.  Then $b$ determines a cut
  in $\mathbb{Q} \times M$.  If $b$ lies in the middle of some coset $q \cdot
  a + M$, then $b$ realizes $r_{q,\xi}$ for some cut $\xi \in
  \overline{M}$.  Otherwise, $b$ realizes $p_\xi$ for some $\xi \in
  \Cut(\mathbb{Q})$.
\end{proof}
\begin{remark} \label{overlap}
  There is some overlap between the $p_\xi$ and $r_{q,\xi}$.
  Specifically, $p_{q^+} = r_{q,+\infty}$ and $p_{q^-} =
  r_{q,-\infty}$.
\end{remark}
\begin{remark}
  In the notation of Definition~\ref{definition:construction_types}, the types $r_{q,m}$,
  $r_{q,m^+}$, and $r_{q,m^-}$ are $f|_p$, $f^+|_p$, and
  $f^-|_p$ respectively, where $f(x) = qx + m$ and $p$ is the 1-type at $+\infty$.
\end{remark}

When we need to make the dependence on $M$ explicit, we write
$p_{\xi}$ and $r_{q,\xi}$ as $p^M_{\xi}$ and $r^M_{q,\xi}$.
\begin{proposition} \label{heirs}
  Let $M'$ be an elementary extension of $M$.
  \begin{enumerate}
  \item \label{first-case} For $\xi \in \Cut(\mathbb{Q})$, $p^{M'}_\xi$ is the unique
    heir of $p^M_\xi$.
  \item \label{second-case} For $q \in \mathbb{Q}$ and $m \in M$,
    \begin{itemize}
    \item $r^{M'}_{q,m}$ is the unique heir of $r^M_{q,m}$.
    \item $r^{M'}_{q,m^+}$ is the unique heir of $r^M_{q,m^+}$.
    \item $r^{M'}_{q,m^-}$ is the unique heir of $r^M_{q,m^-}$.
    \end{itemize}
  \end{enumerate}
\end{proposition}
\begin{proof}
  Let $p'(x,y)$ be an heir of $p^M_\xi$ over $M'$.  Then $p'(x,y)$ contains the formulas
  \begin{gather*}
    \{x > t : t \in M'\} \\
    \{y > qx + s : q \in \mathbb{Q}, ~ q < \xi, ~ s \in M'\} \\
    \{y < qx + s : q \in \mathbb{Q}, ~ q > \xi, ~ s \in M'\}
  \end{gather*}
  because it is an heir.  By Lemma~\ref{p-type}, $p'$ must be
  $p^{M'}_\xi$.  This proves part (\ref{first-case}), and part
  (\ref{second-case}) is similar.
\end{proof}
A type $p$ is definable if and only if it has a unique heir
\cite[Theorem~11.7]{Poizat}, in which case the unique heir is the
canonical extension using the same defining scheme as $p$.  Therefore
Proposition~\ref{heirs} implies that $p_\xi$, $r_{q,m}$, $r_{q,m^+}$,
and $r_{q,m^-}$ are definable types, and specifies their canonical
extensions.

\begin{proposition} \label{classify2}
  Let $p(x,y)$ be a 2-type over $M$ containing the formulas $x > t$
  for all $t \in M$.  Suppose $p(x,y)$ has a definable basis.  Then
  $p(x,y)$ is $r_{q,m}$, $r_{q,m^+}$, or $r_{q,m^-}$ for some $q \in
  \mathbb{Q}$ and $m \in M$.
\end{proposition}
\begin{proof}
  By Lemma~\ref{sdef}.\ref{sdef1}, the type $p$ is definable.
  By Lemma~\ref{classify} and Remark~\ref{overlap}, $p$ must be one of the following:
  \begin{itemize}
  \item $p_\xi$ for some $\xi \in \Cut(\mathbb{Q})$.
  \item $r_{q,\xi}$ for some $q \in \mathbb{Q}$ and $\xi \in \overline{M}
    \setminus \{+\infty,-\infty\}$.
  \end{itemize}
  In the second case, there are restrictions on $\xi$.  The set
  \begin{equation*}
    \{s \in M : (y \le qx + s) \in p(x,y)\} = \{s \in M : s \le \xi\}
  \end{equation*}
  must be definable.  By o-minimality, $\xi$ must be one of
  $+\infty, -\infty, m, m^+,$ or $m^-$ for some $m \in M$.  By
  Remark~\ref{overlap}, we conclude that $p$ is one of the following:
  \begin{itemize}
  \item $p_\xi$ for some $\xi \in \Cut(\mathbb{Q})$
  \item $r_{q,m}$ or $r_{q,m^+}$ or $r_{q,m^-}$ for some $q \in \mathbb{Q}$
    and $m \in M$.
  \end{itemize}
  Let $M'$ be an $\aleph_1$-saturated elementary extension of $M$.  By Proposition~\ref{heirs},
  \begin{align*}
    p^{M'}_\xi &\text{ is the unique heir of } p^M_\xi \\
    r^{M'}_{q,m} &\text{ is the unique heir of } r^M_{q,m} \\
    r^{M'}_{q,m^+} &\text{ is the unique heir of } r^M_{q,m^+} \\
    r^{M'}_{q,m^-} &\text{ is the unique heir of } r^M_{q,m^-}.    
  \end{align*}
  By Lemma~\ref{sdef}.\ref{sdef2}, we may replace $M$ with $M'$, and
  therefore assume that $M$ is $\aleph_1$-saturated.

  Now we apply Lemma~\ref{sdef}.\ref{sdef3} to rule out the case of
  $p_\xi$.  We break into cases depending on $\xi$:
  \begin{itemize}
  \item Suppose $p = p_{+\infty}$.  Define
    \begin{align*}
      V_{q,t} &= \{(x,y) : x > t \text{ and } y > qx\} &&
      \text{ for } q \in \mathbb{Q}, ~ t \in M \\
      \mathcal{V}_q &= \{V_{q,t} : t \in M\} &&\text{ for } q \in \mathbb{Q}.
    \end{align*}
    Each family $\mathcal{V}_q$ is a definable downward directed
    family.  One can verify that $\mathcal{V}_q$ is finer than
    $\mathcal{V}_{q'}$ if and only if $q \ge q'$.  Moreover, the union
    $\bigcup_{n = 1}^\infty \mathcal{V}_n$ generates $p_{+\infty}$.
    By Lemma~\ref{sdef}.\ref{sdef3}, $p_{+\infty}$ does not have a
    definable basis.
  \item The case of $p = p_{-\infty}$ is handled similarly.
  \item Suppose $p = p_\alpha$ for some irrational $\alpha \in \mathbb{R}
    \setminus \mathbb{Q}$.  Take two sequences of rational numbers
    \begin{align*}
      q_1 &> q_2 > q_3 > \cdots \\
      q_1' &< q_2' < q_3' < \cdots
    \end{align*}
    with $\lim\limits_{i \to \infty} q_i = \lim\limits_{i \to \infty} q'_i =
    \alpha$.  For each $i$, let
    \begin{align*}
      V_{i,t} &= \{(x,y) : x > t \text{ and } q'_i x < y < q_i x \} \\
      \mathcal{V}_i &= \{V_{i,t} : t \in M\}.
    \end{align*}
    Again, $\mathcal{V}_i$ is finer than $\mathcal{V}_j$ if and only if $i \ge j$, and
    $\bigcup_{i = 1}^\infty \mathcal{V}_i$ generates $p_\alpha$.
    Therefore $p_\alpha$ does not have a definable basis, by
    Lemma~\ref{sdef}.\ref{sdef3}.
  \item Suppose $p = p_{q^-}$ for some rational $q \in \mathbb{Q}$.  Take a
    sequence of rational numbers
    \begin{equation*}
      q_1' < q_2' < q_3' < \cdots
    \end{equation*}
    with $\lim\limits_{i \to \infty} q'_i = q$.  For each $i$, let
    \begin{align*}
      V_{i,s,t} &= \{(x,y) : x > t \text{ and } q'_i x < y < q x + s \} \\
      \mathcal{V}_i &= \{V_{i,s,t} : s,t \in M\}.
    \end{align*}
    As before, $\mathcal{V}_i$ is finer than $\mathcal{V}_j$ if and only if $i \ge j$,
    and $\bigcup_{i = 1}^\infty \mathcal{V}_i$ generates $p_{q^-}$.
    So $p_{q^-}$ does not have a definable basis.
  \item The case of $p = p_{q^+}$ is handled similarly.
  \end{itemize}
  Therefore, we have ruled out all cases where $p = p_\xi$, and it
  must be that $p = r_{q,\xi}$, where $\xi$ is $m, m^+,$ or $m^-$ for
  some $m \in M$.
\end{proof}

\subsection{The counterexample}

Define
\begin{align*}
  \mathbb{R}^+ = \{x \in \mathbb{R} : x > 0\}. \\
  \mathbb{R}^- = \{x \in \mathbb{R} : x < 0\}.
\end{align*}
\begin{proposition}\label{analysis}
~
  \begin{enumerate}
  \item \label{part1} In the structure $(\mathbb{R}^+,\cdot,\le)$, let
    $p(x,y)$ be a 2-type containing $x > t$ for all $t \in \mathbb{R}^+$.
    Suppose $p(x,y)$ has a definable basis.  Then there is some $q \in
    \mathbb{Q}$ and $c \in \mathbb{R}^+$ such that $p(x,y)$ is generated by one of the following
    three families:
    \begin{align*}
      \{x > t : t \in \mathbb{R}^+\} &\cup \{y = x^q \cdot c\} \\
      \{x > t : t \in \mathbb{R}^+\} &\cup \{x^q \cdot c < y < x^q \cdot c \cdot \mu &&: \mu > 1\} \\
      \{x > t : t \in \mathbb{R}^+\} &\cup \{x^q \cdot c \cdot \mu^{-1} < y <
      x^q \cdot c &&: \mu > 1\}
    \end{align*}
  \item \label{part2} In the structure $(\mathbb{R},\cdot,\le)$, let $p(x,y)$
    be a 2-type containing $y > 0$ and $x > t$ for
    all $t \in \mathbb{R}^+$.  Suppose $p(x,y)$ has a definable basis.  Then there is some $q \in \mathbb{Q}$ and $c \in
    \mathbb{R}^+$ such that $p(x,y)$ is generated by one of the following three families:
    \begin{align*}
      \{x > t : t \in \mathbb{R}^+\} &\cup \{y = x^q \cdot c\} \\
      \{x > t : t \in \mathbb{R}^+\} &\cup \{x^q \cdot c < y < x^q \cdot c \cdot \mu &&: \mu > 1\} \\
      \{x > t : t \in \mathbb{R}^+\} &\cup \{x^q \cdot c \cdot \mu^{-1} < y <
      x^q \cdot c &&: \mu > 1\}
    \end{align*}
  \item \label{part3} In the structure $(\mathbb{R},\cdot,\le)$, let $p(x,y)$
    be a 2-type containing $y < 0$ and $x > t$ for
    all $t \in \mathbb{R}^+$.  Suppose $p(x,y)$ has a definable basis.  Then there is some $q \in \mathbb{Q}$ and $c \in
    \mathbb{R}^-$ such that $p(x,y)$ is generated by one of the following three families:
    \begin{align*}
      \{x > t : t \in \mathbb{R}^+\} &\cup \{y = x^q \cdot c\} \\
      \{x > t : t \in \mathbb{R}^+\} &\cup \{x^q \cdot c > y > x^q \cdot c \cdot \mu &&: \mu > 1\} \\
      \{x > t : t \in \mathbb{R}^+\} &\cup \{x^q \cdot c \cdot \mu^{-1} > y >
      x^q \cdot c &&: \mu > 1\}
    \end{align*}
  \end{enumerate}
\end{proposition}
\begin{proof}
  ~
  \begin{enumerate}
  \item Apply Proposition~\ref{classify2} to the densely ordered
    abelian group $(\mathbb{R}^+,\cdot,\le)$.
  \item This follows from Part~\ref{part1}.  The two structures
    $(\mathbb{R},\cdot,\le)$ and $(\mathbb{R}^+,\cdot,\le)$ are bi-interpretable,
    and we can interpret $(\mathbb{R}^+,\cdot,\le)$ as the definable set
    $\mathbb{R}^+$ in $(\mathbb{R},\cdot,\le)$.
  \item This follows by pushing forward Part~\ref{part2} along
    the definable bijection
    \begin{align*}
      \mathbb{R}^+ \times \mathbb{R}^+ & \to \mathbb{R}^+ \times \mathbb{R}^- \\
      (x,y) & \mapsto (x,-y). \qedhere
    \end{align*}
  \end{enumerate}
\end{proof}
Let $C$ be the 4-ary predicate on $\mathbb{R}$ defined by
\begin{equation*}
  C(x,y,z,w) \iff \left(x < y < 0 < z < w \text{ and } \frac{x}{y} = \frac{w}{z}\right).
\end{equation*}
\begin{lemma}
  The two structures $(\mathbb{R},\le,\cdot)$ and $(\mathbb{R},\le,C)$ have the same
  definable sets.\footnote{However, $(\mathbb{R},\le,C)$ has fewer
  $\emptyset$-definable sets than $(\mathbb{R},\le,\cdot)$.  This will be seen implicitly in the proof of Lemma~\ref{target} below.}
\end{lemma}
\begin{proof}
  We claim that the following four structures all have the same
  definable sets:
  \begin{align*}
    M_1 &= (\mathbb{R},\le,\cdot) \\
    M_2 &= (\mathbb{R},\le,\cdot,1,-1,C) \\
    M_3 &= (\mathbb{R},\le,1,-1,C) \\
    M_4 &= (\mathbb{R},\le,C).
  \end{align*}
  First of all, $M_2$ is a definitional expansion of $M_1$, so they
  have the same definable sets.  It is an easy exercise that the
  relation $x \cdot y = z$ is definable in $M_3$.  Therefore $M_2$ is a
  definitional expansion of $M_3$.  Finally, $M_3$ is an expansion of
  $M_4$ by naming parameters.
\end{proof}
Because the definable sets in $(\mathbb{R},\le,C)$ are the same as in
$(\mathbb{R},\le,\cdot)$, the following makes sense, and follows immediately
from Proposition~\ref{analysis}.\ref{part3}:
\begin{proposition}\label{analysis-c}
  In the structure $(\mathbb{R},\le,C)$, let $p(x,y)$ be a 2-type containing
  $y < 0$ and $x > t$ for all $t \in \mathbb{R}^+$.  Suppose $p(x,y)$ has a
  definable basis.  Then there is some $q \in \mathbb{Q}$ and $c \in \mathbb{R}^-$
  such that $p(x,y)$ is generated by one of the following three families:
  \begin{align*}
    \{x > t : t \in \mathbb{R}^+\} &\cup \{y = x^q \cdot c\} \\
    \{x > t : t \in \mathbb{R}^+\} &\cup \{x^q \cdot c > y > x^q \cdot c
    \cdot \mu &&: \mu > 1\} \\
    \{x > t : t \in \mathbb{R}^+\} &\cup \{x^q \cdot c \cdot \mu^{-1} > y >
    x^q \cdot c &&: \mu > 1\}
  \end{align*}  
\end{proposition}
\begin{lemma}\label{target}
  None of the three types in Proposition~\ref{analysis-c} is
  $\emptyset$-definable in the structure $(\mathbb{R},\le,C)$.
\end{lemma}
\begin{proof}
  Let $\sigma : \mathbb{R} \to \mathbb{R}$ be the function
  \begin{equation*}
    \sigma(x) = 
    \begin{cases}
      2x & x < 0 \\
      x & x \ge 0.
    \end{cases}
  \end{equation*}
  Then $\sigma$ is an automorphism of $(\mathbb{R},\le,C)$: it clearly
  preserves $\le$, and it preserves $C$ because if $x < y < 0 < z <
  w$, then
  \begin{equation*}
    \frac{x}{y} = \frac{w}{z} \iff \frac{2x}{2y} = \frac{w}{z} \iff
    \frac{\sigma(x)}{\sigma(y)} = \frac{\sigma(w)}{\sigma(z)}.
  \end{equation*}
  Suppose $p(x,y)$ is the type generated by one of the following families, for some $q \in \mathbb{Q}$
  and $c \in \mathbb{R}^-$:
  \begin{align*}
    \{x > t : t \in \mathbb{R}^+\} &\cup \{y = x^q \cdot c\} \\
    \{x > t : t \in \mathbb{R}^+\} &\cup \{x^q \cdot c > y > x^q \cdot c
    \cdot \mu &&: \mu > 1\} \\
    \{x > t : t \in \mathbb{R}^+\} &\cup \{x^q \cdot c \cdot \mu^{-1} > y >
    x^q \cdot c &&: \mu > 1\}
  \end{align*}
  Then the pushforward along $\sigma$ is, respectively, the type generated by one of these families:
  \begin{align*}
    \{x > t : t \in \mathbb{R}^+\} &\cup \{y = 2 x^q \cdot c\} \\
    \{x > t : t \in \mathbb{R}^+\} &\cup \{2 x^q \cdot c > y > 2 x^q \cdot c
    \cdot \mu &&: \mu > 1\} \\
    \{x > t : t \in \mathbb{R}^+\} &\cup \{2 x^q \cdot c \cdot \mu^{-1} > y >
    2 x^q \cdot c &&: \mu > 1\}
  \end{align*}
  
  The new type $\sigma_* p$ is never equal to $p$, because $c\neq 2c$.  Therefore $p$ is not fixed by some automorphism over
  $\emptyset$, and $p$ is not $\emptyset$-definable.
\end{proof}

\begin{theorem}
  The structure $(\mathbb{R},\le,C)$ is o-minimal.  Let $\mathcal{S}$ be the
  definable family of sets $\{(x,y) : y < 0, ~ x > t\}$ for $t \in
  \mathbb{R}$.  Then $\mathcal{S}$ is a $\emptyset$-definable downward directed
  family.  If $\mathcal{S}'$ is a definable downward directed family finer than
  $\mathcal{S}$ that generates a 2-type $p(x,y)$,
  then $\mathcal{S}'$ and $p$ are \emph{not} $\emptyset$-definable.
\end{theorem}
\begin{proof}
  The structure $(\mathbb{R},\le,C)$ is o-minimal because, for example, it is
  a reduct of RCF.  The definable family $\mathcal{S}$ is clearly
  downward directed.  Because $\mathcal{S}'$ is finer than $\mathcal{S}$, the type
  $p(x,y)$ contains the formulas $y < 0$ and $x > t$ for all $t \in
  \mathbb{R}$.  Therefore $p(x,y)$ must be one of the three types listed in
  Proposition~\ref{analysis-c}.  By Lemma~\ref{target}, $p(x,y)$ is
  not $\emptyset$-definable. If $\mathcal{S}'$ were
  $\emptyset$-definable, then $p$ would be as well.
\end{proof}

\bibliography{mybib_types_transversals_compactness}

\newcommand{\etalchar}[1]{$^{#1}$}
\providecommand{\noopsort}[1]{}\newcommand{\SortNoop}[1]{}
\begin{thebibliography}{{\SortNoop{Dries}}vdD98}

\bibitem[ADH{\etalchar{+}}16]{vc_density}
Matthias Aschenbrenner, Alf Dolich, Deirdre Haskell, Dugald Macpherson, and
  Sergei Starchenko.
\newblock Vapnik-{C}hervonenkis density in some theories without the
  independence property, {I}.
\newblock {\em Trans. Amer. Math. Soc.}, 368(8):5889--5949, 2016.

\bibitem[AF11]{aschen_fischer_11}
Matthias Aschenbrenner and Andreas Fischer.
\newblock Definable versions of theorems by {K}irszbraun and {H}elly.
\newblock {\em Proc. Lond. Math. Soc. (3)}, 102(3):468--502, 2011.

\bibitem[AG21]{andujar_thesis}
Pablo And\'ujar~Guerrero.
\newblock {\em Definable Topological Spaces in O-minimal Structures}.
\newblock PhD thesis, Purdue University, 2021.

\bibitem[AGTW21]{atw1}
Pablo And\'ujar~Guerrero, Margaret Thomas, and Erik Walsberg.
\newblock Directed sets and topological spaces definable in o-minimal
  structures.
\newblock {\em J. London Math. Soc. (2)}, 104(3):989--1010, 2021.

\bibitem[BK18]{boxall_kestner_18}
Gareth Boxall and Charlotte Kestner.
\newblock The definable {$(P,Q)$}-theorem for distal theories.
\newblock {\em J. Symb. Log.}, 83(1):123--127, 2018.

\bibitem[CK12]{cher_kap_12}
Artem Chernikov and Itay Kaplan.
\newblock Forking and dividing in {${\rm NTP}_2$} theories.
\newblock {\em J. Symbolic Logic}, 77(1):1--20, 2012.

\bibitem[CS12]{cotter_star_12}
Sarah Cotter and Sergei Starchenko.
\newblock Forking in {VC}-minimal theories.
\newblock {\em J. Symbolic Logic}, 77(4):1257--1271, 2012.

\bibitem[Dol04]{dolich04}
Alfred Dolich.
\newblock Forking and independence in o-minimal theories.
\newblock {\em J. Symbolic Logic}, 69(1):215--240, 2004.

\bibitem[{\SortNoop{Dries}}vdD98]{dries98}
Lou {\SortNoop{Dries}}~van~den Dries.
\newblock {\em Tame Topology and O-minimal Structures}, volume 248 of {\em
  London Mathematical Society Lecture Note Series}.
\newblock Cambridge University Press, Cambridge, 1998.

\bibitem[For]{fornasiero}
Antongiulio Fornasiero.
\newblock Definable compactness for topological structures.
\newblock In preparation.

\bibitem[HL16]{hruloeser16}
Ehud Hrushovski and Fran\c{c}ois Loeser.
\newblock {\em Non-archimedean tame topology and stably dominated types},
  volume 192 of {\em Annals of Mathematics Studies}.
\newblock Princeton University Press, Princeton, NJ, 2016.

\bibitem[JL10]{john_las_10}
H.~R. Johnson and M.~C. Laskowski.
\newblock Compression schemes, stable definable families, and o-minimal
  structures.
\newblock {\em Discrete Comput. Geom.}, 43(4):914--926, 2010.

\bibitem[Joh18]{johnson14}
Will Johnson.
\newblock Interpretable sets in dense o-minimal structures.
\newblock {\em J. Symb. Log.}, 83(4):1477--1500, 2018.

\bibitem[KM00]{kar_mac_00}
Marek Karpinski and Angus Macintyre.
\newblock Approximating volumes and integrals in o-minimal and {$p$}-minimal
  theories.
\newblock In {\em Connections between model theory and algebraic and analytic
  geometry}, volume~6 of {\em Quad. Mat.}, pages 149--177. Dept. Math., Seconda
  Univ. Napoli, Caserta, 2000.

\bibitem[Mat04]{matousek04}
Ji\v{r}\'{\i} Matou\v{s}ek.
\newblock Bounded {VC}-dimension implies a fractional {H}elly theorem.
\newblock {\em Discrete Comput. Geom.}, 31(2):251--255, 2004.

\bibitem[Mir71]{mirsky71}
L.~Mirsky.
\newblock A dual of {D}ilworth's decomposition theorem.
\newblock {\em Amer. Math. Monthly}, 78:876--877, 1971.

\bibitem[MS94]{mark_stein_94}
David Marker and Charles~I. Steinhorn.
\newblock Definable types in {$O$}-minimal theories.
\newblock {\em J. Symbolic Logic}, 59(1):185--198, 1994.

\bibitem[Pil87]{pillay87}
Anand Pillay.
\newblock First order topological structures and theories.
\newblock {\em J. Symbolic Logic}, 52(3):763--778, 1987.

\bibitem[Pil88]{pillay88}
Anand Pillay.
\newblock On groups and fields definable in {$o$}-minimal structures.
\newblock {\em J. Pure Appl. Algebra}, 53(3):239--255, 1988.

\bibitem[Poi00]{Poizat}
Bruno Poizat.
\newblock {\em A Course in Model Theory}.
\newblock Universitext. Springer, 2000.

\bibitem[PP07]{pet_pillay_07}
Ya'acov Peterzil and Anand Pillay.
\newblock Generic sets in definably compact groups.
\newblock {\em Fund. Math.}, 193(2):153--170, 2007.

\bibitem[PS87]{pillay_stein_87}
Anand Pillay and Charles Steinhorn.
\newblock On {D}edekind complete {$o$}-minimal structures.
\newblock {\em J. Symbolic Logic}, 52(1):156--164, 1987.

\bibitem[PS99]{pet_stein_99}
Ya'acov Peterzil and Charles Steinhorn.
\newblock Definable compactness and definable subgroups of o-minimal groups.
\newblock {\em J. London Math. Soc. (2)}, 59(3):769--786, 1999.

\bibitem[Sim15]{simon15}
Pierre Simon.
\newblock Invariant types in {NIP} theories.
\newblock {\em J. Math. Log.}, 15(2):1550006, 26, 2015.

\bibitem[SS14]{simon_star_14}
Pierre Simon and Sergei Starchenko.
\newblock On forking and definability of types in some {DP}-minimal theories.
\newblock {\em J. Symb. Log.}, 79(4):1020--1024, 2014.

\bibitem[Tho12]{thomas12}
Margaret E.~M. Thomas.
\newblock Convergence results for function spaces over o-minimal structures.
\newblock {\em J. Log. Anal.}, 4(1), 2012.

\bibitem[vdD05]{dries03}
Lou van~den Dries.
\newblock Limit sets in o-minimal structures.
\newblock In {\em O-minimal Structures, Lisbon 2003}, Proceedings of a Summer
  School by the European Research and Training Network, RAAG, pages 172--215,
  2005.

\bibitem[Wal15]{walsberg15}
Erik Walsberg.
\newblock On the topology of metric spaces definable in o-minimal expansions of
  fields.
\newblock {\em arXiv e-prints}, 2015.
\newblock arXiv:1510.07291.

\bibitem[Wal19]{walsberg19}
Erik Walsberg.
\newblock The {M}arker-{S}teinhorn theorem via definable linear orders.
\newblock {\em Notre Dame J. Form. Log.}, 60(4):701--706, 2019.

\end{thebibliography}
\bibliographystyle{alpha}
\end{document}